
\def\ra{\rightarrow}

\def\fl{\forall}
\def\ify{\infty}
\def\lgl{\langle}
\def\nb{\nabla}
\def\op{\oplus}
\def\ot{\otimes}

\def\part{\partial}
\def\rgl{\rangle}
\def\sbs{\subset}
\def\sm{\simeq}
\def\ts{\times}
\def\wdg{\wedge}

\def\wt{\widetilde}

\def\a{\alpha}
\def\b{\beta}
\def\d{\delta}
\def\g{\gamma}

\def\lb{\lambda}
\def\om{\omega}
\def\s{\sigma}
\def\t{\theta}
\def\ve{\varepsilon}
\def\vp{\varphi}

\def\D{\Delta}
\def\G{\Gamma}

\def\Om{\Omega}

\def\Diff{\mathop{\rm Diff}\nolimits}

\font\tenbb=msbm10
\font\sevenbb=msbm7
\font\fivebb=msbm5
\newfam\bbfam
\textfont\bbfam=\tenbb \scriptfont\bbfam=\sevenbb
\scriptscriptfont\bbfam=\fivebb
\def\bb{\fam\bbfam}

\def\Cb{{\bb C}}

\def\Rb{{\bb R}}

\def\Zb{{\bb Z}}

\def\Ac{{\cal A}}

\def\Hc{{\cal H}}
\def\Lc{{\cal L}}

\def\Uc{{\cal U}}

\def\build#1_#2^#3{\mathrel{
\mathop{\kern 0pt#1}\limits_{#2}^{#3}}}

\def\hfl#1#2{\smash{\mathop{\hbox to 6mm{\rightarrowfill}}
\limits^{\scriptstyle#1}_{\scriptstyle#2}}}

\def\hfll#1#2{\smash{\mathop{\hbox to 6mm{\leftarrowfill}}
\limits^{\scriptstyle#1}_{\scriptstyle#2}}}

\def\boxit#1#2{\setbox1=\hbox{\kern#1{#2}\kern#1}%
\dimen1=\ht1 \advance\dimen1 by #1 \dimen2=\dp1 \advance\dimen2 by #1
\setbox1=\hbox{\vrule height\dimen1 depth\dimen2\box1\vrule}%
\setbox1=\vbox{\hrule\box1\hrule}%
\advance\dimen1 by .4pt \ht1=\dimen1
\advance\dimen2 by .4pt \dp1=\dimen2 \box1\relax}

\def\semi{\mathop{>\!\!\!\triangleleft}}

\catcode`\@=11
\def\displaylinesno #1{\displ@y\halign{
\hbox to\displaywidth{$\@lign\hfil\displaystyle##\hfil$}&
\llap{$##$}\crcr#1\crcr}}

\def\ldisplaylinesno #1{\displ@y\halign{
\hbox to\displaywidth{$\@lign\hfil\displaystyle##\hfil$}&
\kern-\displaywidth\rlap{$##$}
\tabskip\displaywidth\crcr#1\crcr}}
\catcode`\@=12

\baselineskip=14pt

\hsize=118mm 
\hoffset=20mm
\vsize=215mm
\voffset=15mm

\def\Ac{\cal A}
\def\Bc{\cal B}
\def\Hc{\cal H}
\def\Lc{\cal L}
\def\Uc{\cal U}

\font\tenbb=msbm10
\font\sevenbb=msbm7
\font\fivebb=msbm5
\newfam\bbfam
\textfont\bbfam=\tenbb \scriptfont\bbfam=\sevenbb
\scriptscriptfont\bbfam=\fivebb
\def\bb{\fam\bbfam}

\def\Cb{{\bb C}}
\def\Nb{{\bb N}}
\def\Rb{{\bb R}}
\def\Zb{{\bb Z}}

\def\a{\alpha}
\def\b{\beta}
\def\d{\delta}
\def\D{\Delta}
\def\g{\gamma}
\def\G{\Gamma}
\def\lb{\lambda}
\def\L{\Lambda}
\def\om{\omega}
\def\Om{\Omega}
\def\s{\sigma}
\def\t{\theta}
\def\ve{\varepsilon}
\def\vp{\varphi}

\def\bsh{\backslash}
\def\ify{\infty}
\def\nb{\nabla}
\def\part{\partial}
\def\sm{\simeq}
\def\ts{\times}
\def\wdg{\wedge}

\def\lgl{\langle}
\def\rgl{\rangle}

\def\fl{\forall}
\def\op{\oplus}
\def\ot{\otimes}
\def\sbs{\subset}

\def\wt{\widetilde}

\def\ra{\rightarrow}
\def\hra{\hookrightarrow}

\def\build#1_#2^#3{\mathrel{
\mathop{\kern 0pt#1}\limits_{#2}^{#3}}}

\def\semi{\mathop{>\!\!\!\triangleleft}}

\def\xx{\vrule height 0.5em depth 0.2em width 0.5em}


\centerline{\bf Hopf algebras, cyclic cohomology and}

\centerline{\bf the transverse index theorem}

\vglue 1cm

\centerline{A. Connes and H. Moscovici}

\vglue 1cm

\noindent {\bf Introduction}
\bigskip

In this paper we present the solution of
a longstanding internal problem of 
noncommutative geometry, namely the computation of 
the index of transversally 
elliptic operators on foliations.

\smallskip

\noindent The spaces of leaves of foliations are basic examples of
noncommutative spaces and already exhibit most of the 
features of the general theory. 
The index problem for longitudinal elliptic operators
is easy to formulate in the presence of a
transverse measure, cf.~[Co] 
[M-S], and in general it leads 
to the construction (cf.~[C-S]) of a natural map from the 
geometric group to the $K$-theory of the leaf space,
i.e. the $K$-theory of the associated $C^*$-algebra. 
This ``assembly map'' $\mu$ is known in many cases to exhaust the 
$K$-theory  of 
the $C^*$-algebra but property $T$ in the group context and its analogue for 
foliations provide conceptual obstructions to 
tentative proofs of its surjectivity in 
general. One way to test the $K$-group, $K (C^* (V,F)) = K (V/F)$ for short, 
is to use its natural pairing with the $K$-{\it homology} group of $C^* (V,F)$. Cycles in the latter represent 
``abstract elliptic operators'' on $V/F$ and the 
explicit construction for {\it general} foliations of such cycles
is already quite an elaborate problem. 
The delicate point is that we do not want to assume any special 
property of the foliation such as, for instance, the existence of a 
holonomy invariant transverse metric as in Riemannian foliations. 
Equivalently, we do not 
want to restrict in anyway the holonomy pseudogroup of the foliation.

\smallskip

\noindent In [Co1] [H-S] [C-M] a general solution was given to 
the construction of 
transversal elliptic operators for foliations. 
The first step ([Co1]) consists in 
passing by a Thom isomorphism to the total space of the bundle of transversal 
metrics. This first step is a geometric adaptation of the reduction of an 
arbitrary factor of type III to a crossed product of a factor of type II by a 
one-parameter 
group of automorphisms. Instead of only taking care of the volume 
distorsion (as 
in the factor case) of the involved elements of the pseudogroup, 
it takes care 
of their full Jacobian. The second step ([H-S]) consisted in realizing that 
while the standard theory of {\it elliptic} pseudodifferential operators is
too restrictive to allow the construction of 
the desired $K$-homology cycle, 
it suffices to replace it by its refinement to {\it hypoelliptic} operators. 
This was used in [C-M] in order to 
construct a {\it differential} (hypoelliptic) operator $Q$, 
solving the general construction of the $K$-cycle.

\smallskip

\noindent One then arrives at a well posed general index problem. The index 
defines a map: $K(V/F) \ra \Zb$ which is 
simple to compute for those elements 
of 
$K(V/F)$ in the range of the assembly map. The problem is to provide a general 
formula for the cyclic cocycle ${\rm ch}_* (D)$, which computes the index by 
the 
equality
$$
\lgl {\rm ch}_* (D) , {\rm ch}^* (E) \rgl = {\rm Index} \, D_E \qquad \fl \, E 
\in K(V/F) , \leqno (1)
$$
where the Chern character ${\rm ch}^* (E)$ belongs to the cyclic homology of 
$V/F$. We showed in [C-M] that the spectral triple given by the algebra 
${\Ac}$ 
of the foliation, together with the operator $D$ in Hilbert space ${\Hc}$ 
actually fulfills the hypothesis of a general abstract index theorem, holding 
at 
the operator theoretic level. It gives a ``local'' formula for the cyclic 
cocycle ${\rm ch}_* (D)$ in terms of certain residues 
that extend the ideas of the 
Wodzicki-Guillemin-Manin residue as well as of the Dixmier trace. 
Adopting the notation 
${\displaystyle \int \!\!\!\!\!\! -}$ for such a residue, the general formula 
gives the components $\vp_n$ of the cyclic cocycle $\vp = {\rm ch}_* (D)$ as 
universal {\it finite} linear combinations of expressions which have the 
following general form
$$
\int \!\!\!\!\!\! - \, a^0 \, [D,a^1]^{(k_1)} \ldots [D,a^n]^{(k_n)} \, \vert 
D 
\vert^{-n-2 \vert k \vert} \, , \qquad \fl \, a^j \in {\Ac} , \leqno (2)
$$
where for an operator $T$ in $\Hc$ the symbol $T^{(k)}$ means the $k^{\rm th}$ 
iterated commutator of $D^2$ with $T$.

\smallskip

\noindent It was soon realized that, although the general index 
formula easily 
reduces to the local form of the Atiyah-Singer index theorem 
when $D$ is say a Dirac 
operator on a manifold, the actual explicit computation of all the terms (2) 
involved in the cocycle ${\rm ch}_* (D)$ is a rather formidable task. As an 
instance of this let us mention that even in the case of codimension one 
foliations, the printed form of the explicit computation of the cocycle takes 
around one hundred pages. Each step in the computation is straightforward but 
the explicit computation for higher values of $n$ is clearly impossible 
without 
a new organizing principle which allows to bypass them.

\bigskip

\noindent In this paper we shall adapt and develop the theory of cyclic cohomology to 
Hopf algebras and show that this provides exactly the missing organizing principle, thus allowing to perform the computation for arbitrary values of $n$. We shall 
construct for each value of $n$ a specific Hopf algebra ${\Hc} (n)$, 
show that 
it acts on the $C^*$-algebra of the transverse frame bundle 
of any codimension 
$n$ foliation $(V,F)$ and that the index computation takes place within the 
cyclic cohomology of ${\Hc} (n)$. 
We compute this cyclic cohomology explicitly 
as Gelfand-Fuchs cohomology. While the link between cyclic cohomology and 
Gelfand-Fuchs cohomology was already known ([Co]), the novelty
consists in the fact that the entire differentiable 
transverse structure is now
captured by the action of the Hopf algebra ${\Hc}(n)$, 
thus reconciling our approach to noncommutative geometry to a more group 
theoretical one, in the spirit of the Klein program.

\vglue 1cm

\noindent {\bf I. Notations}

\smallskip

We let $M$ be an $n$-dimensional smooth manifold (not necessarily connected or 
compact but assumed to be oriented). Let us first fix the notations for the 
frame bundle of $M$, $F(M)$, in local coordinates
$$
x^{\mu} \quad \mu = 1 , \ldots , n \quad \hbox{for} \quad x \in U \sbs M \, . 
\leqno (1)
$$
We view a frame with coordinates $x^{\mu} , y_j^{\mu}$ as the 1-jet of the map
$$
j : \Rb^n \ra M \, , \ j(t) = x + yt \qquad \fl \, t \in \Rb^n \leqno (2)
$$
where $(y t)^{\mu} = y_i^{\mu} \, t^i \qquad \fl \, t = (t^i) \in \Rb^n$.

\smallskip

\noindent Let $\vp$ be a local (orientation preserving) diffeomorphism of $M$, 
it acts on $F(M)$ by
$$
\vp , j \ra \vp \circ j = \wt{\vp} \, (j) \leqno (3)
$$
which replaces $x$ by $\vp (x)$ and $y$ by $\vp' (x) \, y$ where
$$
\vp' (x)^{\a}_{\b} = \part_{\b} \, \vp (x)^{\a} \quad \hbox{where} \quad \vp 
(x) 
= (\vp (x)^{\a}) \, . \leqno (4)
$$
We restrict our attention to orientation preserving frames $F^+ (M)$, and in 
the 
one dimensional case $(n=1)$ we take the notation
$$
y = e^{-s} \, , \ s \in \Rb \, . \leqno (5)
$$
In terms of the coordinates $x,s$ one has,
$$
\wt{\vp} (s,x) = (s-\log \vp' (x) , \vp (x)) \leqno (6)
$$
and the invariant measure on $F$ is $(n=1)$
$$
{dx \, dy \over y^2} = e^s \, ds \, dx \, . \leqno (7)
$$
One has a canonical right action of $GL^+ (n,\Rb)$ on $F^+$ which is given by
$$
(g,j) \ra j \circ g \, , \ g \in GL^+ (n,\Rb) \, , \ j \in F^+ \leqno (8)
$$
it replaces $y$ by $yg$, $(yg)_j^{\mu} = y_i^{\mu} \, g_j^i \quad \fl \, g \in 
GL^+ (n,\Rb)$ and $F^+$ is a $GL^+ (n,\Rb)$ principal bundle over $M$.

\smallskip

\noindent We let $Y_i^j$ be the vector fields on $F^+$ generating the action 
of $GL^+ (n,\Rb)$,
$$
Y_i^j = y_i^{\mu} \, {\part \over \part \, y_j^{\mu}} = y_i^{\mu} \, 
\part_{\mu}^j \, . \leqno (9)
$$
In the one dimensional case one gets a single vector field,
$$
Y = -\part_s \, . \leqno (10)
$$
The action of ${\rm Diff}^+$ on $F^+$ preserves the $\Rb^n$ valued 1-form on 
$F^+$,
$$
\a^j = (y^{-1})_{\b}^j \, dx^{\b} \, . \leqno (11)
$$
One has $y_j^{\mu} \, \a^j = dx^{\mu}$ and $dy_j^{\mu} \wdg \a^j + y_j^{\mu} \, 
d\a^j =0$.

\smallskip

\noindent Given an affine torsion free connection $\G$, the associated one form 
$\om$,
$$
\om_j^{\ell} = (y^{-1})_{\mu}^{\ell} \, (dy_j^{\mu} + \G_{\a , \b}^{\mu} \, 
y_j^{\a} \, dx^{\b}) \leqno (12)
$$
is a 1-form on $F^+$ with values in {\bf GL(n)} the Lie algebra of $GL^+ 
(n,\Rb)$. The $\G_{\a , \b}^{\mu}$ only depend on $x$ but not on $y$, moreover 
one has,
$$
d\a^j = \a^k \wdg \om_k^j = -\om_k^j \wdg \a^k \leqno (13)
$$
since $\G$ is {\it torsion free}, i.e. $\G_{\a , \b}^{\mu} = \G_{\b , 
\a}^{\mu}$.

\smallskip

\noindent The natural horizontal vector fields $X_i$ on $F^+$ associated to the 
connection $\G$ are,
$$
X_i = y_i^{\mu} (\part_{\mu} - \G_{\a , \mu}^{\b} \, y_j^{\a} \, \part_{\b}^j) 
\, , \leqno (14)
$$
they are characterized by
$$
\lgl \a^j , X_i \rgl = \d_i^j \quad \hbox{and} \quad \lgl \om_k^{\ell} , X_i 
\rgl = 0 \, . \leqno (15)
$$
For $\psi \in {\rm Diff}^+$, the one form $\wt{\psi}^* \, \om$ is still a 
connection 1-form for a new affine torsion free connection $\G'$. The new 
horizontal vector fields $X'_i$ are related to the old ones by
$$
X'_i = \wt{\vp}_* \, X_i \circ \wt{\psi} \, , \ \vp = \psi^{-1} \, . \leqno 
(16)
$$
When $\om$ is the trivial flat connection $\G = 0$ one gets
$$
\G' = \psi' (x)^{-1} \, d\psi' (x) \, , \ \G_{\a , \b}^{\mu} = (\psi' 
(x)^{-1})_{\rho}^{\mu} \, \part_{\b} \, \part_{\a} \, \psi^{\rho} (x) \, . 
\leqno (17)
$$

\vglue 1cm

\noindent {\bf II. Crossed product of $F(M)$ by $\G$ and action of ${\Hc} (n)$}

\smallskip

We let $M$ be gifted with a {\it flat} affine connection $\nb$ and let $\G$ be 
a 
pseudogroup of local diffeomorphisms, preserving the orientation,
$$
\psi : {\rm Dom} \, \psi \ra {\rm Range} \, \psi \leqno (1)
$$
where both the domain, ${\rm Dom} \, \psi$ and range, ${\rm Range} \, \psi$ are 
open sets of $M$. By the functoriality of the construction of $F(M)^+$,
 $\forall \psi \in \G$ we let 
$\wt{\psi}$ be the corresponding local diffeomorphism of $F^+ (M)$.

\smallskip

\noindent We let ${\Ac} = C_c^{\ify} (F^+) \semi \G$ be the crossed product of 
$F^+$ by the action of $\G$ on $F^+$. It can be described directly as 
$C_c^{\ify} (G)$ where $G$ is the etale smooth groupoid,
$$
G = F^+ \semi \G \, , \leqno (2)
$$
an element $\g$ of $G$ being given by a pair $(x,\vp)$, $x \in {\rm Range} \, 
\vp$, while the composition is,
$$
(x,\vp) \circ (y,\psi) = (x,\vp \circ \psi) \quad \hbox{if} \quad y \in {\rm 
Dom} \, \vp \quad \hbox{and} \quad \vp (y) = x \, . \leqno (3)
$$
In practice we shall generate the crossed product ${\Ac}$ as the linear span of 
monomials,
$$
f \, U_{\psi}^* \, , \ f \in C_c^{\ify} ({\rm Dom} \, \psi) \leqno (4)
$$
where the star indicates a contravariant notation. The multiplication rule is
$$
f_1 \, U_{\psi_1}^* \ f_2 \, U_{\psi_2}^* = f_1 (f_2 \circ \wt{\psi}_1) \, 
U_{\psi_2 \psi_1}^* \leqno (5)
$$
where by hypothesis the support of $f_1 (f_2 \circ \wt{\psi}_1)$ is a compact 
subset of
$$
{\rm Dom} \, \psi_1 \cap \psi_1^{-1} \, {\rm Dom} \, \psi_2 \sbs {\rm Dom} \, 
\psi_2 \, \psi_1 \, . \leqno (6)
$$
The canonical action of $GL^+ (n,\Rb)$ on $F(M)^+$ commutes with the action of 
$\G$ and thus extends canonically to the crossed product $\Ac$. At the Lie 
algebra level, this yields the following derivations of $\Ac$,
$$
Y_{\ell}^j (f \, U_{\psi}^*) = (Y_j^{\ell} \, f) \, U_{\psi}^* \, . \leqno (7)
$$
Now the flat connection $\nb$ also provides us with associated horizontal 
vector 
fields $X_i$ on $F^+ (M)$ (cf. section I) which we extend to the crossed 
product 
$\Ac$ by the rule,
$$
X_i (f \, U_{\psi}^*) = X_i (f) \, U_{\psi}^* \, . \leqno (8)
$$
Now, of course, unless the $\psi$'s are {\it affine}, the $X_i$ do not commute 
with the action of $\psi$, but using (16) and (17) of section I we can compute 
the corresponding commutator and get,
$$
X_i - U_{\psi} \, X_i \, U_{\psi}^* = -\g_{ij}^k \, Y_k^j \leqno (9)
$$
where the functions $\g_{ij}^k$ are,
$$
\g_{ij}^k = y_i^{\mu} \, y_j^{\a} \, (y^{-1})_{\b}^k \, \G_{\a \mu}^{\b} \, 
, \leqno (10)
$$
$$ 
  \ \G_{\a , \b}^{\mu} = (\psi' 
(x)^{-1})_{\rho}^{\mu} \, \part_{\b} \, \part_{\a} \, \psi^{\rho} (x) \, . 
$$

It follows that, for any $a,b \in {\Ac}$ one has
$$
X_i (ab) = X_i (a) \, b + a \, X_i (b) + \d_{ij}^k (a) \, Y_k^j (b) \leqno (11)
$$
where the linear operators $\d_{ij}^k$ in $\Ac$ are defined by,
$$
\d_{ij}^k (f \, U_{\psi}^*) = \g_{ij}^k \, f \, U_{\psi}^* \, . \leqno (12)
$$
To prove (11) one takes $a=f_1 \, U_{\psi_1}^*$, $b = f_2 \, U_{\psi_2}^*$ and 
one computes $X_i (ab)$ $=$ $X_i (f_1 \, U_{\psi_1}^* \ f_2 \, U_{\psi_2}^*)$ 
$=$ $X_i (f_1 \, U_{\psi_1}^* \ f_2 \, U_{\psi_1}) \, U_{\psi_2 \psi_1}^*$ $=$ 
$X_i (f_1) \, U_{\psi_1}^* \, f_2 \, U_{\psi_2}^* + f_1 (X_i \, U_{\psi_1}^*$ 
$-$ $U_{\psi_1}^* \, X_i) \, f_2 \, U_{\psi_2}^* + f_1 \, U_{\psi_1}^* \, X_i 
(f_2 \, U_{\psi_2}^*)$.

\smallskip

\noindent One then uses (9) to get the result.

\smallskip

\noindent Next the $\g_{ij}^k$ are characterized by the equality
$$
\wt{\psi}^* \, \om - \om = \g_{jk}^{\ell} \, \a^k = \g \, \a \leqno (13)
$$
where $\a$ is the canonical $\Rb^n$-valued one form on $F^+ (M)$ (cf.~I).

\smallskip

\noindent The equality $\wt{\psi_2 \psi_1}^* \, \om - \om = \wt{\psi_1}^* 
(\wt{\psi_2}^* \, \om - \om) + (\wt{\psi_1}^* \, \om - \om)$ together with the 
invariance of $\a$ thus show that the $\g_{ij}^k$ form a 1-cocycle, so that 
each 
$\d_{ij}^k$ is a derivation of the algebra $\Ac$,
$$
\d_{ij}^k (ab) = \d_{ij}^k (a) \, b + a \, \d_{ij}^k (b) \, . \leqno (14)
$$
Since the connection $\nb$ is flat the commutation relations between the 
$Y_j^{\ell}$ and the $X_i$ are those of the affine group,
$$
\Rb^n \semi GL^+ (n,\Rb) \, . \leqno (15)
$$
The commutation of the $Y_j^{\ell}$ with $\d_{ab}^c$ are easy to compute since 
they correspond to the tensorial nature of the $\d_{ab}^c$. The $X_i$ however 
do 
not have simple commutation relations with the $\d_{ab}^c$, and one lets
$$
\d_{ab,i_1 \ldots i_n}^c = [X_{i_1} , \ldots [X_{i_n} , \d_{ab}^c] \ldots ] \, 
. 
\leqno (16)
$$
All these operators acting on $\Ac$ are of the form,
$$
T \, (f \, U_{\psi}^*) = h \, f \, U_{\psi}^* \leqno (17)
$$
where $h = h^{\psi}$ is a function depending on $\psi$.

\smallskip

\noindent In particular they all commute pairwise,
$$
[\d_{ab,i_1 \ldots i_n}^c , \d_{a'b',i'_1 , \ldots , i'_m}^{c'} ] = 0 \, . 
\leqno (18)
$$
It follows that the linear space generated by the $Y_j^{\ell}$, $X_i$, 
$\d_{ab,i_1 \ldots i_n}^c$ forms a Lie algebra and we let $\Hc$ be the 
corresponding envelopping algebra. We endow $\Hc$ with a coproduct in such a 
way 
that its action on $\Ac$,
$$
h,a \ra h(a) \, , \ h \in {\Hc} \, , \ a \in {\Ac} \leqno (19)
$$
satisfies the following rule,
$$
h(ab) = \sum \, h_{(0)} \, (a) \, h_{(1)} \, (b) \qquad \fl \, a,b \in {\Ac} 
\quad \hbox{where} \quad \D h = \sum \, h_{(0)} \ot h_{(1)} \, . \leqno (20)
$$
One gets from the above discussion the equalities
$$
\D \, Y_i^j = Y_i^j \ot 1 + 1 \ot Y_i^j \leqno (21)
$$
$$
\D X_i = X_i \ot 1 + 1 \ot X_i + \d_{ij}^k \ot Y_k^j \leqno (22)
$$
$$
\D \d_{ij}^k = \d_{ij}^k \ot 1 + 1 \ot \d_{ij}^k \, . \leqno (23)
$$
These rules, together with the equality
$$
\D (h_1 \, h_2) = \D h_1 \, \D h_2 \qquad \fl \, h_j \in {\Hc} \leqno (24)
$$
suffice to determine completely the coproduct in $\Hc$. As we shall see $\Hc$ 
has 
an antipode $S$, we thus get a Hopf algebra ${\Hc} (n)$ which only depends 
upon 
the integer $n$ and which acts on any crossed product,
$$
{\Ac} = C_c^{\ify} (F) \semi \G \leqno (25)
$$
of the frame bundle of a {\it flat} manifold $M$ by a pseudogroup $\G$ of 
local diffeomorphisms.

\smallskip

\noindent We shall devote a large portion of this paper
to the understanding of the structure of the 
Hopf algebra ${\Hc} (n)$ as well as of its cyclic cohomology. 
For notational simplicity we shall concentrate 
on the case $n=1$ 
but all the results are proved in 
such a way as to extend in a straightforward manner to the general case.

\smallskip

\noindent To end this section we shall show that provided we replace $\Ac$ by a 
Morita equivalent algebra we can bypass the flatness condition of the manifold 
$M$.

\smallskip

\noindent To do this we start with an arbitrary manifold $M$ (oriented) and we 
consider a locally finite open cover $(U_{\a})$ of $M$ by domains of local 
coordinates. On $N = \coprod \, U_{\a}$, the disjoint union of the open sets 
$U_{\a}$, one has a natural pseudogroup $\G_0$ of diffeomorphisms which satisfy
$$
\pi \, \psi (x) = \pi (x) \qquad \fl \, x \in {\rm Dom} \, \psi \leqno (26)
$$
where $\pi : N \ra M$ is the natural projection.

\smallskip

\noindent Equivalently one can consider the smooth etale groupoid which is the 
graph of the equivalence relation $\pi (x) = \pi (y)$ in $N$,
$$
G_0 = \{ (x,y) \in N \ts N \, ; \ \pi (x) = \pi (y) \} \, . \leqno (27)
$$
One has a natural Morita equivalence,
$$
C_c^{\ify} (M) \sm C_c^{\ify} (G_0) = C_c^{\ify} (N) \semi \G_0 \leqno (28)
$$
which can be concretely realized as the reduction of $C_c^{\ify} (G_0)$ by the 
idempotent,
$$
e \in C_c^{\ify} (G_0) = C_c^{\ify} (N) \semi \G_0 \, , \ e^2 = e \, , \leqno 
(29)
$$
associated to a partition of unity in $M$ subordinate to the cover $(U_{\a})$,
$$
\sum \, \vp_{\a} (x)^2 = 1 \, , \ \vp_{\a} \in C_c^{\ify} (U_{\a}) \leqno (30)
$$
by the formula,
$$
e (u,\a ,\b) = \vp_{\a} (u) \, \vp_{\b} (u) \, . \leqno (31)
$$
We have labelled the pair $(x,y) \in G_0$ by $u = \pi (x) = \pi (y)$ and the 
indices $\a ,\b$ so that $x \in U_{\a}$, $y \in U_{\b}$.

\smallskip

\noindent This construction also works in the presence of a 
pseudogroup $\G$ of local
diffeomorphisms of $M$ since there is a 
corresponding pseudogroup $\G'$ on $N$ 
containing $\G_0$ and such that, with the above projection $e$,
$$
(C_c^{\ify} (N) \semi \G')_e \sm C_c^{\ify} (M) \semi \G \, . \leqno (32)
$$
Now the manifold $N$ is obviously flat and the above 
construction of the action 
of the Hopf algebra ${\Hc} (n)$ gives an action on ${\Ac}' = C_c^{\ify} (N) 
\semi \G'$, $({\Ac}')_e = {\Ac} = C_c^{\ify} (M) \semi \G$.

\vglue 1cm

\noindent {\bf III. One dimensional case, the Hopf algebras ${\Hc}_n$}

\smallskip

We first define a bialgebra by generators and relations. As an algebra we 
view $\Hc$ as the envelopping algebra of the Lie algebra which is the 
linear span of $Y$, $X$, $\d_n$, $n \geq 1$ with the relations,
$$
[Y,X] = X , [Y,\d_n] = n \, \d_n , [\d_n , \d_m] = 0 \quad \fl \, n , m 
\geq 1 , [X,\d_n] = \d_{n+1} \quad \fl \, n \geq 1 \, . \leqno (1)
$$
We define the coproduct $\D$ by
$$
\D \, Y = Y \ot 1 + 1 \ot Y \ , \ \D \, X = X \ot 1 + 1 \ot X + \d_1 \ot Y 
\ , \ \D \, \d_1 = \d_1 \ot 1 + 1 \ot \d_1 \leqno (2)
$$
and with $\D \, \d_n$ defined by induction using (1).
\smallskip

\noindent One checks that the presentation (1) is preserved by $\D$, so 
that $\D$ extends to an algebra homomorphism,
$$
\D : {\Hc} \ra {\Hc} \ot {\Hc} \leqno (3)
$$
and one also checks the coassociativity.

\smallskip

For each $n$ we let ${\Hc}_n$ be the algebra generated by $\d_1 , \ldots 
, \d_n$,
$$
{\Hc}_n = \{ P (\d_1 , \ldots , \d_n) \ ; \ P \ \hbox{polynomial in} \ n 
\ \hbox{variables} \} \, . \leqno (4)
$$
We let ${\Hc}_{n,0}$ be the ideal,
$$
{\Hc}_{n,0} = \{ P ; P(0) = 0 \} \, . \leqno (5)
$$
By induction on $n$ one proves the following

\medskip

\noindent {\bf Lemma 1.} {\it For each $n$ there exists $R_{n-1} \in 
{\Hc}_{n-1,0} \ot {\Hc}_{n-1,0}$ such that $\D \, \d_n = \d_n \ot 1 + 1 
\ot \d_n + R_{n-1}$.}

\medskip

\noindent {\it Proof.} It holds for $n=1$, $n=2$. 
Assuming that it holds for $n$ one has $\D \, \d_{n+1} = \D \, [X,\d_n] 
= [\D \, X , \D \, \d_n] = [X \ot 1 + 1 \ot X + \d_1 \ot Y , \d_n \ot 1 
+ 1 \ot \d_n + R_{n-1}] = \d_{n+1} \ot 1 + 1 \ot \d_{n+1} + [ X \ot 1 + 
1 \ot X , R_{n-1}] + \d_1 \ot [Y,\d_n] + [\d_1 \ot Y , R_{n-1}] = 
\d_{n+1} \ot 1 + 1 \ot \d_{n+1} + R_n$ where
$$
R_n = [X \ot 1 + 1 \ot X , R_{n-1}] + n \, \d_1 \ot \d_n + [\d_1 \ot Y , 
R_{n-1}] \, . \leqno (6)
$$
Since $[X, {\Hc}_{n-1,0}] \sbs  {\Hc}_{n,0}$ and $[Y, {\Hc}_{n-1,0}] 
\sbs  {\Hc}_{n-1,0} \sbs {\Hc}_{n,0}$, one gets that $R_n \in 
{\Hc}_{n,0}$.~\xx

\medskip

\noindent For each $k \leq n$ we introduce a linear form $Z_{k,n}$ on 
${\Hc}_n$
$$
\lgl Z_{k,n} , P \rgl = \left( {\part \over \part \, \d_k} \, P \right) 
(0) \, . \leqno (7)
$$
One has by construction,
$$
\lgl Z_{k,n} , PQ \rgl = \lgl Z_{k,n} , P \rgl \, Q (0) + P(0) \, \lgl 
Z_{k,n} , Q \rgl \leqno (8)
$$
and moreover $\ve$, $\lgl \ve , P \rgl = P(0)$ is the counit in 
${\Hc}_n$,
$$
\lgl L \ot \ve , \D \, P \rgl = \lgl \ve \ot L , \D \, P \rgl = \lgl L,P 
\rgl \qquad \fl \, P \in {\Hc}_n \, . \leqno (9)
$$
(Check both sides on a monomial $P = \d_1^{a_1} \ldots \d_n^{a_n}$.)

\smallskip

\noindent Thus in the dual agebra ${\Hc}_n^*$ one can write (8) as
$$
\D \, Z_{k,n} = Z_{k,n} \ot 1 + 1 \ot Z_{k,n} \, . \leqno (10)
$$
Moreover the $Z_{k,n}$ form a basis of the linear space of solutions of 
(10) and we need to determine the Lie algebra structure determined by the 
bracket.

\smallskip

\noindent We let for a better normalization,
$$
Z'_{k,n} = (k+1) ! \, Z_{k,n} \, . \leqno (11)
$$

\medskip

\noindent {\bf Lemma 2.} {\it One has $[Z'_{k,n} , Z'_{\ell,n}] = (\ell 
- k) \, Z'_{k+\ell , n}$ if $k + \ell \leq n$ and $0$ if $k + \ell > 
n$.}

\medskip

\noindent {\it Proof.} Let $P = \d_1^{a_1} \ldots \d_n^{a_n}$ be a 
monomial. We need to compute $\lgl \D \, P , Z_{k,n} \ot Z_{\ell , n} - 
Z_{\ell , n} \ot Z_{k,n} \rgl$. One has
$$
\D \, P = (\d_1 \ot 1 + 1 \ot \d_1)^{a_1} \, (\d_2 \ot 1 + 1 \ot \d_2 + 
R_1)^{a_2} \ldots (\d_n \ot 1 + 1 \ot \d_n + R_{n-1})^{a_n} \, .
$$
We look for the terms in $\d_k \ot \d_{\ell}$ or $\d_{\ell} \ot \d_k$ 
and take the difference. The latter is non zero only if all $a_j = 0$ 
except $a_q = 1$. Moreover since $R_m$ is homogeneous of degree $m+1$ 
one gets $q=k+\ell$ and in particular $[Z'_{k,n} , Z'_{\ell ,n}] = 0$ if 
$k + \ell > n$. One then computes by induction using (6) the bilinear 
part of $R_m$. One has $R_1^{(1)} = \d_1 \ot \d_1$, and from (6)
$$
R_n^{(1)} = [(X \ot 1 + 1 \ot X) , R_{n-1}^{(1)}] + n \, \d_1 \ot \d_n 
\, . \leqno (12)
$$
This gives
$$
R_{n-1}^{(1)} = \d_{n-1} \ot \d_1 + C_n^1 \, \d_{n-2} \ot \d_2 + \ldots 
+ C_n^{n-2} \, \d_1 \ot \d_{n-1} \, . \leqno (13)
$$
Thus the coefficient of $\d_k \ot \d_{\ell}$ is $C_{k+\ell}^{\ell - 1}$ 
and we get
$$
[Z_{k,n} , Z_{\ell , n}] = (C_{k+\ell}^{\ell - 1} - C_{k+\ell}^{k - 1}) 
\, Z_{k+\ell , n} \, . \leqno (14)
$$
One has ${(k+1) ! \, (\ell + 1) ! \over (k+\ell + 1) !} \, 
(C_{k+\ell}^{\ell - 1} - C_{k+\ell}^{k - 1}) = {\ell (\ell + 1) - k 
(k+1) \over k + \ell + 1} = \ell - k$ thus using (11) one gets the 
result.~\xx

\medskip

\noindent For each $n$ we let ${\Ac}_n^1$ be the Lie algebra of vector 
fields
$$
f(x) \, \part / \part x \quad , \quad f(0) = f' (0) = 0 \leqno (15)
$$
modulo $x^{n+2} \, \part$.

\smallskip

\noindent The elements $Z_{k,n} = {x^{k+1} \over (k+1) !} \, \part / 
\part x$ are related by (11) to $Z'_{k,n} = x^{k+1} \, \part / \part x$ 
which satisfy the Lie algebra of lemma 2.

\smallskip

\noindent Thus ${\Ac}_n^1$ is the Lie algebra of jets of order $(n+1)$ 
of vector fields which vanish to order 2 at 0.

\medskip

\noindent {\bf Proposition 3.} {\it The Hopf algebra ${\Hc}_n$ is the 
dual of the envelopping agebra ${\Uc} ({\Ac}_n^1)$, ${\Hc}_n = {\Uc} 
({\Ac}_n^1)^*$.}

\medskip

\noindent {\it Proof.} This follows from the Milnor-Moore theorem.~\xx

\medskip

\noindent Since the ${\Ac}_n^1$ form a projective system of Lie 
algebras, with limit the Lie algebra ${\Ac}^1$ of formal vector fields 
which vanish to order 2 at 0, the inductive limit ${\Hc}^1$ of the Hopf 
algebras ${\Hc}_n$ is,
$$
{\Hc}^1 = {\Uc} ({\Ac}^1)^* \, . \leqno (16)
$$
The Lie algebra ${\Ac}^1$ is a {\it graded} Lie algebra, with one 
parameter group of automorphisms,
$$
\a_t \, (Z_n) = e^{nt} \, Z_n \leqno (17)
$$
which extends to ${\Uc} ({\Ac}^1)$ and transposes to ${\Uc} ({\Ac}^1)^*$ 
as
$$
\lgl [Y,P] , a \rgl = \left\lgl P , {\part \over \part \, t} \, \a_t \, 
(a)_{t=0} \right\rgl \qquad \fl \, P \in {\Hc}^1 \ , \ a \in {\Uc} 
({\Ac}^1) \, . \leqno (18)
$$
Indeed $(\a^*_t)$ is a one parameter group of automorphisms of ${\Hc}^1$ 
such that
$$
\a^*_t \, (\d_n) = e^{nt} \, \d_n \, . \leqno (19)
$$
One checks directly that $\a^*_t$ is compatible with the coproduct on 
${\Hc}^1$ and that the corresponding Lie algebra automorphism is (17).

\smallskip

\noindent Now (cf. [Dix] 2.1.11) we take the basis of ${\Uc} ({\Ac}^1)$ 
given by the monomials,
$$
Z_n^{a_n} \, Z_{n-1}^{a_{n-1}} \ldots Z_2^{a_2} \, Z_1^{a_1} \ , \ a_j 
\geq 0 \, . \leqno (20)
$$
To each $L \in {\Uc} ({\Ac}^1)^*$ one associates (cf. [Dix] 2.7.5) the 
formal power series
$$
\sum \, {L (Z_n^{a_n} \ldots Z_1^{a_1}) \over a_n ! \ldots a_1 !} \, 
x_1^{a_1} \ldots x_n^{a_n} \, , \leqno (21)
$$
in the commuting variables $x_j$, $j \in \Nb$.

\smallskip

\noindent It follows from [Dix] 2.7.5 that we obtain in this way an 
isomorphism of the algebra of polynomials $P (\d_1 , \ldots , \d_n)$ on 
the algebra of polynomials in the $x_j$'s. To determine the formula for 
$\d_n$ in terms of the $x_j$'s, we just need to compute
$$
\lgl \d_n , Z_n^{a_n} \ldots Z_1^{a_1} \rgl \, . \leqno (22)
$$
Note that (22) vanishes unless $\sum \, j \, a_j = n$.

\smallskip

\noindent In particular, for $n=1$, we get
$$
\rho \, (\d_1) = x_1 \leqno (23)
$$
where $\rho$ is the above isomorphism.

\smallskip

\noindent We determine $\rho \, (\d_n)$ by induction, using the 
derivation
$$
D(P) = \sum \, \d_{n+1} \, {\part \over \part \, \d_n} \, (P) \leqno 
(24)
$$
(which corresponds to $P \ra [X,P]$).

\smallskip

\noindent One has by construction,
$$
\lgl \d_n , a \rgl = \lgl \d_{n-1} , D^t (a) \rgl \qquad \fl \, a \in 
{\Uc} ({\Ac}^1) \leqno (25)
$$
where $D^t$ is the transpose of $D$.

\smallskip

\noindent By definition of $Z_n$ as a linear form (7) one has,
$$
D^t \, Z_n = Z_{n-1} \ , \ n \geq 2 \ , \ D^t \, Z_1 = 0 \, . \leqno 
(26)
$$
Moreover the compatibility of $D^t$ with the coproduct of ${\Hc}^1$ is
$$
D^t (ab) = D^t (a) \, b + a \, D^t (b) + (\d_1 \, a) \, \part_t \, b 
\qquad \fl \, a,b \in {\Uc} ({\Ac}^1) \leqno (27)
$$
where $a \ra \d_1 \, a$ is the natural action of the algebra ${\Hc}^1$ on 
its dual
$$
\lgl P , \d_1 \, a \rgl = \lgl P \, \d_1 , a \rgl \qquad \fl \, P \in 
{\Hc}^1 \, , \ a \in {\Uc} ({\Ac}^1) \, . \leqno (28)
$$
To prove (27) one pairs both side with $P \in {\Hc}^1$. The $\ell . h . 
s .$ gives $\lgl P , D^t (ab) \rgl = \lgl \D \, [X,P] , a \ot b \rgl = 
\lgl [X \ot 1 + 1 \ot X + \d_1 \ot Y , \D \, P ] , a \ot b \rgl$. The 
terms in $[X \ot 1 , \D \, P ]$ yield $\lgl \D \, P , D^t \, a \ot b 
\rgl$ and similarly for $[1 \ot X , \D \, P ]$. The term $[\d_1 \ot Y , 
\D \, P ]$ yield $\lgl \D \, P , \d_1 \, a \ot \part_t \, b \rgl$ thus 
one gets (27).

\medskip

\noindent {\bf Lemma 4.} {\it When restricted to ${\Uc} ({\Ac}^2)$, 
$D^t$ is the unique derivation, with values in ${\Uc} ({\Ac}^1)$ 
satisfying (26), moreover
$$
D^t (Z_n^{a_n} \ldots Z_2^{a_2} \, Z_1^{a_1}) = D^t (Z_n^{a_n} \ldots 
Z_2^{a_2}) \, Z_1^{a_1} + Z_n^{a_n} \ldots Z_2^{a_2} \, {a_1 (a_1 - 1) 
\over 2} \, Z_1^{a_1 - 1} \, .
$$
}

\smallskip

\noindent {\it Proof.} The equality $\D \, \d_1 = \d_1 \ot 1 + 1 \ot 
\d_1$ shows that $a \ra \d_1 \, a$ is a derivation of ${\Uc} ({\Ac}^1)$. 
One has $\d_1 \, Z_n = 0$ for $n \ne 1$ so that $\d_1 = 0$ on ${\Uc} 
({\Ac}^2)$ and the first statement follows from (27) and (26). The 
second statement follows from,
$$
D^t (Z_1^m) = {m (m-1) \over 2} \, Z_1^{m-1} \leqno (29)
$$
which one proves by induction on $m$ using (27).~\xx

\medskip

\noindent Motivated by the first part of the lemma, we enlarge the Lie algebra 
${\Ac}^1$ by adjoining an element $Z_{-1}$ such that,
$$
[Z_{-1} , Z_n] = Z_{n-1} \qquad \fl \, n \geq 2 \, , \leqno (30)
$$
we then define $Z_0$ by
$$
[Z_{-1} , Z_1] = Z_0 \ , \ [Z_0 , Z_k] = k \, Z_k \, . \leqno (31)
$$
The obtained Lie algebra $\Ac$, is the Lie algebra of formal vector 
fields with $Z_0 = x \, {\part \over \part \, x}$, $Z_{-1} = {\part \over 
\part \, x}$ and as above $Z_n = {x^{n+1} \over (n+1)!} \, {\part \over 
\part \, x}$.

\smallskip

\noindent We can now compare $D^t$ with the bracket with $Z_{-1}$. They 
agree on ${\Uc} ({\Ac}^2)$ and we need to compute $[Z_{-1} , Z_1^m ]$. 
One has
$$
[Z_{-1} , Z_1^m ] = {m (m-1) \over 2} \, Z_1^{m-1} + m \, Z_1^{m-1} 
\, Z_0 \leqno (32)
$$
($[Z_{-1} , Z_1^m \, Z_1] = \left({m (m-1) \over 2} \, Z_1^{m-1} + m \, 
Z_1^{m-1} \, Z_0 \right) \, Z_1 + Z_1^m \, Z_0 =$ $\left({m (m-1) \over 
2} + m \right)$ $Z_1^m + (m+1) \, Z_1^m \, Z_0$).

\smallskip

\noindent Thus, if one lets $\Lc$ be the left ideal in ${\Uc} ({\Ac})$ 
generated by $Z_{-1}$, $Z_0$ we get,

\medskip

\noindent {\bf Proposition 5.} {\it The linear map $D^t : {\Uc} ({\Ac}^1) 
\ra {\Uc} ({\Ac}^1)$ is uniquely determined by the equality $D^t (a) = 
[Z_{-1} , a]$ {\rm mod} $\Lc$.}

\medskip

\noindent {\it Proof.} For each monomial $Z_n^{a_n} \ldots Z_1^{a_1}$ one 
has $D^t (a) - [Z_{-1} , a] \in {\Lc}$, so that this holds for any $a \in 
{\Uc} ({\Ac}^1) $. Moreover, using the basis of ${\Uc} ({\Ac})$ given by 
the $Z_n^{a_n} \ldots Z_1^{a_1} \, Z_0^{a_0} \, Z_{-1}^{a_{-1}}$ we see 
that ${\Uc} ({\Ac})$ is the direct sum ${\Lc} \op {\Uc} ({\Ac}^1)$.~\xx

\medskip

\noindent (The linear span of the $Z_n^{a_n} \ldots  Z_0^{a_0} \, 
Z_{-1}^{a_{-1}}$ with $a_0 + a_{-1} > 0$ is a left ideal in ${\Uc} 
({\Ac})$ since the product $(Z_m^{b_m} \ldots Z_{-1}^{b_{-1}})$ 
$(Z_n^{a_n} \ldots  Z_0^{a_0} \, Z_{-1}^{a_{-1}})$ can be expressed by 
decomposing $(Z_m^{b_m} \ldots Z_{-1}^{b_{-1}} \, Z_n^{a_n} \ldots 
Z_1^{a_1})$ as a sum of monomials $Z_q^{c_q} \ldots Z_1^{c_1} \, 
Z_0^{c_0}$ $Z_{-1}^{c_{-1}}$ which are then multiplied by $Z_0^{a_0} \, 
Z_{-1}^{a_{-1}}$ which belongs to the augmentation ideal of $\Uc$ (Lie 
algebra of $Z_0$, $Z_1$).)

\smallskip

\noindent We now define a linear form $L_0$ on ${\Uc} ({\Ac})$ by
$$
L_0 (Z_n^{a_n} \ldots Z_1^{a_1} \, Z_0^{a_0} \, Z_{-1}^{a_{-1}}) = 0 \ 
\hbox{unless} \ a_0 = 1 \, , \ a_j = 0 \quad \fl \, j \, , \leqno (33)
$$
and $L_0 (Z_0) = 1$.

\medskip

\noindent {\bf Proposition 6.} {\it For any $n \geq 1$ one has}
$$
\lgl \d_n , a \rgl = L_0 ([ \underbrace{\ldots}_{n \, {\rm times}} 
[Z_{-1} , a ] \ldots ]) \qquad \fl \, a \in {\Uc} ({\Ac}^1) \, .
$$

\smallskip

\noindent {\it Proof.} Let us first check it for $n=1$. We let $a = 
Z_n^{a_n} \ldots Z_1^{a_1}$. Then the degree of $a$ is $\sum \, j \, a_j$ 
and $L_0 ([Z_{-1} , a ]) \ne 0$ requires $\sum \, j \, a_j = 1$ so that 
the only possibility is $a_1 = 1$, $a_j = 0 \quad \fl \, j$. In this case 
one gets $L_0 ([Z_{-1} , Z_1]) = L_0 (Z_0) = 1$. Thus by (23) we get the 
equality of Proposition 6 for $n=1$.

\smallskip

\noindent For the general case note first that $\Lc$ is stable under 
right multiplication by $Z_{-1}$ and hence by the derivation $[Z_{-1} , 
\cdot]$. Thus one has
$$
(D^t)^n \, (a) = [Z_{-1} , \ldots [Z_{-1} , a] \ldots ] \ \hbox{mod} \ 
{\Lc} \qquad \fl \, a \in {\Uc} ({\Ac}^1) \, . \leqno (34)
$$
Now for $a \in {\Lc}$ one has $L_0 ([Z_{-1} , a ]) = 0$. Indeed writing 
$a = ( Z_n^{a_n} \ldots Z_1^{a_1})$ $(Z_0^{a_0} \, Z_{-1}^{a_{-1}}) = bc$ 
with $b \in {\Uc} ({\Ac}^1)$, $c = Z_0^{a_0} \, Z_{-1}^{a_{-1}}$, one has 
$[Z_{-1} , a] = [Z_{-1} , b] \, c + b \, [Z_{-1} , c]$. Since $b \in 
{\Uc} ({\Ac}^1)$ and $[Z_{-1} , c]$ has strictly negative degree one has 
$L_0 (b \, [Z_{-1} , c]) = 0$. Let $Z_n^{b_n} \ldots Z_1^{b_1} \, 
Z_0^{b_0}$ be a non zero component of $[Z_{-1} , b]$, then unless all 
$b_i$ are 0 it contributes by 0 to $L_0 ([Z_{-1} , b] \, c)$. But 
$[Z_{-1} , b] \in {\Uc} ({\Ac}^0)_0$ has no constant term. Thus one has
$$
L_0 ([Z_{-1} , a]) = 0 \qquad \fl \, a = Z_n^{a_n} \ldots Z_1^{a_1} \, 
Z_0^{a_0} \, Z_{-1}^{a_{-1}} \leqno (35)
$$
except if all $a_j = 0$, $j \ne 1$ and $a_1 = 1$. $L_0 ([Z_{-1} , Z_1]) = 
1$.

\smallskip

\noindent Using (25) one has $\lgl \d_n , a \rgl = \lgl \d_1 , 
(D^t)^{n-1} \, (a) \rgl$ and the lemma follows.~\xx

\medskip

\noindent One can now easily compute the first values of $\rho \, 
(\d_n)$, $\rho \, (\d_1) = x_1$, $\rho \, (\d_2) = x_2 + {x_1^2 \over 
2}$, $\rho \, (\d_3) = x_3 + x_2 \, x_1 + { x_1^3 \over 2}$, $\rho \, 
(\d_4) = x_4 + x_3 \, x_1 + {2} \, x_2^2 + {2} \, x_2 \, 
x_1^2 + {3 \over 4} \, x_1^4$.

\smallskip

\noindent The affine structure provided by the $\d_n$ has the following 
compatibility with left multiplication in ${\Uc} ({\Ac}^1)$,

\medskip

\noindent {\bf Proposition 7.} a) {\it One has $R_{n-1} = \sum \, 
R_{n-1}^k \ot \d_k$, $R_{n-1}^k \in {\Hc}_{n-1,0}$.} 

\noindent b) {\it For fixed 
$a_0 \in {\Uc} ({\Ac}^1)$ there are $\lb_n^k \in \Cb$ such that}
$$
\lgl \d_n , (a_0 \, a) \rgl = \lgl \d_n , a_0 \rgl \, \ve (a) + \sum 
\lb_n^k \, \lgl \d_k , a \rgl \, .
$$

\smallskip

\noindent {\it Proof.} a) By induction using (6). b) Follows, using 
$\lb_n^k = \lgl R_{n-1}^k , a_0 \rgl$.~\xx

\medskip

\noindent The antipode $S$ in ${\Uc} ({\Ac}^1)$ is the unique 
antiautomorphism such that
$$
S \, Z_n = - Z_n \qquad \fl \, n \, . \leqno (36)
$$
It is non trivial to express in terms of the coordinates $\d_n$.

\smallskip

\noindent In fact if we use the basis $Z_j$ of ${\Ac}^1$ but in reverse 
order to construct the map $\rho$ we obtain a map $\wt{\rho}$ whose first 
values are $\wt{\rho} \, (\d_1) = z_1$, $\wt{\rho} \, (\d_2) = z_2 + 
{z_1^2 \over 2}$, $\wt{\rho} \, (\d_3) = z_3 + 3 \, z_1 \, z_2 + {1 \over 
2} \, z_1^3$, $\wt{\rho} \, (\d_4) = z_4 + {2} \, z_2^2 + 6 \, 
z_1 \, z_3 + {9} \, z_1^2 \, z_2 + {3 \over 4} \, z_1^4$.

\smallskip

\noindent One has $\lgl \d_n , S \, (Z_m^{a_m} \ldots Z_1^{a_1}) \rgl = 
(-1)^{\sum a_j} \, \lgl \d_n , Z_1^{a_1} \ldots Z_m^{a_m} \rgl$ so 
that $\rho \, (S^t \, \d_n)$ $= \sum \lgl \d_n , S \, (Z_m^{a_m} \ldots 
Z_1^{a_1}) \rgl \, x_1^{a_1} \ldots x_m^{a_m} = \sum (-1)^{\sum a_j} 
\, \lgl \d_n , Z_1^{a_1} \ldots Z_m^{a_m} \rgl \, x_1^{a_1} \ldots$ 
$x_m^{a_m} = \wt{\rho} \, (\d_n)$ with $z_j = -x_j$ in the latter 
expression. 

\smallskip

\noindent Thus $\rho \, (S^t \, \d_1) = -x_1$, $\rho \, (S^t \, \d_2) = 
-x_2 + {x_1^2 \over 2}$, $\rho \, (S^t \, \d_3) = -x_3 + 3 \, x_1 \, x_2 
- { x_1^3 \over 2}$, $\rho \, (S^t \, \d_4) = -x_4 + {2} \, 
x_2^2 + 6 \, x_1 \, x_3 - {9} \, x_1^2 \, x_2 + {3\over 4} \, 
x_1^4$. We thus get
$$
S^t \, \d_1 = -\d_1 \ , \ S^t \, \d_2 = -\d_2 + \d_1^2 \ , \ S^t \, \d_3 
= -\d_3 + 4 \, \d_1 \, \d_2 - 2 \, \d_1^3 \ , \ \ldots \leqno (37)
$$
The antipode $S$ is characterized abstractly as the inverse of the 
element $L(a) = a$ in the algebra of linear maps $L$ from ${\Uc} 
({\Ac}^1)$ to ${\Uc} ({\Ac}^1)$ with the product
$$
(L_1 * L_2) (a) = \sum \, L_1 (a_{(1)}) \, L_2 (a_{(2)}) \qquad \D \, a = 
\sum \, a_{(1)} \ot a_{(2)} \ , \ a \in {\Uc} \, . \leqno (38)
$$
Thus one has
$$
\sum \, (S^t \, \d_{n,(1)}) \, \d_{n,(2)} = 0 \qquad \fl \, n \ , \ \D \, 
\d_n = \sum \, \d_{n,(1)} \ot \d_{n,(2)} \leqno (39)
$$
writing $S^t \, \d_n = - \d_n + P_n$ where $P_n (\d_1 , \ldots , 
\d_{n-1})$ is homogeneous of degree $n$, this allows to compute $S^t \, 
\d_n$ by induction on $n$.

\smallskip
\smallskip
\noindent {\it Remark}. Note that the Schwartzian
expression $\s = \d_2 - {1 \over 2} \, 
\d_1^2$ is uniquely characterized by
$$
\rho \, (\s) = x_2 \, ;\leqno (40)
$$
thus, $\rho^{-1} (x_n)$ can be regarded as higher analogues of the Schwartzian. 

\smallskip
\smallskip
Let us now describe in a conceptual manner the action of the Hopf 
algebra ${\Hc}_{\ify}$ on the crossed product,
$$
{\Ac} = C_c^{\ify} (F) \semi \G \leqno (41)
$$
of the frame bundle of a 1-manifold $M$ by the 
pseudogroup $\G$ associated to ${\rm Diff}^+  \, M$. 
We are given a flat connection $\nb$ on $M$, which we view as a 
$GL(1)$-equivariant section,
$$
\g : F \ra J(M) \leqno (42)
$$
from $F$ to the space of jets $ J(M) = 
\{ \Rb^n \build \ra_{}^{j} M \}$. For $\a \in F$, 
$\g \, (\a)$ is the jet
$$
\g \, (\a) = (\exp \, \nb) \circ \a \leqno (43)
$$
where $\exp \, \nb$ is the exponential map associated to the connection 
$\nb$. We let $G$ be the groupoid $F \semi \G$ and let $(\vp , \a) \in G$ 
with
$$
s \, (\vp , \a) = \a \in F \ , \ r \, (\vp , \a) = \wt{\vp} \, \a \in F \ 
, \ \vp \in \G \, . \leqno (44)
$$
We let $G ({\Ac}^1) \sbs {\Uc} ({\Ac}^1)$ (completed $I$-adically, $I =$ 
augmentation ideal) be the group like elements,
$$
\psi \ , \ \D \, \psi = \psi \ot \psi \, . \leqno (45)
$$
We then have a canonical homomorphism $\g$ from $G$ to $G ({\Ac}^1)$ given 
by
$$
\g \, (\vp , \a) = \g \, (\wt{\vp} \, \a)^{-1} \circ \vp \circ \g \, (\a) 
\leqno (46)
$$
where we identify $G ({\Ac}^1)$ with the group of germs of diffeomorphisms 
by the equality
$$
A \, f = f \circ \vp^{-1} \qquad \fl \, A \in G ({\Ac}^1) \ , \ f \ 
\hbox{function on} \ \Rb \, .
$$

\smallskip

\noindent {\bf Theorem 8.} {\it For any $f \in C_c^{\ify} (G)$ and $n \in 
\Nb$ one has,}
$$
(\d_n \, f) (g) = \d_n \, (\g \, (g)^{-1}) \, f(g) \qquad \fl \, g \in G 
\, .
$$

\smallskip

\noindent {\it Proof.} We first define a representation $\pi$ of ${\Ac}^1$ 
in the Lie algebra of vector fields on $F(\Rb)$ preserving the 
differential form $e^s \, ds \, dx$, $y = e^{-s}$,
$$
\pi \, (Z_n) = -{x^n \over n!} \, \part_s + {x^{n+1} \over (n+1)!} \, 
\part_x \, . \leqno (47)
$$
(One has $i_{Z_n} (e^s \, ds \, dx) = -{x^n \over n!} \, e^s \, dx - 
{x^{n+1} \over (n+1)!} \, e^s \, ds$ which is closed.) Let then $H$ be the 
function of $F(\Rb)$ given by
$$
H (s,x) = s \, . \leqno (48)
$$
By construction the representation $\pi$ is in fact representing $\Ac$, 
and moreover for any $a \in {\Uc} ({\Ac})$ one has,
$$
L_0 (a) = - (\pi (a) \, H) (0) \, . \leqno (49)
$$
Indeed, for $a = Z_0$ the r.h.s. is 1 and given a monomial $Z_n^{a_n} 
\ldots Z_1^{a_1} \, Z_0^{a_0} \, Z_{-1}^{a_{-1}}$, it vanishes if $a_{-1} 
> 0$ or if $a_0 > 1$ and if $a_{-1} = 0$, $a_0 = 0$. If $a_{-1} = 0$, $a_0 
= 1$ the only case in which it does not vanish is $a_j = 0 \quad \fl \, j 
> 0$.

\smallskip

\noindent One has $\pi \, (Z_{-1}) = \part_x$ and it follows from 
Proposition 6 that,
$$
\lgl \d_n , a \rgl = - (\part_x^n \, \pi (a) \, H) (0) \, . \leqno (50)
$$
Now if $a = A \in G ({\Ac}^1)$ we have, with $\psi = \vp^{-1}$, that
$$
(\pi (a) \, f) (s,x) = f \, (s - \log \psi' (x) , \psi (x)) \qquad \fl \, 
f \ \hbox{function on} \ F(\Rb) \leqno (51)
$$
and we thus have,
$$
\lgl \d_n , A \rgl = (\part_x^n \log \psi' (x))_{x=0} \, . \leqno (52)
$$
We now consider $F(M)$ with the same notations in local coordinates, i.e. 
$(y,x)$ with $y = e^{-s}$. In crossed product terms we have,
$$
\d_n (f \, U_{\vp}) = f \, \g_n \, U_{\vp} \ , \ \g_n (y,x) = y^n \, 
\part_x^n (\log \psi' (x)) \ , \ \psi = \vp^{-1} \, . \leqno (53)
$$
Now $U_{\vp}$, as a function on $G$ is the characteristic function of the 
set $\{ (\vp , \a) ; \a \in F \}$ and one has $\g (\vp , \a)$, for $\a = 
(y,x)$, given by
$$
t \ra (\vp \, (x+y \, t) - \vp \, (x)) / y \, \vp' (x) = \g \, (\vp , \a) (t) 
\, . \leqno (54)
$$

\vglue 1cm

\noindent {\bf  IV.  The dual algebra ${\Hc}^*$ }

\smallskip
\noindent To understand the dual algebra ${\Hc}^*$ we associate to $L \in 
{\Hc}^*$, viewed as a linear form on $\Hc$, assumed to be continuous in 
the $I$-adic topology, the function with values in ${\Uc} ({\Ac}^1)$,
$$
f(s,t) \ , \ \lgl f(s,t) , P \rgl = \lgl L , P \, e^{tX} \, e^{sY} \rgl 
\qquad \fl \, P \in {\Hc}^1 \, . \leqno (1)
$$
We shall now write the product in ${\Hc}^*$ in terms of the functions 
$f(s,t)$. We first recall the expansional formula,
$$
e^{A+B} = \sum_{n=0}^{\ify} \ \int_{\sum u_j = 1 , u_j \geq 0} e^{u_0 A} 
\, B^{u_1 A} \, B \ldots e^{u_n A} \, \Pi \, du_j \, . \leqno (2)
$$
We use this formula to compute $\D \, e^{tX}$, say with $t > 0$,
$$
\D \, e^{tX} = \sum_{n=0}^{\ify} \ \int_{0 \leq s_1 \leq \ldots \leq s_n 
\leq t} \Pi \, ds_i \, \d_1 (s_1) \ldots \d_1 (s_n) \, e^{tX} \ot Y (s_1) 
\ldots Y (s_n) \, e^{tX} \leqno (3)
$$
where $\d_1 (s) = e^{sX} \, \d_1 \, e^{-sX}$, $Y (s) = e^{sX} \, Y \, 
e^{-sX} = Y - sX$. One has,
$$
(Y - s_1 X) \, e^{tX} \, e^{sY} = (\part_s + (t - s_1) \, \part_t) \, 
e^{tX} \, e^{sY} \, . \leqno (4)
$$
(Since $e^{tX} \, Y \, e^{sY} = (Y - tX) \, e^{tX} \, e^{sY}$.)

\smallskip

\noindent We thus get the following formula for the product, $(t > 0)$,
$$
\matrix{
(f_1 \, f_2) (s,t) = {\displaystyle \sum_{n=0}^{\ify} \ \int_{0 \leq s_1 
\leq \ldots \leq s_n \leq t}} \Pi \, ds_i \, \d_1 (s_1) \ldots && \cr
\cr
\d_n (s_n) \, f_1 (s,t) \ (\part_s + (t - s_n) \, \part_t) \ldots (\part_s 
+ (t - s_1) \, \part_t) \, f_2 (s,t) \, . && \cr
} \leqno (5)
$$
We apply this by taking for $f_1$ the constant function,
$$
f_1 (s,t) = \vp \in G ({\Ac}^1) \sbs {\Uc} ({\Ac}^1) \leqno (6)
$$
while we take the function $f_2$ to be scalar valued.

\smallskip

\noindent One has $\d_1 (s) = e^{sX} \, \d_1 \, e^{-s X} = {\displaystyle 
\sum_{n=0}^{\ify}} \ \d_{n+1} \, {s^n \over n!}$, and its left action on 
${\Uc} ({\Ac}^1)$ is given, on group like elements $\vp$ by
$$
\d_1 (s) \, \vp = \lgl \d_1 (s) , \vp \rgl \, \vp \, . \leqno (7)
$$
(Using $\lgl \d_1 (s) \, \vp , P \rgl = \lgl \vp , P \, \d_1 (s) \rgl = 
\lgl \D \, \vp , P \ot \d_1 (s) \rgl = \lgl \vp , P \rgl \, \lgl \vp , 
\d_1 (s) \rgl$.) Moreover by (50) section III, one has
$$
\lgl \d_1 (s), \vp \rgl = - \sum \ {s^n \over n!} \, \part_x^{n+1} (\pi 
(\vp) \, H)_0 \leqno (8)
$$
while, with $\psi = \vp^{-1}$, one has $(\pi (\vp) \, H) (s,x) = s - \log 
\psi' (x)$, so that (8) gives
$$
\lgl \d_1 (s) , \vp \rgl = \sum \ {s^n \over n!} \, \part_x^n \left( 
{\psi'' \over \psi'} \right) (x)_{x=0} = \left( {\psi'' \over \psi'} 
\right) (s) , \ \psi = \vp^{-1} \, . \leqno (9)
$$
Thus we can rewrite (5) as $(t>0)$
$$
\matrix{
(\vp \, f) (s,t) = {\displaystyle \sum_{n=0}^{\ify} \ \int_{0 \leq s_1 
\leq \ldots \leq s_n \leq t}} \Pi \, ds_i \, {\displaystyle \prod_1^n \left( 
{\psi'' \over \psi'} \right)} (s_i) && \cr 
\cr
(\part_s + (t - s_n) \, \part_t) \ldots (\part_s + (t - s_1) \, \part_t) 
\, f(s,t) \, . && \cr
} \leqno (10)
$$
We first apply this formula to $f(s,t) = f(s)$, independent of $t$, we get
$$
\sum_{n=0}^{\ify} \, {1 \over n!} \left( \int_0^t \left( {\psi'' \over 
\psi'} \right) (s) \, ds \right)^n \, \part_s^n \, f(s) = f (s + \log 
\psi' (t)) \, .
$$
We then apply it to $f (s,t) = t$. The term $(\part_s + (t - s_n) \, 
\part_t) \ldots (\part_s + (t - s_1) \, \part_t) \, f$ gives $(t - s_n)$, 
thus we get,
$$
\matrix{
{\displaystyle \sum_{n=0}^{\ify} \ \int_{0 \leq s_n \leq t} {\psi'' (s_n) 
\over \psi' (s_n)} \, {1 \over (n-1)!} \left( \int_0^{s_n} \left( {\psi'' 
\over \psi'} \right) (u) \, du \right)^{n-1} (t-s_n) \, ds_n + t} && \cr
\cr
= t + {\displaystyle \int_0^t {\psi'' (s) \over \psi' (s)}} \ (\exp \log 
\psi' (s)) (t-s) \, ds && \cr
\cr
= t + {\displaystyle \int_0^t} \psi'' (s) (t-s) \, ds = t + \psi (t) - t 
\build{\psi' (0)}_{\Vert \atop 1}^{} - \build{\psi (0)}_{\Vert \atop 0}^{} 
\, . && \cr
}
$$
Thus in general we get,
$$
(\vp \, f) (s,t) = f (s + \log \psi' (t) , \psi (t)) \ , \ \psi = \vp^{-1} 
\, . \leqno (11)
$$

\vglue 1cm

\noindent {\bf V. Hopf algebra ${\Hc} (G)$ associated to a 
matched pair of subgroups}

\smallskip
In this section we recall a basic construction of Hopf algebras 
([K],[B-S],[M]). 
We let $G$ be a finite group, $G_1$, $G_2$ be subgroups of $G$ such that,
$$
G = G_1 \, G_2 \, ,  G_1 \cap G_2 = 1 \, ,  \leqno (1)
$$
i.e. we assume that any $g \in G$ admits a unique decomposition as
$$
g = k \, a \quad , \quad k \in G_1 \, , \ a \in G_2 \, . \leqno (2)
$$
Since $G_1 \cong G / G_2$ one has a natural left action of $G$ on $G_1$ which 
for $g \in G_1$ coincides with the left action of $G_1$ on itself, it is 
given by
$$
g (k) = \pi_1 (g k) \qquad \fl \, g \in G \, , \ k \in G_1 \, , \leqno (3)
$$
where $\pi_j : G \ra G_j$ are the two projections.

\smallskip

\noindent For $g \in G_1$ one has $g (k) = g k$ while for $a \in G_2$ one 
has,
$$
a (1) = 1 \qquad \fl \, a \in G_2 \, . \leqno (4)
$$
(Since $\pi_1 \, (a) = 1 \quad \fl \, a \in G_2$.)

\smallskip

\noindent Since $G_2 \cong G_1\setminus G$, one has a right action of $G$ on $G_2$ 
which restricted to $G_2 \sbs G$ is the right action of $G_2$ on itself,
$$
a \cdot g = \pi_2 \, (a g) \qquad \fl \, a \in G_2 \, , \ g \in G \, . \leqno 
(5)
$$
As above one has
$$
1 \cdot k = 1 \qquad \fl \, k \in G_1 \, . \leqno (6)
$$

\smallskip

\noindent {\bf Lemma 1.} a) {\it For $a \in G_2$, $k_1 , k_2 \in G_1$ one 
has $a (k_1 \, k_2) = a (k_1) ((a \cdot k_1) (k_2))$.}

\smallskip

\noindent b) {\it For $k \in G_1$, $a_1 , a_2 \in G_2$ one has $(a_1 \, 
a_2) \cdot k = (a_1 \cdot a_2 (k)) (a_2 \cdot k)$.}

\medskip

\noindent {\it Proof.} a) One has $a \, k_1 = k'_1 \, a'$ with $k'_1 = a 
(k_1)$, $a' = a \cdot k_1$. Then $(a \, k_1) \, k_2 = k'_1 \, a' \, k_2 = 
k'_1 \, k'_2 \, a''$ with $k'_2 = a' (k_2)$. Thus $k'_1 \, k'_2 = a (k_1 
\, k_2)$ which is the required equality.

\smallskip

\noindent b) One has $a_2 \, k = k' \, a'_2$ with $k' = a_2 (k)$, $a'_2 = 
a_2 \cdot k$. Then $a_1 \, a_2 \, k = a_1 (k' \, a'_2) = (a_1 \, k') \, 
a'_2 = k'' \, a'_1 \, a'_2$ where $a'_1 = a_1 \cdot k'$, thus $(a_1 \, 
a_2) \cdot k = a'_1 \, a'_2$ as required.~\xx

\medskip

\noindent One defines a Hopf algebra ${\Hc}$ as follows. As an algebra 
${\Hc}$ is the crossed product of the algebra of functions $h$ on $G_2$ by 
the action of $G_1$. Thus elements of ${\Hc}$ are of the form $\sum \, h_k 
\, X_k$ with the rule,
$$
X_k \, h \, X_k^{-1} = k (h) \, , \ k(h) \, (a) = h (a \cdot k) \qquad \fl 
\, a \in G_2 , k \in G_1 \, . \leqno (7)
$$
The coproduct $\D$ is defined as follows,
$$
\D \, \ve_c = \sum_{ba=c} \ve_a \ot \ve_b \ , \ \ve_c (g) = 1 \ \hbox{if} \ 
g = c \ \hbox{and} \ 0 \ \hbox{otherwise}. \leqno (8)
$$
$$
\D \, X_k = \sum_{k'} \ h_{k'}^k \, X_k \ot X_{k'} \ , \ h_{k'}^k (a) = 1 \ 
\hbox{if} \ k' = a(k) \ \hbox{and} \ 0 \ \hbox{otherwise}. \leqno (9)
$$
One first checks that $\D$ defines a covariant representation. The equality 
(8) defines a representation of the algebra of functions on $G_2$. Let us 
check that (9) defines a representation of $G_1$. First, for $k=1$ one gets 
by (4) that $\D \, X_1 = X_1 \ot X_1$. One has $\D \, X_{k_1} \ \D \, 
X_{k_2} = {\displaystyle \sum_{k'_1 , k'_2}} \, h_{k'_1}^{k_1} \, X_{k_1} \ 
h_{k'_2}^{k_2} \, X_{k_2} \ot X_{k'_1 k'_2} = {\displaystyle \sum_{k'_1 , 
k'_2}} \, h_{k'_1}^{k_1} \, k_1 (h_{k'_2}^{k_2}) \, X_{k_1 k_2} \ot X_{k'_1 
k'_2}$.

\smallskip

\noindent For $a \in G_2$ one has $(h_{k'_1}^{k_1} \, k_1 (h_{k'_2}^{k_2})) 
(a) \ne 0$ only if $k'_1 = a(k_1)$, $k'_2 = (a \cdot k_1) (k_2)$. Thus 
given $a$ there is only one term in the sum to contribute, and by  lemma 
1 , one then has $k'_1 \, k'_2 = a (k_1 \, k_2)$, thus,
$$
\D \, X_{k_1} \ \D \, X_{k_2} = \sum \, h_{k''}^{k_1 k_2} \, X_{k_1 k_2} 
\ot X_{k''} = \D \, X_{k_1 k_2} \, . \leqno (10)
$$
Next, one has $\D \, X_k \, \D \, \ve_c = {\displaystyle \sum_{k'} \, 
\sum_{ba=c}} \ h_{k'}^k \, X_k \, \ve_a \ot X_{k'} \, \ve_b$, and $X_k \, 
\ve_a = \ve_{a \cdot k^{-1}} \, X_k$, so that
$$
\D \, X_k \, \D \, \ve_c = \sum_{k'} \, \sum_{ba=c} \ h_{k'}^k \, \ve_{a 
\cdot k^{-1}} \, X_k \ot \ve_{b \cdot k'^{-1}} \, X_{k'} \, . \leqno (11)
$$
One has $\D \, \ve_{c \cdot k^{-1}} \, \D \, X_k = {\displaystyle \sum_{k'} \, 
\sum_{b' a' = c \cdot k^{-1}}} \, h_{k'}^k \, \ve_{a'} \, X_k \ot \ve_{b'} \, 
X_{k'}$. 
In (11), given $a,b$ with $ba=c$ the only $k'$ that appears is $k' = (a 
\cdot k^{-1}) (k)$. But $(a \cdot k^{-1}) (k) = a (k^{-1})^{-1}$ and $b 
\cdot k'^{-1} = b \cdot a (k^{-1})$ so that by lemma 1, $(b \cdot a 
(k^{-1})) (a \cdot k^{-1}) = (ba) \cdot k^{-1} = c \cdot k^{-1}$. Thus one 
gets
$$
\D \, X_k \, \D \, \ve_c = \D \, \ve_{c \cdot k^{-1}} \, \D \, X_k \leqno 
(12)
$$
which shows that $\D$ defines an algebra homomorphism. 

\smallskip

\noindent To show that the 
coproduct $\D$ is coassociative let us identify the dual algebra ${\Hc}^*$ 
with the crossed product,
$$
(G_1)_{\rm space} \semi G_2 \, . \leqno (13)
$$
For $f \, U_a^* \in {\Hc}^*$ we define the pairing with ${\Hc}$ by
$$
\lgl h \, X_k , f \, U_a^* \rgl = f(k) \, h(a) \leqno (14)
$$
while the crossed product rules are
$$
\matrix{
U_a^* \, f = f^a \, U_a^* \, , \ f^a (k) = f(a(k)) \qquad \fl \, k \in G_1 
\cr
\cr
U_{ab}^* = U_b^* \, U_a^* \qquad \fl \, a,b \in G_2 \, . \cr
} \leqno (15)
$$
What we need to check is
$$
\lgl \D \, h \, X_k , f \, U_a^* \ot g \, U_b^* \rgl = \lgl h \, X_k , f \, 
U_a^* \, g \, U_b^* \rgl \, . \leqno (16)
$$
We can assume that $h = \ve_c$ so that $\D \, h \, X_k = {\displaystyle 
\sum_{ba=c}} \, \ve_a \, X_k \ot \ve_b \, X_{a(k)}$. The left hand side of 
(16) is then $f(k) \, g(a(k))$ or 0 according to $ba=c$ or $ba \ne c$ which 
is the same as the right hand side.

\smallskip

\noindent Let us now describe the antipode $S$. The counit is given by
$$
\ve \, (h \, X_k) = h \, (1) \qquad 1 \ \hbox{unit of} \ G_2 \, . \leqno 
(17)
$$
We can consider the Hopf subalgebra ${\Hc}^1$ of ${\Hc}$ given by the $h \, 
X_k$, for $k=1$. The antipode $S^1$ of ${\Hc}^1$ is given by (group case)
$$
(S^1 \, h) (a) = h \, (a^{-1}) = \wt h \, (a) \, . \leqno (18)
$$
Thus it is natural to expect,
$$
S \, (\ve_a \, X_k) = X_{a(k)}^{-1} \, \ve_{a^{-1}} \, . \leqno (19)
$$
One needs to check that given $c \in G_2$ one has
$$
\sum_{ba=c} \ S \, (\ve_a \, X_k) \, \ve_b \, X_{a(k)} = \sum \, \ve_a \, 
X_k \, S \, (\ve_b \, X_{a(k)}) = \ve \, (\ve_c \, X_k) \, . \leqno (20)
$$
The first term is ${\displaystyle \sum_{ba=c}} \, X_{a(k)}^{-1} \, 
\ve_{a^{-1}} \, \ve_b \, X_{a(k)}$ which is 0 unless $c=1$. When $c=1$ it 
is equal to 1 since the $(a^{-1}) \cdot a(k) = (a \cdot k)^{-1}$ label 
$G_2$ when $a$ varies in $G_2$.

\smallskip

\noindent Similarly the second term gives $\ve_a \, X_k \, X_{b(a(k))}^{-1} 
\, \ve_{b^{-1}}$ which is non zero only if $a \cdot k \, c (k)^{-1} = 
b^{-1}$, i.e. $a \cdot k = b^{-1} \cdot c(k)$, i.e. if $c \cdot k = 1$ 
(since by lemma 1, one has $(b^{-1} \, c) \cdot k = ((b^{-1}) \cdot c(k)) 
(c \cdot k))$. Thus $c=1$ and the sum gives 1.

\smallskip

\noindent Let us now compute the ${\Hc}$-bimodule structure of ${\Hc}^*$.

\medskip

\noindent {\bf Lemma 2.} a) {\it The left action of $h \, X_k \in {\Hc}$ on 
$f \, U_a^* \in {\Hc}^*$ is given by $(h \, X_k) \cdot (f \, U_a^*) = h_a \, 
X_k (f) \, U_a^*$, where for $k_1 \in G_1$, $h_a (k_1) = h (a \cdot k_1)$ 
and $X_k (f) (k_1) = f(k_1 \, k)$.}

\noindent b) {\it The right action of $h \, X_k \in {\Hc}$ on $f \, U_a^* 
\in {\Hc}^*$ is given by $(f \, U_a^*) \cdot (h \, X_k) = h(a) \, L_{k^{-1}} 
(f) \, U_{a \cdot k}^*$, where for $k_1 \in G_1$, $(L_{k^{-1}} \, f) (k_1) = f 
(k \, k_1)$.}

\medskip

\noindent {\it Proof.} a) By definition $\lgl (h \, X_k) \cdot f \, U_a^* , 
h_0 \, X_{k_0} \rgl = \lgl f \, U_a^* , h_0 \, X_{k_0} \, h \, X_k \rgl$. 
Thus one has to check that $f( k_0 \, k) \, h_0 (a) \, h(a \cdot k_0) = h_0 
(a) \, (h_a \, X_k (f)) (k_0)$ which is clear.

\noindent b) One has $\lgl (f \, U_a^*) \cdot h \, X_k , h_0 \, X_{k_0} \rgl 
= \lgl f \, U_a^* , h \, X_k \, h_0 \, X_{k_0} \rgl = f(k \, k_0) \, h(a)$ 
$h_0 (a \cdot k)$ while $\lgl h(a) \, L_{k^{-1}} (f) \, U_{a \cdot k}^* , 
h_0  \, X_{k_0} \rgl = h(a) \, f(k \, k_0) \, h_0 (a \cdot k)$.~\xx

\vglue 1cm

\noindent {\bf VI.  Duality between ${\Hc}$ and $C_c^{\ify} (G_1) \semi G_2$, 
$G 
= G_1 \, G_2 = {\rm Diff} \, \Rb$}

\smallskip

Let, as above, ${\Hc}$ be the Hopf algebra generated by $X$, $Y$, $\d_n$. 
While there is a formal group $G({\Ac}^1)$ associated to the subalgebra 
${\Ac}^1$ of the Lie algebra of formal vector fields, there is no such group 
associated to ${\Ac}$ itself. As a substitute for this let us take,
$$
G = {\rm Diff} \, (\Rb) \, . \leqno (1)
$$
(We take smooth ones but restrict to real analytic if necessary).

\smallskip

\noindent We let $G_1 \sbs G$ be the subgroup of affine diffeomorphisms,
$$
k(x) = ax + b \qquad \fl \, x \in \Rb \leqno (2)
$$
and we let $G_2 \sbs G$ be the subgroup,
$$
\vp \in G \, , \ \vp (0) = 0 \, , \ \vp' (0) = 1 \, . \leqno (3)
$$
Given $\vp \in G$ it has a unique decomposition $\vp = k \, \psi$ where $k 
\in G_1$, $\psi \in G_2$ and one has,
$$
a = \vp' (0) \, , \ b = \vp (0) \, , \ \psi (x) = {\vp (x) - \vp (0) \over 
\vp' (0)} \, . \leqno (4)
$$
The left action of $G_2$ on $G_1$ is given by applying (4) to $x \ra \vp 
(ax+b)$, for $\vp \in G_2$. This gives
$$
b' = \vp (b) \, , \ a' = a \, \vp' (b) \leqno (5)
$$
which is the natural action of $G_2$ on the frame bundle $F(\Rb)$. Thus,

\medskip

\noindent {\bf Lemma 1.} {\it The left action of $G_2$ on $G_1$ coincides 
with the action of $G_2$ on $F(\Rb)$.}

\medskip

\noindent Let us then consider the right action of $G_1$ on $G_2$. In fact 
we consider the right action of $\vp_1 \in G$ on $\vp \in G_2$, it is given 
by
$$
(\vp \cdot \vp_1) (x) = {\vp (\vp_1 (x)) - \vp (\vp_1 (0)) \over \vp' (\vp_1 
(0)) \, \vp'_1 (0)} \, . \leqno (6)
$$

\smallskip

\noindent {\bf Lemma 2.} a) {\it The right action of $G$ on $G_2$ is affine 
in the coordinates $\d_n$ on $G_2$.}

\noindent b) {\it When restricted to $G_1$ it coincides with the action of 
the Lie algebra $X$, $Y$.}

\medskip

\noindent {\it Proof.} a) By definition one lets
$$
\d_n (\psi) = (\log \, \psi')^{(n)} \, (0) \, . \leqno (7)
$$
With $\psi = \vp \cdot \vp_1$, the first derivative $\psi' (x)$ is $\vp' 
(\vp_1 (x)) \, \vp'_1 (x) / \vp' (\vp_1 (0)) \, \vp'_1 (0)$ so that up to a 
constant one has,
$$
\log \, \psi' (x) = (\log \, \vp') \, (\vp_1 (x)) + \log \, \vp'_1 (x) \, . 
\leqno (8)
$$
Differentiating $n$ times the equality (8) proves a).

\noindent To prove b) let $\vp_1 (x) = ax+b$ while, up to a constant,
$$
(\log \, \vp') (x) = \sum_1^{\ify} \ {\d_n^0 \over n!} \, x^n \quad , \quad 
\d_n^0 = \d_n (\vp) \, . \leqno (9)
$$
Then the coordinates $\d_n = \d_n (\vp \cdot \vp_1)$ are obtained by 
replacing $x$ by $ax+b$ in (9), which gives
$$
\d_n = a^n \, \d_n^0 \ \hbox{if} \ b=0 \, , \ {\part \over \part b} \, \d_n 
= \d_{n+1}^0 \ \hbox{at} \ b=0 \, , \ a=1 \, . \leqno (10)
$$
\hfill \xx

\smallskip

\noindent We now consider the discrete crossed product of $C_c^{\ify} (G_1)$ 
by $G_2$, i.e. the algebra of finite linear combinations of terms
$$
f \, U_{\psi}^* \ , \ f \in C_c^{\ify} (G_1) \ , \ \psi \in G_2 \leqno (11)
$$
where the algebraic rules are 
$$
U_{\psi}^*  \ f = (f \circ \psi) \, U_{\psi}^* \, . \leqno (12)
$$
We want to define a pairing between the (envelopping) algebra ${\Hc}$ and 
the crossed product $C_c^{\ify} (G_1) \semi G_2$, by the equality
$$
\lgl h \, X_k , f \, U_{\psi}^* \rgl = h (\psi) \, f(k) \qquad \fl \, k \in 
G_1 \, , \ \psi \in G_2 \, . \leqno (13)
$$
In order to make sense of (13) we need to explain how we write an element of 
${\Hc}$ in the form $h \, X_k$.

\smallskip

\noindent Given a polynomial $P (\d_1 , \ldots , \d_n)$, we want to view it 
as a function on $G_2$ in such a way that the left action of that function 
$h$ given by lemma 2 of section V coincides with the multiplication of 
$U_{\psi}^*$ by
$$
P (\g_1 , \ldots , \g_n) \, , \ \g_j = \left( {\part \over \part x} 
\right)^j \, \log \, \psi' (x) \, e^{-js} \, , \ k = (e^{-s} , x ) \in G_1 
\, . \leqno (14)
$$
The formula of lemma 2 of section V  gives the multiplication by
$$
h (\psi \cdot k) \leqno (15)
$$
which shows that with $\d_n$ defined by (7) one has,
$$
h = P (\d_1 , \ldots ,\d_n) \, . \leqno (16)
$$
We then need to identify the Lie algebra generated by $X$, $Y$ with the Lie 
algebra ${\bf G}_1$ of $G_1$ (generated by $Z_{-1}$, $Z_0$) in such a way 
that the left action of the latter coincides with
$$
X \, f \, U_{\psi}^* = \left( e^{-s} \, {\part \over \part x} \, f \right) 
\, U_{\psi}^* \ , \ Y \, f \, U_{\psi}^* = - \left( {\part \over \part s} \, 
f \right) \, U_{\psi}^* \qquad (k = (e^{-s} , x)) \, . \leqno (17)
$$
The formula of lemma 2 of section V  gives $X_k \, f \, U_{\psi}^* = (X_k \, f) 
\, 
U_{\psi}^*$ with $(X_k \, f)$ $(k_1) = f (k_1 \, k)$. One has $f(k) = f (s,x)$ 
for $k = (e^{-s} , x)$, i.e. $k(t) = e^{-s} \, t + x$. With $k_1 = (e^{-s_1} , 
x_1)$ and $k (\ve) = (e^{\ve} , 0)$ one gets
$$
{\part \over \part \ve} \, f (k_1 \, k)_{\ve = 0} = - {\part \over \part 
s_1} \, f(k_1) = ( Yf) (k_1) \leqno (18)
$$
so that $Y$ corresponds to the one parameter subgroup $(e^{\ve} , 0)$ of 
$G_1$. With $k(\ve) = (1,\ve)$ one has
$$
{\part \over \part \ve} \, f (k_1 \, k_{\ve})_{\ve = 0} = (e^{-s} \, 
\part_{x_1} \, f) (k_1) = (Xf) (k_1) \leqno (19)
$$
so that $X$ corresponds to the one parameter subgroup $(1,\ve)$ of $G_1$. 
Now the element $e^{tX} \, e^{sY}$ of $G_1$ considered in section IV  is given 
by
$$
k = e^{tX} \, e^{sY} = (e^s , t) \leqno (20)
$$
which has the effect of changing $s$ to $-s$ in our formulas and thus 
explains the equality (11) of section IV.

\smallskip

\noindent This gives a good meaning to (13) as a pairing between ${\Hc}$ and 
the crossed product $C_c^{\ify} (G_1) \semi G_2$.

\vglue 1cm

\noindent {\bf VII.  Hopf algebras  and 
cyclic cohomology}

\smallskip

Let us first make sense of the right action of ${\Hc}$ on $C_c^{\ify} (G_1) 
\semi G_2$. We use the formula of lemma 2.b) section V
$$
(f \, U_{\psi}^*) \cdot (h \, X_k) = h(\psi) \, L_{k^{-1}} (f) \, U_{\psi 
\cdot k}^* \leqno (1)
$$
where $(L_{k^{-1}} (f)) (k_1) = f (k \, k_1)$ for $k_1 \in G_1$.

\smallskip

\noindent For the action of functions $h = P (\d_1 , \ldots , \d_n)$ we see 
that the difference with the left action is that we multiply $U_{\psi}^*$ by 
a {\it constant}, namely $h (\psi)$. Next, since we took a {\it discrete} 
crossed product to get $C_c^{\ify} (G_1) \semi G_2$, we can only act by the 
same type of elements on the right, i.e. by elements in
$$
\hbox{$\wt{\Hc}$ = {\it finite linear combinations of}} \ h \, X_k \, . \leqno (2)
$$
The algebra $\wt{\Hc}$ has little in common with ${\Hc}$, but both are 
multipliers of the smooth crossed product by $G_1$.

\smallskip

\noindent In fact, $\wt{\Hc}$ acts on both sides on $C_c^{\ify} (G_1) \semi 
G_2$ but only the left action makes sense at the Lie algebra level, i.e. as 
an action of ${\Hc}$.

\smallskip

\noindent The coproduct $\D$ is not defined for $\wt{\Hc}$ since $\D 
(e^{tX})$ cannot be written in $\wt{\Hc} \ot \wt{\Hc}$. Thus there is a 
problem to make sense of the right invariance property of an $n$-cochain,
$$
\vp (x^0 , \ldots , x^n) \ ; \ x^j \in C_c^{\ify} (G_1) \semi G_2 \leqno (3)
$$
which we would usually write as
$$
\sum \, \vp (x^0 \, y_{(0)} , \ldots , x^n \, y_{(n)}) = \ve (y) \, \vp (x^0 , 
\ldots , x^n) \leqno (4)
$$
for $\D^n \, y = \sum y_{(0)} \ot \ldots \ot y_{(n)} \ , \ y \in \wt{\Hc}$.

\smallskip

\noindent In fact it is natural to {\it require} as part of the right 
invariance invariance property of the cochain, that it possesses the right 
continuity property in the variables $\psi_j \in G_2$ so that the 
integration required in the coproduct formula (3) section IV, does make sense.

\smallskip

\noindent This problem does not arise for $n=0$, in which case we define the 
functional,
$$
\tau_0 (f \, U_{\psi}^*) = 0 \ \hbox{if} \ \psi \ne 1 \, , \ \tau_0 (f) = 
\int f(s,x) \, e^s \, ds \, dx \, , \leqno (5)
$$
where we used $k = (e^{-s} , x) \in G_1$.

\smallskip

\noindent One has $(f \circ \psi) (s,x) = f (s - \log \, \psi' (x) , \psi 
(x))$ by (5) section VI, so that $\tau_0$ is a {\it trace} on the algebra 
$C_c^{\ify} (G_1) \semi G_2 = {\Hc}_*$.

\smallskip

\noindent Let us compute $\tau_0 ((f \, U_{\psi}^*) (h \, X_k))$ and compare 
it with $\ve (h \, X_k) \, \tau_0 (f \, U_{\psi}^*)$. First $f \, U_{\psi}^* 
\, h \, X_k = h(\psi) (L_{k^{-1}} \, f) \, U_{\psi \cdot k}^*$ so that 
$\tau_0$ vanishes unless $\psi \cdot k = 1$ i.e. unless $\psi = 1$. We can 
thus assume that $\psi = 1$. Then we just need to compare $\tau_0 (L_{k^{-1}} 
\, f)$ with $\tau_0 (f)$. For $k = (e^{-s_1} , x_1)$ one has $k \circ 
(e^{-s} , x) (t) = e^{-s_1} (e^{-s} \, t + x) + x_1 = e^{-(s + s_1)} \, t + 
(e^{-s_1} \, x + x_1)$. This corresponds to $\psi (x) = e^{-s_1} \, x + x_1$ 
and preserves $\tau_0$. Thus

\medskip

\noindent {\bf Lemma 1.} {\it $\tau_0$ is a right invariant trace on 
${\Hc}_* = C_c^{\ify} (G_1) \semi G_2$.}

\medskip

\noindent Let us now introduce a bilinear pairing between ${\Hc}^{\ot 
(n+1)}$ and ${\Hc}_*^{\ot (n+1)}$ by the formula,
$$
\lgl y_0 \ot \ldots \ot y_n \ , \ x^0 \ot x^1 \ot \ldots \ot x^n \rgl = 
\tau_0 (y_0 (x^0) \ldots y_n (x^n))  \leqno (6)
$$
$\fl \, y_j \in {\Hc} \ , \ x^k \in {\Hc}_*$.

\smallskip

\noindent This pairing defines a corresponding weak topology and we let

\medskip

\noindent {\bf Definition 2.} {\it An $n$-cochain $\vp \in C^n$ on the 
algebra ${\Hc}_*$ is} right invariant {\it iff it is in the range of the 
above pairing.}

\medskip

\noindent We have a natural linear map $\t$ from ${\Hc}^{\ot (n+1)}$ to 
right invariant cochains on ${\Hc}_*$, given by
$$
\t (y_0 \ot \ldots \ot y_n) (x^0 , \ldots , x^n) = \lgl y_0 \ot \ldots 
\ot y_n , x^0 \ot \ldots \ot x^n \rgl \leqno (7)
$$
and we investigate the subcomplex of the cyclic complex of ${\Hc}_*$ given 
by the range of $\t$.

\smallskip

\noindent It is worthwhile to lift the cyclic operations at the level of $$
\bigoplus_{n=0}^{\ify} \, {\Hc}^{\ot (n+1)}
$$
and consider $\t$ as a morphism of $\L$-modules.

\smallskip

\noindent Thus let us recall that the basic operations in the cyclic complex 
of an algebra are given on cochains $\vp (x^0 , \ldots , x^n)$ by,
$$
\matrix{
(\d_i \, \vp) (x^0 , \ldots , x^n) &=& \vp (x^0 , \ldots , x^i \, x^{i+1} , 
\ldots , x^n) \hfill &i=0,1,\ldots , n-1 \hfill \cr
\cr
(\d_n \, \vp) (x^0 , \ldots , x^n) &=& \vp (x^n \, x^0 , x^1 , \ldots , 
x^{n-1}) \hfill \cr
\cr
(\s_j \, \vp) (x^0 , \ldots , x^n) &=& \vp (x^0 , \ldots , x^j , 1 , x^{j+1} 
, \ldots , x^n) \hfill &j=0,1,\ldots , n \hfill \cr
\cr
(\tau_n \, \vp) (x^0 , \ldots , x^n) &=& \vp (x^n , x^0 , \ldots , x^{n-1}) 
\, . \hfill \cr
} \leqno (8)
$$
These operations satisfy the following relations
$$
\matrix{
\tau_n \, \d_i = \d_{i-1} \, \tau_{n-1} &1 \leq i \leq n , &\tau_n \, \d_0 = 
\d_n \hfill \cr
\cr
\tau_n \, \s_i = \s_{i-1} \, \tau_{n+1} &1 \leq i \leq n , &\tau_n \, \s_0 = 
\s_n \, \tau_{n+1}^2 \cr
\cr
\tau_n^{n+1} = 1_n \, . \hfill \cr
} \leqno (9)
$$
In the first line $\d_i : C^{n-1} \ra C^n$. In the second line $\s_i : 
C^{n+1} \ra C^n$. Note that $(\s_n \, \vp) (x^0 , \ldots , x^n) = \vp (x^0 , 
\ldots , x^n ,1)$, $(\s_0 \, \vp) (x^0 , \ldots , x^n) = \vp (x^0 , 1 , x^1 , 
\ldots ,$ $x^n)$. The map $\t$ maps ${\Hc}^{\ot (n+1)}$ to $C^n$ thus there is 
a shift by 1 in the natural index $n$. We let
$$
\d_i (h^0 \ot h^1 \ot \ldots \ot h^i \ot \ldots \ot h^{n-1}) = h^0 \ot \ldots 
\ot h^{i-1} \ot \D \, h^i \ot h^{i+1} \ot \ldots \ot h^{n-1} \leqno (10)
$$
and this makes sense for $i=0,1,\ldots ,n-1$.

\smallskip

\noindent One has $\tau_0 (h^0 (x^0) \ldots h^i (x^i \, x^{i+1}) \, h^{i+1} 
(x^{i+2}) \ldots h^{n-1} (x^n)) = \sum \, \tau_0 (h^0 (x^0) \ldots$ 
$h_{(0)}^i \, (x^i) \, h_{(1)}^i (x^{i+1}) \ldots h^{n-1} (x^n))$
$$
\d_n (h^0 \ot h^1 \ot \ldots \ot h^{n-1}) = \sum \, h_{(1)}^0 \ot h^1 \ot 
\ldots \ot h^{n-1} \ot h_{(0)}^0 \leqno (11)
$$
which is compatible with $h^0 (x^n \, x^0) = \sum \, h_{(0)}^0 \, (x^n) \, 
h_{(1)}^0 \, (x^0)$, together with the trace property of $\tau_0$.

\smallskip

\noindent With $\ve : {\Hc} \ra \Cb$ the counit, we let
$$
\s_j (h^0 \ot \ldots \ot h^{n+1}) = h^0 \ot \ldots \ot \ve (h^{j+1}) \, 
h^{j+2} \ldots \ot h^{n+1} \quad j=0,1,\ldots ,n \leqno (12)
$$
which corresponds to $h^j (1) = \ve (h^j) \, 1$.

\smallskip

\noindent Finally we let $\tau_n$ act on ${\Hc}^{\ot (n+1)}$ by
$$
\tau_n (h^0 \ot \ldots \ot h^n) = h^1 \ot h^2 \ot \ldots \ot h^{n-1} \ot h^n 
\ot h^0 \leqno (13)
$$
which corresponds to $\tau_0 (h^0 (x^n) \, h^1 (x^0) \ldots h^n (x^{n-1})) = 
\tau_0 (h^1 (x^0) \ldots h^0 (x^n))$. One checks that with these operations 
${\Hc}_{\natural}$ is a $\L$-module where $\L$ is the cyclic category. To the 
relations (9) one has to add the relations of the simplicial $\D$, namely,
$$
\d_j \, \d_i = \d_i \, \d_{j-1} \ \hbox{for} \ i < j \, , \ \s_j \, \s_i = 
\s_i \, \s_{j+1} \qquad i \leq j \leqno (14)
$$
$$
\s_j \, \d_i = \left\{ \matrix{
\d_i \, \s_{j-1} \hfill &i < j \hfill \cr
1 \hfill &\hbox{if} \ i=j \ \hbox{or} \ i = j+1 \cr
\d_{i-1} \, \s_j \hfill &i > j+1 \, . \hfill \cr
} \right.
$$
The small category $\L$ is best defined as a quotient of the following 
category $E \, \L$. The latter has one object $(\Zb , n)$ for each $n$ and the 
morphisms $f : (\Zb , n) \ra (\Zb , m)$ are non decreasing maps, $(n,m \geq 
1)$
$$
f : \Zb \ra \Zb \ , \ f(x+n) = f(x)+m \qquad \fl \, x \in \Zb \, . \leqno (15)
$$
In defining $\L$ (cf. [Co]) one uses homotopy classes of non decreasing maps 
from $S^1$ to $S^1$ of degree 1, mapping $\Zb / n$ to $\Zb / m$. Given such a 
map we can lift it to a map satisfying (15). Such an $f$ defines uniquely a 
homotopy class downstairs and, if we replace $f$ by $f+km$, $k \in \Zb$ the 
result downstairs is the same. When $f(x) = a \, (m) \quad \fl \, x$, one can 
restrict $f$ to $\{0,1, \ldots , n-1 \}$

$$
Figure
$$

\noindent then $f(j)$ is either $a$ or $a+m$ which labels the various choices. 
One has $\L = (E \, \L) / \Zb$.

\smallskip

\noindent We recall that $\d_i$ is the injection that misses $i$, while $\s_j$ 
is the surjection which identifies $j$ with $j+1$.

\medskip

\noindent {\bf Proposition 3.} {\it ${\Hc}_{\natural}$ is a $\L$-module and $\t$ 
is a $\L$-module morphism to the $\L$-module $C^* ({\Hc}_*)$ of cochains on 
${\Hc}_* = C_c^{\ify} (G_1) \semi G_2$.}

\medskip

\noindent This is clear by construction.

\smallskip

\noindent Now the definition of ${\Hc}_{\natural}$ only involves ((10) 
$\ldots$ (13)) the coalgebra structure of $\Hc$, it is thus natural to compare 
it with the more obvious duality which pairs ${\Hc}^{\ot (n+1)}$ with 
${\Hc}_*^{\ot (n+1)}$ namely,
$$
\lgl h^0 \ot h^1 \ot \ldots \ot h^n , x^0 \ot \ldots \ot x^n \rgl = \prod_0^n 
\ \lgl h^j , x^j \rgl \, . \leqno (16)
$$
One has $\lgl h^i , x^i \, x^{i+1} \rgl = \lgl \D \, h^i , x^i \ot x^{i+1} 
\rgl$, so that the rules (10) and (11) are the correct ones. One has $\lgl 
h^{j+1} , 1 \rgl = \ve (h^{j+1})$ so that (12) is right. Finally (13) is also 
right.

\smallskip

\noindent This means that $C^* ({\Hc}_*) = {\Hc}_{\natural}^{**}$ as 
$\L$-modules. Thus,
$$
\t : {\Hc}_{\natural} \ra {\Hc}_{\natural}^{**} \leqno (17)
$$
is a cyclic morphism.

\smallskip

\noindent To understand the algebraic nature of $\t$, let us compute it in the 
simplest cases first. We first take ${\Hc} = \Cb \, G_1$ where $G_1$ is a 
finite group, and use the Hopf algebra ${\Hc} (G)$ for $G = G_1$, $G_2 = \{ e 
\}$. Thus as an algebra it is the group ring $\sum \, \lb_k \, X_k$, $k \in 
G_1$. As a right invariant trace on ${\Hc}^*$ we take
$$
\tau_0 (f) = \sum_{G_1} \, f(k) \ , \ \fl \, f \in {\Hc}^* \, . \leqno (18)
$$
The pairing $\lgl h , f \rgl$, for $h = \sum \, \lb_k \, X_k$, $f \in {\Hc}^*$ 
is given by $\sum \, \lb_k \, f(k)$. The left action $h(f)$ is given by lemma 
2 section V, i.e.
$$
h(f) (x) = \sum \, \lb_k \, f(xk) \, . \leqno (19)
$$
Thus the two pairings are, for (16): ${\displaystyle \sum_{k_i}} \, h^0 (k_0) 
\, f_0 (k_0) \ldots h^n (k_n) \, f_n (k_n)$ and for (6), ${\displaystyle 
\sum_{k,k_i}} \, h^0 (k_0) \, f_0 (k \, k_0) \ldots h^n (k_n) \, f_n (k \, 
k_n)$. 

\smallskip

\noindent Thus at the level of the $h^0 \ot \ldots \ot h^n$ the map $\t$ is 
just the sum of the left translates,
$$
\sum_{G_1} \, (L_g \ot L_g \ot \ldots \ot L_g) \, . \leqno (20)
$$
Next, we take the dual case, $G_1 = \{ e \}$, $G_2 = G$, with $G$ finite as 
above. Then ${\Hc}$ is the algebra of functions $h$ on $G_2$, and the dual 
${\Hc}^*$ is the group ring of $G_2^{\rm op} \sm G_2$, with generators 
$U_g^*$, $g \in G_2$. For a trace $\tau_0$ on this group ring, the right 
invariance under ${\Hc}$ means that $\tau_0$ is the regular trace,
$$
\tau_0 \left( \sum \, f_g \, U_g^* \right) = f_e \, . \leqno (21)
$$
This has a natural normalization, $\tau_0 (1) = 1$, for which we should expect 
$\t$ to be an idempotent. The pairing between $h$ and $f$ is,
$$
\lgl h , f \rgl = \sum \, h(g) \, f_g \, . \leqno (22)
$$
Thus the two pairings (16) and (6) give respectively, for (16) $\sum \, h_0 
(g_0) \, f_0(g_0)$ $h_1 (g_1) \, f_1 (g_1) \ldots h_n (g_n) \, f_n (g_n)$ and 
for (6), knowing that $h (U_g^*) = h(g) \, U_g^*$, i.e. $h \, \sum \, f_g \, 
U_g^* = \sum \, h(g) \, f_g \, U_g^*$ one gets,
$$
\sum_{g_n \ldots g_1 g_0 = 1} \, h_0 (g_0) \, f_0 (g_0) \, h_1 (g_1) \, f_1 
(g_1) \ldots h_n (g_n) \, f_n (g_n) \, . \leqno (23)
$$
Thus, at the level of ${\Hc}^{\ot n+1}$ the map $\t$ is exactly the 
localisation on the conjugacy class of $e$.

\smallskip

\noindent These examples clearly show that in general ${\rm Ker} \, \t \ne \{ 
0 \}$. Let us compute in our case how $\tau_0$ is modified by the {\it left} 
action of ${\Hc}$ on ${\Hc}_*$. By lemma 2 section V one has $(h \, X_k) (f \, 
U_{\psi}^*) = h_{\psi} \, X_k (f) \, U_{\psi}^*$ and $\tau_0$ vanishes unless 
$\psi = 1$. In this case $h_{\psi}$ is the constant $h(1) = \ve (h)$, while
$$
X_k (f) (k_1) = f (k_1 \, k) \, . \leqno (24)
$$
Thus we need to compare $\int f((e^{-s} , x) (a,b)) \, e^s \, ds \, dx$ with 
its value for $a=1$, $b=0$. With $k^{-1} = (a^{-1} , -b/a)$ the right 
multiplication by $k^{-1}$ transforms $(y,x)$ to $(y' , x')$ with $y' = y \, 
a^{-1}$, $x' = x - y \, b/a$, so that ${dy' \over y'^2} \wdg dx' = a \, {dy 
\over y^2} \wdg dx$,
$$
\tau_0 ((h \, X_k) \cdot f \, U_{\psi}^*) = \ve (h) \, \d (k) \, \tau_0 (f \, 
U_{\psi}^*) \leqno (25)
$$
where the module $\d$ of the group $G_1$ is,
$$
\d (a,b) = a \, . \leqno (26)
$$
In fact we view $\d$ as a character of $\Hc$, with
$$
\d (h \, X_k) = \ve (h) \, \d (k) \, . \leqno (27)
$$
(Note that $1 \cdot k = 1$ for all $k \in G_1$ so that (27) defines a 
character of $\Hc$.) Thus in our case we have a (non trivial) character of 
$\Hc$ such that
$$
\tau_0 (y(x)) = \d (y) \, \tau_0 (x) \qquad \fl \, x \in {\Ac} \, , \ y \in 
{\Hc} \, . \leqno (28)
$$
In fact we need to write the invariance property of $\tau_0$ as a formula for 
integrating by parts. To do this we introduce the twisted antipode,
$$
\wt S (y) = \sum \, \d (y_{(0)}) \, S (y_{(1)}) \ , \ y \in {\Hc} \, , \ \D \, 
y = \sum \, y_{(0)} \ot y_{(1)} \, . \leqno (29)
$$
One has $\wt S (y) = S (\s (y))$ where $\s$ is the automorphism obtained by 
composing $(\d \ot 1) \circ \D : {\Hc} \ra {\Hc}$. One can view $\wt S$ as $\d 
* S$ in the natural product (cf.(38) sectionV) on the algebra of linear maps 
from the coalgebra 
${\Hc}$ to the algebra ${\Hc}$. Since $S$ is the inverse of the identity map, 
i.e. $I * S = S * I = \ve$, one has $(\d * S) * I = \d$, i.e.
$$
\sum \, \wt S (y_{(0)}) \, y_{(1)} = \d (y) \qquad \fl \, y \in {\Hc} \, . 
\leqno (30)
$$
The formula that we need as a working hypothesis on $\tau_0$ is,
$$
\tau_0 (y \, (a) \, b) = \tau_0 (a \, \wt S (y) \, (b)) \qquad \fl \, a,b \in 
{\Ac} \, , \ y \in {\Hc} \, . \leqno (31)
$$
Using this formula we shall now determine ${\rm Ker} \, \t$ purely 
algebraically. We let $h = \sum \, h_0^i \ot \ldots \ot h_n^i \in {\Hc}^{\ot 
(n+1)}$, we associate to $h$ the following element of ${\Hc}^{\ot (n)}$:
$$
t(h) = \sum_i (\D^{n-1} \, \wt S (h_0^i)) \, h_1^i \ot \ldots \ot h_n^i \leqno 
(32)
$$
where we used both the coproduct of ${\Hc}$ and the product of ${\Hc}^{\ot 
(n)}$ to perform the operations.

\medskip

\noindent {\bf Lemma 3.} {\it $h \in {\rm Ker} \, \t$ iff $t(h) = 0$.}

\medskip

\noindent {\it Proof.} Let us first show that $h - 1 \ot t(h) \in {\rm Ker} \, 
\t$ for any $h$. One can assume that $h = h^0 \ot \ldots \ot h^n$. Using (31) 
one has $\tau_0 (h^0 (x^0) \, h^1 (x^1) \ldots h^n (x^n)) = \tau_0 (x^0 (\wt S 
\, h^0) (h^1 (x^1) \ldots h^n (x^n))) = \tau_0 (x^0 (\wt S \, h^0)_{(0)} \, 
h^1 (x^1) (\wt S \, h^0)_{(1)} \, h^2 (x^2)$ $\ldots$ \break $(\wt S \, 
h^0)_{(n-1)} \, h^n (x^n))$. 

\smallskip

\noindent It follows that if $t(h) = 0$ then $h \in {\rm Ker} \, \t$. 
Conversely, let us show that if $ 1 \ot t(h) \in {\rm Ker} \, \t$ then $t(h) = 
0$. We assume that the Haar measure $\tau_0$ is faithful i.e. that
$$
\tau_0 (ab) = 0 \qquad \fl \, a \in {\Ac} \quad \hbox{implies} \ b = 0 \, . 
\leqno 
(33)
$$
Thus, with $\wt h = t(h)$, $\wt h = \sum \, \wt{h}_1^i \ot \ldots \ot 
\wt{h}_n^i$, one has
$$
\sum \, \wt{h}_1^i (x^1) \ldots \wt{h}_n^i (x^n) = 0 \qquad \fl \, x^j \in 
{\Ac} \, . \leqno (34)
$$
Applying the unit $1 \in {\Hc}$ to both sides we get,
$$
\sum \, \lgl \wt{h}_1^i , x^1 \rgl \ldots \lgl \wt{h}_n^i , x^n \rgl = 0 
\qquad \fl \, x^j \in {\Ac} \leqno (35)
$$
which implies that $\wt h = 0$ in ${\Hc}^{\ot n}$.~\xx

\medskip

\noindent {\bf Definition 4.} {\it The cyclic module $C^* ({\Hc})$ of a Hopf 
algebra $\Hc$ is the quotient of ${\Hc}_{\natural}$ by the kernel of $t$.}

\medskip

Note that to define $t$ we needed the module $\d : {\Hc} \ra \Cb$, but that 
any reference to analysis has now disappeared in the definition of $C^* 
({\Hc})$.

\smallskip

\noindent Note also that the construction of $C^* ({\Hc})$ uses in an 
essential way both the coproduct and the product of ${\Hc}$. As we shall see 
it provides a working definition of the analogue of Lie algebra cohomology in 
general. (We did assume however that $\tau_0$ was a trace. 
This is an unwanted 
restriction which should be removed by making use of the modular theory.)

\smallskip

\noindent When ${\Hc} = {\Uc} ({\bf G})$ is the 
(complex) envelopping algebra of a (real) Lie 
algebra ${\bf G}$, there is a natural interpretation of the Lie algebra cohomology,
$$
H^* ({\bf G} , \Cb) = H^* ({\Uc} ({\bf G}) , \Cb) \leqno (36)
$$
where the right hand side is the Hochschild cohomology with coefficients in 
the ${\Uc} ({\bf G})$-bimodule $\Cb$ obtained using the augmentation.
In general, 
given a Hopf algebra ${\Hc}$ we can dualise the construction of the Hochschild 
complex $C^n ({\Hc}^* , \Cb)$ where $\Cb$ is viewed as a bimodule on ${\Hc}^*$ 
using the augmentation, i.e. the counit of ${\Hc}^*$. This gives the following 
operations: ${\Hc}^{\ot (n-1)} \ra {\Hc}^{\ot n}$, defining a cosimplicial 
space
$$
\matrix{
\d_0 (h^1 \ot \ldots \ot h^{n-1}) \hfill &= &1 \ot h^1 \ot \ldots \ot h^{n-1} 
\, , \hfill \cr
\d_j (h^1 \ot \ldots \ot h^{n-1}) \hfill &= &h^1 \ot \ldots \ot \D \, h^j \ot 
\ldots \ot 
h^{n-1} \, , \hfill \cr
\d_n (h^1 \ot \ldots \ot h^{n-1}) \hfill &= &h^1 \ot \ldots \ot h^{n-1} \ot 1 
\, , \hfill \cr
\s_i (h^1 \ot \ldots \ot h^{n+1}) \hfill &= & h^1 \ot \ldots \ve (h^{i+1}) \ot 
\ldots \ot h^{n+1} \, , \hfill &0 \leq i \leq n \, . \cr
} \leqno (37)
$$

\smallskip

\noindent {\bf Proposition 5.} {\it The map $t$ is an isomorphism of 
cosimplicial spaces.}

\medskip

\noindent {\it Proof.} Modulo ${\rm Ker} \, t = {\rm Ker} \, \t$ any element 
of 
${\Hc}^{\ot (n+1)}$ is equivalent to an element of the form $\sum \, 1 \ot 
h_i^1 
\ot \ldots \ot h_i^n = \xi$. One has $t(\xi) = \sum \, h_i^1 \ot \ldots \ot 
h_i^n$. It is enough to show that the subspace $1 \ot {\Hc}^{\ot *}$ is a 
cosimplicial subspace isomorphic to (37) through $t$. Thus we let $h^0 = 1$ in 
the definition (10) of $\d_1$ and (11) of $\d_n$ and check that they give 
(37). Similarly for $\s_i$.~\xx

\medskip

\noindent This shows that the underlying cosimplicial space of the cyclic 
module $C^* ({\Hc})$ is a standard object of homological algebra attached to 
the coalgebra ${\Hc}$ together with $1$, $\D 1 = 1 \ot 1$. The essential new 
feature, due to the Hopf algebra structure is that this cosimplicial space 
carries a cyclic structure. The latter is determined by giving the action of 
$\tau_n$ which is,
$$
\tau_n (h^1 \ot \ldots \ot h^n) = (\D^{n-1} \, \wt S (h_1)) \, h^2 \ot \ldots 
\ot h^n \ot 1 \leqno (38)
$$
where one uses the product in ${\Hc}^{\ot n}$ and the twisted antipode $\wt 
S$. It is nontrivial to check directly that $(\tau_n)^{n+1} = 1$, for instance 
for $n=1$ this means that $\wt S$ is an involution, i.e. ${\wt S}^2 = 1$. Note 
that the antipode $S$ of the Hopf algebra of section III is {\it not} an 
involution, while $\wt S$ is one. The first two cases 
in which we shall compute the cyclic cohomology of  $\Hc$ are the
following.

\medskip

\noindent {\bf Proposition 6.} 1) {\it The periodic cyclic cohomology $H^* 
({\Hc})$, for ${\Hc} = {\Uc} ({\bf G})$ the envelopping algebra of a Lie 
algebra ${\bf G}$ is isomorphic to the Lie algebra homology $H_* ({\bf G} , 
\Cb_{\d})$ where $\Cb_{\d} = \Cb$ is viewed a 
${\bf G}$-module by using the modular function $\d$ of $G$.}

\noindent 2) {\it The periodic cyclic cohomology $H^* ({\Hc})$, for ${\Hc} = 
{\Uc} ({\bf G})_*$, is isomorphic to the Lie algebra cohomology of ${\bf G}$ 
with trivial coefficients, provided $G$ is an affine space in the coordinates 
of ${\Hc}$. This holds in the nilpotent case.}

\medskip

\noindent {\it Proof.} 1) One has a natural inclusion ${\bf G} \sbs {\Uc} 
({\bf G})$. Let us consider the corresponding inclusion of $\L^n \, {\bf G}$ 
in ${\Hc}^{\ot n}$, given by
$$
X_1 \wdg \ldots \wdg X_n \ra \sum \, (-1)^{\s} \, X_{\s (1)} \ot \ldots \ot 
X_{\s (n)} \, . \leqno (41)
$$
Let $b : {\Hc}^{\ot n} \ra {\Hc}^{\ot (n+1)}$ be the Hochschild coboundary, 
one has
$$
{\rm Im} \, b \op \L^n \, {\bf G} = {\rm Ker} \, b \qquad \fl \, n \, . \leqno 
(42)
$$
For $n=1$ one has $b(h) = h \ot 1 - \D h + 1 \ot h$ so that $b(h) = 0$ iff $h 
\in {\bf G}$. For $n=0$, $b(\lb) = \lb - \lb = 0$ so that $b=0$. In general 
the statement (42) only uses the cosimplicial structure, i.e. only the 
coalgebra structure of ${\Hc}$ together with the element $1 \in {\Hc}$. This 
structure is unaffected if we replace the Lie algebra structure of ${\bf G}$ 
by the trivial commutative one. More precisely let us define the linear 
isomorphism,
$$
\pi : S ({\bf G}) \ra {\Uc} ({\bf G}) \ , \ \pi (X^n) = X^n \qquad \fl \, X 
\in {\bf G} \, . \leqno (43)
$$
Then $\D \circ \pi = (\pi \ot \pi) \circ \D_S$ where $\D_S$ is the coproduct 
of $S ({\bf G})$. Indeed it is enough to check this equality on $X^n$, $X \in 
{\bf G}$ and both sides give $\sum \, C_n^k \, X^k \ot X^{n-k}$.

\smallskip

\noindent The result then follows by dualising the homotopy between the 
standard resolution and the Koszul resolution $S(E) \ot \L (E)$ of the module 
$\Cb$ over $S(E)$ for a vector space $E$ .

\smallskip

\noindent Let us then compute $B (X_1 \wdg \ldots \wdg X_n)$. Note that $B_0 
(X_1 \wdg \ldots \wdg X_n)$ corresponds to the functional $ (-1)^{\s} \tau_0 
(X_{\s (1)} 
(x^0) \ldots X_{\s (n)} 
(x^{n-1})) $ which is already cyclic. Thus it is enough to compute 
$B_0$. One has
$$
\wt S (X) = -X + \d (X) \qquad \fl \, X \in {\bf G} \, , \leqno (44)
$$
thus $\D^{n-2} \, \wt S \, X = - \sum \, 1 \ot \ldots \ot X \ot \ldots \ot 1 + 
\d (X) \, 1 \ot 1 \ldots \ot 1$.

\smallskip

\noindent We get $B_0 (X_1 \wdg \ldots \wdg X_n) = {\displaystyle \sum_{\s}} 
\, 
(-1)^{\s} \, \d (X_{\s (1)}) \, X_{\s (2)} \ot \ldots \ot X_{\s (n)} - 
{\displaystyle \sum_{\s}} \, (-1)^{\s}$ $X_{\s (2)} \ot \ldots \ot X_{\s (1)} 
\, X_{\s (j)} \ot \ldots \ot X_{\s (n)} = \sum \, (-1)^{k+1} \, \d (X_k) \, 
X_1 \wdg \ldots \wdg \build{X_k}_{}^{\!\!\vee} \wdg \ldots \wdg X_n + 
{\displaystyle \sum_{i < j}} \, (-1)^{i+j} \, [X_i , X_j] \wdg X_1 \wdg \ldots 
\wdg \build{X_i}_{}^{\!\!\vee} \wdg \ldots \wdg \build{X_j}_{}^{\!\!\vee} \wdg 
\ldots \wdg X_n$. This shows that $B$ leaves $\L^* \, {\bf G}$ invariant and 
coincides there with the boundary map of Lie algebra homology. The situation 
is identical to what happens in computing cyclic cohomology of the algebra of 
smooth 
functions on a manifold.

\smallskip

\noindent 2) The Hochschild complex of ${\Hc}$ is by construction the dual of 
the standard chain complex which computes the Hochschild homology of ${\Uc} 
({\bf G})$ with coefficients in $\Cb$ (viewed as a bimodule using $\ve$). 
Recall 
that in the latter complex the boundary $d_n$ is
$$
\matrix{
d_n (\lb_1 \ot \ldots \ot \lb_n) = \ve (\lb_1) \, \lb_2 \ot \ldots \ot \lb_n - 
\lb_1 \, \lb_2 \ot \ldots \ot \lb_n + \ldots \cr
\cr
+ (-1)^i \, \lb_1 \ot \ldots \ot \lb_i \, \lb_{i+1} \ot \ldots \ot \lb_n + 
\ldots + (-1)^n \, \lb_1 \ot \ldots \ot \lb_{n-1} \, \ve (\lb_n) \, . \cr
} \leqno (45)
$$
One has a homotopy between the complex (45) and the subcomplex $V({\bf G})$ 
obtained by the following map from $\L^* \, {\bf G}$,
$$
X_1 \wdg \ldots \wdg X_n \ra \sum \, (-1)^{\s} \, X_{\s (1)} \ot \ldots \ot 
X_{\s (n)} \in {\Uc} ({\bf G})^{\ot n} \, . \leqno (46)
$$
This gives a subcomplex on which $d_n$ coincides with the boundary in Lie 
algebra homology with trivial coefficients,
$$
d^t (X_1 \wdg \ldots \wdg X_n) = \sum_{i<j} \, (-1)^{i+j} \, [X_i , X_j] \wdg 
X_1 \wdg \ldots \wdg \build{X_i}_{}^{\!\!\vee} \wdg \ldots \wdg 
\build{X_j}_{}^{\!\!\vee} \wdg \ldots \wdg X_n \, . \leqno (47)
$$
It is thus natural to try and dualise the above homotopy to the Hochschild 
complex of $\Hc$. Now a Hochschild cocycle
$$
h = \sum \, h_1^i \ot \ldots \ot h_n^i \in {\Hc}^{\ot n} \leqno (48)
$$
gives an $n$-dimensional group cocycle on $G = \{ g \in {\Uc} ({\bf G})_{\rm 
completed} ; \D g = g \ot g \}$ where,
$$
c (g_1 , \ldots , g_n) = \sum \, \lgl h_1^i , g_1 \rgl \ldots \lgl h_n^i , g_n 
\rgl \, . \leqno (49)
$$
These cocycles are quite special in that they depend polynomially on the 
$g_i$'s. Thus we need to construct a map of cochain complexes from the complex 
of Lie algebra cohomology to the complex of polynomial cocycles and prove that 
it gives an isomorphism in cohomology,
$$
s : \L^n \, ({\bf G})^* \ra {\Hc}^{\ot n} \, . \leqno (50)
$$
If we let $j$ be the restriction of a polynomial cochain to $\L ({\bf G})^*$ 
we expect to have $j \circ s = {\rm id}$ and to have a homotopy,
$$
s \circ j - 1 = dk + kd \, . \leqno (51)
$$
Using the affine coordinates on $G$ we get for $g^0 , \ldots , g^n \in G$, an
affine simplex
$$
\D (g^0 , \ldots , g^n) = \left\{ \sum_0^n \, \lb_i \, g^i \ ; \ \lb_i \in 
[0,1] , \sum \, \lb_i = 1 \right\} \, , \leqno (52)
$$
moreover the right multiplication by $g \in G$ being affine, we have,
$$
\D (g^0 \, g , \ldots , g^n \, g) = \D (g^0 , \ldots , g^n) \, g \, .\leqno 
(53)
$$
The map $s$ is then obtained by the following formula,
$$
(s \, \om) (g^1 , \ldots , g^n) = \g (1, g^1 , g^1 \, g^2 , \ldots , g^1 
\ldots g^n) \leqno (54)
$$
where 
$$
\g (g^0 , \ldots , g^n) = \int_{\D (g^0 , \ldots , g^n)} \, {\bf w}
$$
where ${\bf w}$ is the right invariant form on $G$ associated to $\om \in \L^n 
({\bf G})^*$. To prove the existence of the homotopy (51) we introduce the 
bicomplex of the proof of the Van~Est theorem but we restrict to forms with 
polynomial coefficients: $A^n (G)$ and to group cochains which are polynomial. 
Thus an element of $C^{n,m}$ is an $n$-group cochain
$$
c (g^1 , \ldots , g^n) \in A^m \leqno (55)
$$
and one uses the right action of $G$ on itself to act on forms,
$$
(g \, \om) (x) = \om \, (xg) \, . \leqno (56)
$$
The first coboundary $d_1$ is given by
$$
\matrix{
(d_1 \, c) (g^1 , \ldots , g^{n+1}) = g^1 \, c (g^2 , \ldots ,g^{n+1}) - c 
(g^1 \, g^2 , \ldots , g^{n+1}) + \ldots \cr
\cr
+ (-1)^n \, c (g^1 , \ldots , g^n \, g^{n+1}) + (-1)^{n+1} \, c (g^1 , \ldots 
, g^n) \, . \cr
}\leqno (57)
$$
The second coboundary is simply
$$
(d_2 \, c) (g^1 , \ldots , g^n) = dc (g^1 , \ldots , g^n) \, . \leqno (58)
$$
(One should put a sign so that $d_1 \, d_2 = - d_2 \, d_1$.)

\smallskip

\noindent We need to write down explicitly the homotopies of lines and columns 
in order to check that they preserve the polynomial property of the cochains. 
In 
the affine coordinates $\d_{\mu}$ on $G$ we let
$$
X = \sum \, \d_{\mu} \, {\part \over \part \, \d_{\mu}} \leqno (59)
$$
be the vector field which contracts $G$ to a point.

\smallskip

\noindent Then the homotopy $k_2$ for $d_2$ comes from,
$$
k \, \om = \int_0^1 (i_X \, \om) (t \, \d) \, {dt \over t} \leqno (60)
$$
which preserves the space of forms with polynomial 
coefficients\footnote{$^{(*)}$}{For instance for $G = \Rb$, $\om = f(x) \, dx$ 
one gets $k \om = \a$, $\a (x) = x \int_0^1 f(tx)$ $dt$. In general take $\om 
= P(\d) \, d \, \d_1 \wdg \ldots \wdg d \, \d_k$, then $i_X \, \om$ is of the 
same form and one just needs to know that $\int_0^1 P (t \, \d) \, t^a \, dt$ 
is still a polynomial in $\d$, which is clear.}. The homotopy $k_1$ for $d_1$ 
comes from the structure of induced module, i.e. from viewing a cochain $c 
(g_1 , \ldots , g_n)$ as a function of $x \in G$ with values in $\L^m \, {\bf 
G}^*$,
$$
(k_1 \, c) (g_1 , \ldots , g_{n-1}) (x) = c (x, g_1 , \ldots , g_{n-1}) (e) \, 
. \leqno (61)
$$
One has $(d_1 \, k_1 \, c) (g_1 , \ldots , g_n) (x) = (k_1 \, c) (g_2 , \ldots 
, g_n) (x \, g_1) - (k_1 \, c) (g_1 \, g_2 , \ldots ,$ $g_n) (x) + \ldots + 
(-1)^{n-1} \, (k_1 \, c) (g_1 , \ldots , g_{n-1} \, g_n) (x) + (-1)^n \, (k_1 
\, c) (g_1 , \ldots , g_{n-1})$ $(x) = \{ c (x \, g_1 , g_2 , \ldots , 
g_n) - c (x,g_1 \, g_2 , \ldots , g_n) + \ldots + (-1)^{n-1} \, c (x,g_1 , 
\ldots , g_{n-1}$ $g_n) + (-1)^n \, c (x,g_1 ,\ldots , g_{n-1}) \} (e) = c 
(g_1 
,\ldots , g_n) (x) - (d_1 \, c) (x,g_1 ,\ldots , g_n) (e) = c (g_1 ,\ldots , 
g_n) (x) - (k_1 \, d_1 \, c) (g_1 ,\ldots , g_n) (x)$.

\smallskip

\noindent This homotopy $k_1$ clearly preserves the polynomial cochains. This 
is 
enough to show that the Hochschild cohomology of ${\Hc}$ is isomorphic to the 
Lie algebra cohomology $H^* ({\bf G} , \Cb)$. But it follows from the 
construction of the cocycle $s \, \om$ (54) that
$$
(s \, \om) (g_1 ,\ldots , g_n) = 0 \quad \hbox{if} \quad g_1 \ldots g_n = 1 \, 
, 
\leqno (62)
$$
and this implies (as in the case of discrete groups) that the corresponding 
cocycle is also cyclic.~\xx

\medskip

\noindent Our goal now is to compute the cyclic cohomology of our original 
Hopf 
algebra 
${\Hc}$. This should combine the two parts of Proposition 6. In the first part 
the Hochschild cohomology was easy to compute and the operator $B$ was non 
trivial. In the second part $b$ was non trivial. 
 At the level of $G_1$, one needs to transform the Lie algebra 
cohomology into the Lie algebra homology with coefficients in $\Cb_{\d}$. The 
latter corresponds to invariant currents on $G$ and the natural isomorphism is 
a Poincar\'e duality. For $G_1$ with Lie algebra $\{ X , Y \}$, $[Y,X] = X$, 
$\d (X) = 0$, $\d (Y) = 1$, one gets that $X \wdg Y$ is a 2-dimensional cycle, 
while since $b \, Y = 1$, there is no zero dimensional cycle. For the Lie 
algebra cohomology one checks that there is no 2-dimensional cocycle.
\smallskip

\noindent
We shall start by 
constructing 
an explicit map from the Lie algebra cohomology of ${\bf G} = {\Ac}$, the Lie 
algebra of formal vector fields, to the cyclic cohomology of ${\Hc}$.
As an 
intermediate step in the construction of this map, we shall use the 
following double complex $(C^{n,m} , d_1 , d_2)$. For $0 \leq k \leq \dim G_1$ 
and let $\Om_k (G_1) = \Om_k$ be the space of de~Rham currents on $G_1$, and 
we 
let $\Om_k = \{ 0 \}$ for $k \notin \{ 0, \ldots ,\dim G_1 \}$. We let 
$C^{n,m} 
= \{ 0 \}$ unless $n \geq 0$ and $-\dim G_1 \leq m \leq 0$, and let $C^{n,m}$ 
be the space of totally antisymmetric polynomial maps $\g : G_2^{n+1} \ra 
\Om_{-m}$ such that,
$$
\g (g_0 \, g , \ldots , g_n \, g) = g^{-1} \, \g (g_0 , \ldots ,g_n) \qquad 
\fl 
\, g_i \in G_2 \, , \ g \in G \leqno (63)
$$
where we use the right action of $G$ on $G_2 = G_1 \bsh G$ to make sense of 
$g_i \, g$ and the left action of $G$ on $G_1 = G / G_2$ 
to make sense of $g^{-1} \g $. 
The coboundary $d_1 : C^{n,m} \ra C^{n+1,m}$ is given by
$$
(d_1 \, \g) (g_0 , \ldots ,g_{n+1}) = (-1)^m \sum_{j=0}^{n+1} \ (-1)^j \, \g 
(g_0 , \ldots , \build{g_j}_{}^{\!\!\vee} , \ldots , g_{n+1}) \, . \leqno (64)
$$
The coboundary $d_2 : C^{n,m} \ra C^{n,m+1}$ is the de~Rham boundary,
$$
(d_2 \, \g) (g_0 , \ldots ,g_n) = d^t \, \g (g_0 , \ldots ,g_n) \, . \leqno 
(65)
$$
For $g_0 , \ldots ,g_n \in G_2$, we let $\D (g_0 , \ldots ,g_n)$ be the affine 
simplex with vertices the $g_i$ in the affine coordinates $\d_k$ on $G_2$. 
Since the right action of $G$ on $G_2$ is affine in these coordinates, we have
$$
\D (g_0 \, g , \ldots , g_n \, g) = \D (g_0 , \ldots ,g_n) \, g \qquad \fl \, 
g_i \in G_2 \, , \ g \in G \, . \leqno (66)
$$
Let $\om$ be a left invariant differential form on $G$ associated to a cochain 
of degree $k$ in the complex defining the Lie algebra cohomology of the Lie 
algebra $\Ac$. For each pair of integers $n \geq 0$, $-\dim G_1 
\leq 
m \leq 0$ such that $n+m = k - \dim G_1$, let
$$
\lgl C_{n,m} (g^0 , \ldots ,g^n) , \a \rgl = (-1)^{{m(m+1) \over 2}}\int_{(G_1 
\ts \D (g^0 , \ldots 
,g^n))^{-1}} \pi_1^* (\a) \wdg \om \leqno (67)
$$
for any smooth diferential form $\a$, with compact support on $G_1$ and of 
degree $-m$.

\smallskip

\noindent In this formula we use $G_1 \ts \D (g_0 , \ldots ,g_n)$ as a cycle 
in 
$G$ and we need to show that if $K_i \sbs G_i$ are compact subsets, 
the subset 
of $G$
$$
K = \{ g \in G \, ; \ \pi_1 (g) \in K_1 \, , \ \pi_2 (g^{-1}) \in K_2 \} 
\leqno 
(68)
$$
is compact.

\smallskip

\noindent For $g \in K$ one has $g = ka$ with $k \in K_1$ and $\pi_2 (g^{-1}) 
\in K_2$ so that $\pi_2 (a^{-1} \, k^{-1})$ $\in K_2$. But $\pi_2 (a^{-1} \, 
k^{-1}) = a^{-1} \cdot k^{-1}$ and one has $a^{-1} \in K_2 \cdot K_1^{-1}$ and 
$a \in (K_2 \cdot K_1^{-1})^{-1}$, thus, the required compactness follows from
$$
K \sbs K_1 (K_2 \cdot K_1^{-1})^{-1} \, . \leqno (69)
$$
We let $C^* ({\Ac})$ be the cochain complex defining the Lie algebra 
cohomology 
of $\Ac$ and let $C$ be the map defined by (67).

\medskip

\noindent {\bf Lemma 7.} {\it The map $C$ is a morphism to the total complex 
of 
$(C^{n,m} , d_1 , d_2)$.}

\medskip

\noindent {\it Proof.} Let us first check the invariance condition (63). One 
has (66) so that for $C_{n,m} (g^0 \, g , \ldots , g^n \, g)$ the integration 
takes place on $\{ h \in G \, ; \ \pi_2 (h^{-1}) \in \D (g^0 , \ldots ,g^n) \, 
g \} = \sum'$. Since $\pi_2 (h^{-1}) \, g^{-1} = \pi_2 (h^{-1} \, g^{-1})$ one 
has $\sum' = g^{-1} \sum$ in $G$, with $\sum = \{ h \in G \, ; \ \pi_2 
(h^{-1}) 
\in \D (g^0 , \ldots ,g^n)\}$. One has, with $\b = \pi_1^* (\a) \wdg \om$, the 
equality
$$
\int_{g^{-1} \sum} \ \b = \int_{\sum} \ L_g \, \b \leqno (70)
$$
where $g \ra L_g$ is the natural action of $G$ on forms on $G$ by left 
translation. One has $\pi_1 (gk) = g \, \pi_1 (k)$ so that $L_g \, \pi_1^* 
(\a) 
= \pi_1^* (L_g \, \a)$. Moreover $\om$ is left invariant by hypothesis, so one 
gets,
$$
\lgl C_{n,m} (g^0 \, g , \ldots , g^n \, g) , \a \rgl = \lgl C_{n,m} (g^0 , 
\ldots ,g^n) , L_g \, \a \rgl \leqno (71)
$$
which since $L_g^t = L_{g^{-1}}$, is the invariance condition (63).

\smallskip

\noindent Before we check that $C_{n,m}$ is polynomial in the $g^i$'s, let us 
check that
$$
C (d\om) = (d_1 + d_2) \, C(\om) \, . \leqno (72)
$$
One has $d (\pi_1^* (\a) \wdg \om) = \pi_1^* (d \a) \wdg \om + (-1)^m \, 
\pi_1^* 
(\a) \wdg d \om$ and since
$$
\int_{\pi_2 (g^{-1}) \in \D (g^0 , \ldots ,g^{n+1})} d\b = \sum (-1)^i 
\int_{\pi_2 (g^{-1}) \in \D (g^0 , \ldots , \build{g^i}_{}^{\!\!\vee} , \ldots 
, g^n ,g^{n+1})} \b \, . \leqno (73)
$$
With $\b = \pi_1^* (\a) \wdg \om$, the r.h.s. gives $(-1)^m \, (d_1 \, C) \, 
(g^0 , \ldots , g^{n+1})$, while the l.h.s. gives $\lgl C (g^0 , \ldots , 
g^{n+1}) , d\a \rgl + (-1)^m \lgl C' (g^0 , \ldots , g^{n+1}) , \a \rgl$, with 
$C' = C (d\om)$. Thus we get,
$$
d_1 \, C_{n,m} + d_2 \, C_{n+1,m-1} = C'_{n+1,m} \leqno (74)
$$
provided we use the sign: $(-1)^{{m(m+1) \over 2}}$ in the definition (67) 
of $C$.

\noindent We shall now be more specific on the polynomial expression of $C_{n,m} (g^0 , 
\ldots ,$ $g^n)$ and write this de~Rham current on $G_1$ in the form,
$$
\sum \, c_j \, \rho_j \leqno (75)
$$
where $\rho_0 = 1$, $\rho_1 = ds$, $\rho_2 = e^s \, dx$, $\rho = \rho_1 \wdg 
\rho_2$ form a basis of left $G_1$ invariant forms on $G_1$, while the $c_j$ 
are functions on $G_1$ which are finite linear combinations of finite products 
of the following functions,
$$
k \in G_1 \ra \d_p (g_j \cdot k) \ , \ j \in \{ 0,\ldots ,n \} \, . \leqno 
(76)
$$
The equality (67) defines $C_{n,m}$ as the integration over the fibers for the 
map $\pi_1 : G \ra G_1$ of the product of the smooth form $\om$ by the current 
$\wt \D$ of integration on $(G_1 \ts \D)^{-1} = \{ g \in G \, ; \, \pi_2 
(g^{-1}) \in \D \}$,
$$
C_{n,m} = (\pi_1)_* \, (\om \wdg \wt \D ) \, . \leqno (77)
$$
Thus $c_0$, which is a function, 
is obtained as the integral of $\om \wdg \wt \D$ 
along the fibers, and its value at $1 \in G_1$ is
$$
c_0 (1) = \int_{g^{-1} \in \D} \om \vert G_2 \, . \leqno (78)
$$
To obtain the value of $c_j (1)$ by a similar formula, one can contract by a 
vector field $Z$ on $G_1$ given by left translation, $k \in G_1 \ra  
\part_\ve\,(k(\ve)k)_ {\ve=0}  $,
 $k(\ve) \in G_1$. Let $\wt Z$ be the vector field on 
$G$ given by the same left translation, $g \in G_1 \ra  \part_\ve\,(k(\ve)g)_ 
{\ve=0}  $.
 The equality $\pi_1 (k(\ve) \, g) = k(\ve) \, \pi_1 (g)$ shows that 
$\wt Z$ is a left of $Z$ for the fibration $\pi_1$. It follows that for any 
current $\b$ on $G$ one has
$$
i_Z \, \pi_{1_*} (\b) = \pi_{1_*} (i_{\wt Z} \, \b) \, . \leqno (79)
$$
(One has $\lgl i_Z \, \pi_{1_*} (\b) , \a \rgl = \lgl \pi_{1_*} (\b) , i_Z \, 
\a \rgl = \lgl \b , \pi_1^* (i_Z \, \a) \rgl = \int \pi_1^* (i_Z \, \a) \wdg 
\b 
= \int i_{\wt Z} \, \pi_1^* (\a) \wdg \b = \int \pi_1^* (\a) \wdg i_{\wt Z} \, 
\b$.)

\smallskip

\noindent Next one has $i_{\wt Z} \, (\om \wdg \wt \D ) = (i_{\wt Z} \, \om) 
\ts \wt \D + (-1)^{\part \om} \, \om \wdg i_{\wt Z} \, \wt \D$. The 
contribution of the first term is simple,
$$
\int_{g^{-1} \in \D} (i_{\wt Z} \, \om) \vert G_2 \, . \leqno (80)
$$
To compute the contribution of the second term, one needs to understand the 
current $i_{\wt Z} \, \wt \D$ on $G$. In general, if $\b$ is the current of 
integration on a manifold $M \hra G$ (possibly with boundary), the contraction 
$i_{\wt Z} \, \b$ is obtained as a limit for $\ve \ra 0$ from the manifold $M 
\ts [-\ve , \ve]$ which maps to $G$ by $(x,t) \ra (\exp t \, \wt Z) (x)$. 
Applying this to $\wt \D$ one gets the map from $G_1 \ts \D \ts [-\ve , \ve]$ 
to $G$ given by
$$
(k,a,t) \ra k(t) \, a^{-1} \, k^{-1} \, . \leqno (81)
$$
One has $k(t) \, a^{-1} \, k^{-1} = (ka \, k(t)^{-1})^{-1} = k^1 (a \cdot 
k(t)^{-1}))^{-1}$ where the $a \cdot k^{-1} (t)$ belong to the simplex $\D 
\cdot k(t)^{-1}$ while $k^1 \in G_1$ is arbitrary.

\smallskip

\noindent Thus, if we let $Z'$ be the vector field on $G_2$ given by,
$$
a \ra  \part_\ve\,(a \cdot k^{-1} (\ve) )_ {\ve=0} \leqno (82)
$$
we see that the contribution of the second term is
$$
\int_{g \in \D} i_{Z'} \, \wt{\om} \vert G_2 \ , \ \wt{\om} \, (g) = \om (g^{-1}) 
\, . \leqno (83)
$$
Thus there exists a differential form $\mu_j$ on $G_2$ obtained from $\wt \om$ 
by contraction by suitable vector fields and restriction to $G_2$, such 
that 
$$
c_j (1) = \int_{\D (g^0 , \ldots , g^n)} \mu_j \, . \leqno (84)
$$
The value of $c_j$ at $k \in G_1$ can now be computed using (71) for $g = k 
\in 
G_1$. The forms $\rho_j$ are left invariant under the action of $G_1$ and by 
(71) one has $L_{g^{-1}} \, C (g_0 , \ldots , g_n) = C (g_0 \, g , 
\ldots , g_n 
\, g)$, $g=k$. Thus $c_j (k)$ corresponds to the current $L_{k^{-1}} \, C (g_0 ,  \ldots , g_n)$ evaluated at 1 and one has,
$$
c_j (k) = \int_{\D (g^0 , \ldots , g^n) \cdot k} \mu_j = \int_{\D (g^0 \cdot k 
, \ldots , g^n \cdot k)} \mu_j \, . \leqno (85)
$$
Let us now apply (71) for $g \in G_2$. The point $1 \in G_1$, is fixed by the 
action of $G_2$ while the forms $\rho_j$ vary as follows under the action of 
$G_2$,
$$
L_{g^{-1}} \, \rho_j = \rho_j \ , \ j \ne 1 \ , \ L_{g^{-1}} \, 
\rho_1 = \rho_1 - 
\d_1 (g \cdot k) \, \rho_2 \qquad ({\rm at} \ k \in G_1) \, . \leqno (86)
$$
Thus we see that while $\mu_j$, $j \ne 1$ are right invariant forms on $G_2$, 
the form $\mu_1$ satisfies,
$$
\mu_1 (ag) = \mu_1 (a) - \d_1 (g) \, \mu_2 (a) \, . \leqno (87)
$$
Now in $G_2$ the product rule as well as $g \ra g^{-1}$ are polynomial in the 
coordinates $\d_n$. It follows that the forms $\mu_j$ are polynomial forms in 
these coordinates and that the formula (84) is a polynomial function of the 
$\d_j (g^k)$. Using (85) we obtain the desired form (76) for $c_i$.~\xx

\medskip

\noindent We shall now use the canonical map $\Phi$ of [Co] Theorem~14 p.~220, 
from the bicomplex $(C^{n,m} , d_1 , d_2)$ to the $(b,B)$ bicomplex of the 
algebra ${\Hc}_* = C_c^{\ify} (G_1) \semi G_2$. What we need to prove is that 
the obtained cochains on ${\Hc}_*$ are right invariant in the sense of 
definition 2 above.

\smallskip

\noindent Let us first rewrite the construction of $\Phi$ using the notation 
$f 
U_{\psi}^*$ for the generators of ${\Hc}_*$. As in [Co], we let $\Bc$  be the 
tensor product,
$$
{\Bc} = A^* (G_1) \ot \L \, \Cb \, (G'_2) \leqno (88)
$$
where $A^* (G_1)$ is the algebra of smooth forms with compact support on 
$G_1$, 
while we label the generators of the exterior algebra $\L \, \Cb \, (G'_2)$ as
$\d_{\psi}$, $\psi \in G_2$, with $\d_0 = 0$. We take the crossed product,
$$
C = {\Bc} \semi G_2 \leqno (89)
$$
of $\Bc$ by the product action of $G_2$, so that
$$
\matrix{
U_{\psi}^* \, \om \, U_{\psi} = \psi^* \, \om = \om \circ \psi \hfill &\qquad 
\fl 
\, \om \in A^* (G_1) \, , \cr
\cr
U_{\psi_1}^* \, \d_{\psi_2} \, U_{\psi_1} =\d_{\psi_2 \circ \psi_1} - 
\d_{\psi_1} \hfill &\qquad \fl \, \psi_j \in G_2 \, . \hfill \cr
} \leqno (90)
$$
The differential $d$ in $C$ is given by
$$
d (b \, U_{\psi}^*) = db \, U_{\psi}^* - (-1)^{\part b} \, b \, \d_{\psi} \, 
U_{\psi}^* \leqno (91)
$$
where the first term comes from the exterior differential in $A^* (G_1)$. Thus 
the $\d_{\psi}$ play the role of
$$
- \d_{\psi} = (d \, U_{\psi}^*) \, U_{\psi} = - U_{\psi}^* \, d \, U_{\psi} \, 
. \leqno (92)
$$
A cochain $\g \in C^{n,m}$ in the above bicomplex determines a linear form 
$\wt{\g}$ on $C$, by,
$$
\wt{\g} (\om \ot \d_{g_1} \ldots \d_{g_n}) = \lgl \om , \g (1, g_1 , \ldots , 
g_n) \rgl \ , \ \wt{\g} (b \, U_{\psi}^*) = 0 \quad \hbox{if} \ \psi \ne 1 \ , 
. \leqno (93)
$$
What we shall show is that the following cochains on ${\Hc}_*$ satisfy 
definition 2,
$$
\vp (x^0 , \ldots , x^{\ell}) = \wt{\g} (dx^{j+1} \ldots dx^{\ell} \, x^0 \, 
dx^1 \ldots dx^j) \ , \ x^j \in {\Hc}_* \, . \leqno (94)
$$
We can assume that $\g (1,g_1 , \ldots , g_n) = {\displaystyle 
\prod_{j=1}^{n}} 
P_j (\d (g_j \cdot k) \, \rho_i (k)$ where each $P_j$ is a polynomial (in fact 
monomial) in the $\d_n$.

\smallskip

\noindent We take the $\rho_j$, $j=0,1,2,3$ as a basis of $A^* (G_1)$ viewed 
as 
a module over $C_c^{\ify} (G_1)$ and for $f \in C_c^{\ify} (G_1)$ we write 
$df$ 
as
$$
df = -(Yf) \, \rho_1 + (Xf) \, \rho_2 \ , \ Y = -{\part \over \part s} \ , \ X 
= e^{-s} \, {\part \over \part x} \leqno (95)
$$
which is thus expressed in terms of the left action of ${\Hc}$ on ${\Hc}_* = 
C_c^{\ify} (G_1) \semi G_2$. Moreover, using (86), one has $U_{\psi}^* \, 
\rho_j \, U_{\psi} = \rho_j$, $j \ne 1$ and
$$
U_{\psi}^* \, \rho_1 \, U_{\psi} = \rho_1 - \d_1 (U_{\psi}^*) \, U_{\psi} \, 
\rho_2 \leqno (96)
$$
or in other terms $\rho_1 \, U_{\psi}^* = \d_1 (U_{\psi}^*) \rho_2 + 
U_{\psi}^* 
\, \rho_1$.

\smallskip

\noindent This shows that provided we replace some of the $x_j$'s in (94) by 
the $h_j (x_j)$, $h_j \in {\Hc}$, we can get rid of all the exterior 
differentials $df$ and move all the $\rho_j$'s to the end of the expression 
which becomes,
$$
\wt{\g} \, (f^0 \, \d_{\psi_0} \, U_{\psi_0}^* \, f^1 \, \d_{\psi_1} \, 
U_{\psi_1}^* \ldots f^{\ell} \, \d_{\psi_{\ell}} \, U_{\psi_{\ell}}^* \; 
\rho_i) \leqno (97)
$$
provided we relabel the $x_j$'s in a cyclic way (which is allowed by 
definition 
2) and we omit several $\d_{\psi_j}$.

\smallskip

\noindent To write (97) in the form (6) we can assume that $\psi_{\ell} \ldots 
\psi_1 \, \psi_0 = 1$ since otherwise one gets 0. We first simplify the 
parenthesis using the crossed product rule in $C$ and get,
$$
\matrix{
f^0 \, (f^1 \circ \psi^0) (f^2 \circ \psi^1 \circ \psi^0) \ldots (f^{\ell} 
\circ \psi^{\ell -1} \ldots \psi_0) \cr
\d_{\psi_0} (\d_{\psi_1 \psi_0} - \d_{\psi_0}) \ldots (\d_{\psi_{\ell} \ldots 
\psi_0} - \d_{\psi_{\ell - 1} \ldots \psi_0}) \, \rho_i \, . \cr
} \leqno (98)
$$
Let $f = f^0 (f^1 \circ \psi^0) \ldots (f^{\ell} \circ \psi^{\ell -1} \ldots 
\psi_0)$. When we apply $\wt{\g}$ to (98) we get
$$
\int_{G_1} f(k) \prod_{j=1}^{\ell} P_j (\d (\psi^{j-1} \ldots \psi^0 \cdot k)) 
\, \rho (k) \leqno (99)
$$
where we used the equality $\d_g^2 = 0$ in $\L \, \Cb \, G'_2$. The same 
result holds if we omit several $\d_{\psi_j}$, one just takes $P_j = 1$ in the 
expression (99).  

\smallskip

\noindent We now rewrite (99) in the form,
$$
\tau_0 \, (f \, (P_0 (\d) \, U_{\psi_0}^*) \, U_{\psi_0} \, (P_1 (\d) \, 
U_{\psi_1 \psi_0}^*) \, U_{\psi_1 \psi_0} \ldots (P_{\ell} (\d) \, 
U_{\psi_{\ell -1} \ldots \psi_0}^*) \, U_{\psi_{\ell -1} \ldots \psi_0}) \, . 
\leqno (100)
$$
Let us replace $f$ by $f^0 (f^1 \circ \psi^0) \ldots (f^{\ell} \circ 
\psi^{\ell 
-1} \ldots \psi_0)$ and move the $f^j$ so that they appear without 
composition, 
we thus get
$$
\tau_0 \, (f^0 (P_0 (\d) U_{\psi_0}^*) \, f^1 (U_{\psi_0} \, P_1 (\d)  
U_{\psi_1 \psi_0}^*) \, f^2 (U_{\psi_1 \psi_0} \, P_2 (\d) \, U_{\psi_2 \psi_1 
\psi_0}^*) \ldots ) \, . \leqno (101)
$$
We now use the coproduct rule to rearange the terms, thus
$$
P_1 (\d) \, U_{\psi_1 \psi_0}^* = \sum \ P_1^{(1)} (\d) \, U_{\psi_0}^* \, \ 
P_1^{(2)} (\d) \, U_{\psi_1}^* \leqno (102)
$$
and we can permute $f^1$ with $U_{\psi_0} \, P_1^{(1)} \, (\d) \, 
U_{\psi_0}^*$ 
and use the equality
$$
(P(\d) \, U_{\psi_0}^*) (U_{\psi_0} \, Q \, (\d) \, U_{\psi_0}^*) = (PQ) (\d) 
\, U_{\psi_0}^* \, . \leqno (103)
$$
Proceeding like this we can rewrite (101) in the form,
$$
\tau_0 \, (f^0 \, Q_0 (\d) \, U_{\psi_0}^* \, f^1 \, Q_1 (\d) \, U_{\psi_1}^* 
\ldots f^{\ell} \, Q_{\ell} (\d) \, U_{\psi_{\ell}}^*) \leqno (104)
$$
which shows that the functional (94) satisfies definition 2.

\smallskip

\noindent Thus the map $\Phi$ of [Co] p.~220 together with Lemma 7 gives us a 
morphism $\t$ of complexes from the complex $C^* ({\Ac})$ of the Lie algebra 
cohomology of the Lie algebra $\Ac$ of formal vector fields, to the $(b,B)$ 
bicomplex of the Hopf algebra ${\Hc}$.

\smallskip

\noindent Since the current $c(g^0 , \ldots , g^n)$ is determined by its value 
at $1 \in G_1$, i.e. by
$$
\sum \, c_j (1) \, \rho_j \leqno (105)
$$
we can view the map $C$ as a map from $C^* ({\Ac})$ to the cochains of the 
group cohomology of $G_2$ with coefficients in the module $E$,
$$
E = \L \, {\bf G}_1^* \leqno (106)
$$
which is the exterior algebra on the cotangent space $T_1^* (G_1)$. Since the 
action of $G_2$ on $G_1$ fixes 1 it acts on $T_1^* (G_1)$ and in the basis 
$\rho_i$ the action is given by (86).

\smallskip

\noindent Since $\Ac$ is the direct sum ${\Ac} = {\bf G}_1 \op {\bf G}_2$ of 
the Lie algebras of $G_1$ and $G_2$ viewed as Lie subalgebras of ${\Ac}$ (it 
is 
a direct sum as vector spaces, not as Lie algebras), one has a natural 
isomorphism
$$
\L \, {\Ac}^* \sm \L \, {\bf G}_1^* \ot \L \, {\bf G}_2^* \leqno (107)
$$
of the cochains in $C^* ({\Ac})$ with cochains in $C^* ({\bf G}_2 , E)$, the 
Lie algebra cohomology of ${\bf G}_2$ with coefficients in $E$.

\medskip

\noindent {\bf Lemma 8.} {\it Under the above identifications, 
the map $C$ coincides with the cochain implementation 
$C^* ({\bf G}_2 , E) \ra C^* (G_2 , E)$ of the van Est isomorphism,  
which associates to a 
right invariant form $\mu$ on $G_2$, with values in $E$, the totally 
antisymmetric homogeneous cochain
$$
C(\mu) (g^0 , \ldots ,g^n) = \int_{\D (g^0 , \ldots ,g^n)} \mu \, .
$$
}
\smallskip

\noindent {\it Proof.} By (84) we know that there exists a right invariant 
form, $\mu = \sum \, \mu_j \, \rho_j$ on $G_2$ with values in $E$, such that
$$
c_j (1) = \int_{\D (g^0 , \ldots ,g^n)} \mu_j \, . \leqno (108)
$$
The value of $\mu_j$ at $1 \in G_2$ is obtained by contraction of $\om$ 
evaluated at $1 \in G$, by a suitable element of $\L \, {\bf G}_1$. Indeed 
this follows from (80) and the vanishing of the vector field $Z'$ of (82) at 
$a=1 \in G_2$. Thus the map $\om \ra \mu$ is the isomorphism (107).~\xx

\medskip

\noindent Of course the coboundary $d_1$ in the cochain complex $C^* ({\bf 
G}_2 , E)$ is not equal to the coboundary $d$ of $C^* ({\Ac})$, but it 
corresponds by the map $C$ to the coboundary $d_1$ of the bicomplex $(C^{n,m} 
, d_1 , d_2)$. We should thus check directly that $d_1$ and $d$ anticommute in 
$C^* ({\Ac})$. To see this, we introduce a bigrading in $C^* ({\Ac})$ 
associated to the decomposition ${\Ac} = {\bf G}_1 \op {\bf G}_2$. What we 
need to check is that the Lie algebra cohomology coboundary $d$ transforms an 
element of bidegree $(n,m)$ into a sum of two elements of bidegree $(n+1,m)$ 
and $(n,m+1)$ respectively. It is enough to do that for 1-forms. Let $\om$ be 
of bidegree $(1,0)$, then
$$
d\om (X_1 , X_2) = - \om ([X_1 , X_2]) \qquad \fl \, X_1 , X_2 \in {\Ac} \, . 
\leqno (109)
$$
This vanishes if $X_1 , X_2 \in {\bf G}_2$ thus showing that $d\om$ has no 
component of bidegree $(0,2)$.

\smallskip

\noindent We can then decompose $d$ as $d = d_1 + d_2$ where $d_1$ is of 
bidegree $(1,0)$ and $d_2$ of bidegree $(0,1)$.

\smallskip

\noindent Let us check that $d_1$ is the same as the coboundary of Lie algebra 
cohomology of ${\bf G}_2$ with coefficients in $\L \, {\bf G}_1^*$. Let $\a 
\in \L^m \, {\bf G}_1^*$ and $\om \in \L^n \, {\bf G}_2^*$. The component of 
bidegree $(n+1,m)$ of $d (\a \wdg \om) = (d\a) \wdg \om + (-1)^m \, \a \wdg 
d\om$ is
$$
d_1 \, \a \wdg \om + (-1)^m \, \a \wdg d_1 \, \om \leqno (110)
$$
where $d_1 \, \om$ takes care of the second term in the formula for the 
coboundary in Lie algebra cohomology,
$$
\matrix{
\sum (-1)^{i+1} \, X_i \, \om (X_1 , \ldots , \build{X_i}_{}^{\!\!\vee} , 
\ldots , X_{n+1}) + \cr
\sum_{i<j} (-1)^{i+j} \, \om ([X_i , X_j] , X_1 , \ldots 
\build{X_i}_{}^{\!\!\vee} , \ldots , \build{X_j}_{}^{\!\!\vee} , \ldots , 
X_{n+1}) \, . \cr
} \leqno (111)
$$
Thus, it remains to check that $d_1 \, \a \wdg \om$ corresponds to the first 
term in (111) for the natural action of ${\bf G}_2$ on $E = \L^* \, {\bf 
G}_1$, which can be done directly for $\a$ a 1-form, since $d_1 \, \a$ is the 
transpose of the action of ${\bf G}_2$ on ${\bf G}_1$.

\smallskip

\noindent We can now state the main lemma allowing to prove the surjectivity 
of the map $\t$.

\medskip

\noindent {\bf Lemma 9.} {\it The map $\t$ from $(C^* ({\Ac}) , d_1)$ to the 
Hochschild complex of $\Hc$ gives an isomorphism in cohomology.}

\medskip

\noindent {\it Proof.} Let us first observe that by Lemma 8 and the above 
proof of the van Est theorem, the map $C$ gives an isomorphism in cohomology 
from $(C^* ({\Ac}) , d_1)$ to the complex $(C^* , d_1)$ of polynomial 
$G_2$-group cochains with coefficients in $E$. Thus it is enough to show that 
the map $\Phi$ gives an isomorphism at the level of Hochschild cohomology. 
This of course requires to understand the Hochschild cohomology of the algebra 
${\Hc}_*$ with coefficients in the module $\Cb$ given by the augmentation on 
${\Hc}_*$.

\smallskip

\noindent In order to do this we shall use the abstract version ([C-E] 
Theorem~6.1 p.~349) of the Hochschild-Serre spectral sequence [Ho-Se].

\smallskip

A subalgebra ${\Ac}_1 \sbs {\Ac}_0$ of a augmented algebra ${\Ac}_0$ is 
called {\it normal} iff the right ideal $J$ generated in ${\Ac}_0$ by the 
${\rm Ker} \, \ve$ (of the augmentation $\ve$ of ${\Ac}_1$) is also a left 
ideal.

\smallskip

\noindent When this is so, one lets
$$
{\Ac}_2 = {\Ac}_0 / J \leqno (112)
$$
be the quotient of ${\Ac}_0$ by the ideal $J$. One has then a spectral 
sequence converging to the Hochschild cohomology $H_{{\Ac}_0}^* (\Cb)$ of 
${\Ac}_0$ with coefficients in the module $\Cb$ (using the augmentation $e$) 
and with $E_2$ term given by
$$
H_{{\Ac}_2}^p (H_{{\Ac}_1}^q (\Cb)) \, . \leqno (113)
$$
To prove it one uses the equivalence, for any right ${\Ac}_2$-module $A$ and 
right ${\Ac}_0$-module $B$,
$$
{\rm Hom}_{{\Ac}_2} (A, {\rm Hom}_{{\Ac}_0} ({\Ac}_2 , B)) \sm {\rm 
Hom}_{{\Ac}_0} (A,B) \leqno (114)
$$
where in the left term one views ${\Ac}_2$ as an ${\Ac}_2 - {\Ac}_0$ bimodule 
using the quotient map $\vp : {\Ac}_0 \ra {\Ac}_2$ to get the right action of 
${\Ac}_0$. Also in the right term, one uses $\vp$ to turn $A$ into a right 
${\Ac}_0$-module.

\smallskip

\noindent Replacing $A$ by a projective resolution of ${\Ac}_2$-right modules 
and $B$ by an injective resolution of ${\Ac}_0$-right modules, and using (since 
$ {\Ac}_2 = \Cb \ot_{{\Ac}_1} {\Ac}_0)$ ) the 
equivalence,
$$
{\rm Hom}_{{\Ac}_0} ({\Ac}_2 , B) = {\rm Hom}_{{\Ac}_1} (\Cb , B) \ ,  
 \leqno (115)
$$
one obtains the desired spectral sequence.

\smallskip

\noindent In our case we let ${\Ac}_0 = {\Hc}_*$ and ${\Ac}_1 = C_c^{\ify} 
(G_1)$. These algebras are non unital but the results still apply.
 We first need to prove that ${\Ac}_1$ is a {\it normal} subalgebra of 
${\Ac}_0$. The augmentation $\ve$ on ${\Hc}_*$ is given by,
$$
\ve \, (f \, U_{\psi}^*) = f(1) \qquad \fl \, f \in C_c^{\ify} (G_1) \, , \ 
\psi \in G_2 \, . \leqno (116)
$$
Its restriction to ${\Ac}_1$ is thus $f \ra f(1)$. Thus the ideal $J$ is 
linearly generated by elements $g \, U_{\psi}^*$ where
$$
g \in C_c^{\ify} (G_1) \ , \ g(1) = 0 \, . \leqno (117)
$$
We need to show that it is a left ideal in ${\Ac}_0 = {\Hc}_*$. For this it is 
enough to show that $U_{\psi_1}^* \, g \, U_{\psi_2}^*$ is of the same form, 
but this follows because,
$$
\psi (1) = 1 \qquad \fl \, \psi \in G_2 \, . \leqno (118)
$$
Moreover, with the above 
notations, the algebra ${\Ac}_2$ is the group ring of $G_2$.
We thus obtain by [C-E] loc.~cit., a spectral sequence which converges to the 
Hochschild cohomology of $\Hc$ and whose $E_2$ term is given by the polynomial 
group cohomology of $G_2$ with coefficients in the Hochschild cohomology of 
the coalgebra ${\Uc} ({\bf G}_1)$, which according to proposition 6.1 is given 
by $\L \, {\bf G}_1^*$. \smallskip

\noindent It thus follows that, combining the van Est theorem with the above 
spectral sequence, the map $\t$ gives an isomorphism in Hochschild 
cohomology.~\xx

\medskip

\noindent We can summarize this section by the following result.

\medskip

\noindent {\bf Theorem 10.} {\it The map $\t$ defines an isomorphism from the 
Lie algebra cohomology of $\Ac$ to the periodic cyclic cohomology of $\Hc$.}

\medskip

\noindent This theorem extends to the higher dimensional case, i.e. where the 
Lie algebra $\Ac$ is replaced by the Lie algebra of formal vector fields in 
$n$-dimensions, while $G = {\rm Diff} \, \Rb^n$, with $G_1$ the subgroup of 
affine diffeomorphisms, and $G_2 = \{ \psi ; \psi (0) = 0 , \psi' (0) = {\rm 
id} 
\}$. It also admits a relative version in which one considers the Lie algebra 
cohomology of $\Ac$ relative to $SO (n) \sbs G_0 = GL (n)  \sbs G_1$.

\vfill\eject

\noindent {\bf VIII.  Characteristic classes for actions of 
Hopf algebras }

\smallskip

In the above section we have defined and computed the cyclic cohomology of the 
Hopf algebra ${\Hc}$ as the Lie algebra cohomology of $\Ac$, the Lie algebra 
of formal vector fields.

\smallskip

\noindent The theory of characteristic classes for actions of $\Hc$ extends 
the construction (cf.~[Co2]) of cyclic cocycles from a Lie algebra of 
derivations of a $C^*$ algebra $A$, together with an {\it invariant trace} 
$\tau$ on $A$. At the purely algebraic level, given an algebra $A$ and an 
action of the Hopf algebra $\Hc$ on $A$,
$$
{\Hc} \ot A \ra A \ , \ h \ot a \ra h(a) \leqno (1)
$$
satisfying $h_1 (h_2 \, a) = (h_1 \, h_2) (a) \quad \fl \, h_i \in {\Hc}$ and
$$
h(ab) = \sum \, h_{(0)} \, (a) \, h_{(1)} \, (b) \qquad \fl a,b \in A
$$
we shall say that a trace $\tau$ on $A$ is {\it invariant} iff the following 
holds,
$$
\tau (h(a)b) = \tau (a \, \wt S (h)(b)) \qquad \fl \, a,b  \in A \, , \ h \in 
{\Hc} \, . \leqno (2)
$$
One has the following straightforward,

\medskip

\noindent {\bf Proposition 1.} {\it Let $\tau$ be an ${\Hc}$-invariant trace 
on $A$, then the following defines a canonical map $\gamma:HC^* ({\Hc}) \ra 
HC^* (A)$,
$$
\matrix{
\g (h^1 \ot \ldots \ot h^n) \in C^n (A) \, , \ \g (h^1 \ot \ldots \ot h^n) 
(x^0 , \ldots , x^n) = \cr
\cr
\tau (x^0 \, h^1 (x^1) \ldots h^n (x^n)) \, . \cr
}
$$
}

\medskip

\noindent In the interesting examples the algebra $A$ is a $C^*$ algebra and 
the action of $\Hc$ on $A$ is only densely defined. It is then crucial to know 
that the common domain,
$$
A^{\ify} = \{ a \in A \, ; \ h(a) \in A ,\qquad \fl \, h \in {\Hc} \} \leqno 
(3)
$$
is a subalgebra stable under the holomorphic functional calculus.

\smallskip

\noindent It is clear from the coproduct rule that $A^{\ify}$ is a subalgebra 
of $A$, the question is to show the stability under holomorphic functional 
calculus.

\smallskip

\noindent Our aim is to show that for any action of our Hopf algebra $\Hc$ (of 
section 2) on a $C^*$ algebra $A$ the common domain (3) is stable under 
holomorphic functional calculus. For short we shall say that a Hopf algebra 
$\Hc$ is {\it differential} iff this holds.

\medskip

\noindent {\bf Lemma 2.} {\it Let $\Hc$ be a Hopf algebra, ${\Hc}_1 \sbs \Hc$ 
a Hopf subalgebra. We assume that ${\Hc}_1$ is {\it differential} and that as 
an algebra $\Hc$ is generated by ${\Hc}_1$ and an element $\d \in {\Hc}$, $(\ve 
(\d) = 0)$ such 
that the following holds,} 

\noindent a) {\it For any $h \in {\Hc}_1$, there exist $h_1 , h_2 , h'_1 , h'_2 
\in {\Hc}_1$ 
such that
$$
\d h = h_1 \d + h_2 \ , \ h\d = \d h'_1 + h'_2 \, .
$$
}

\noindent b) {\it There exists $R \in {\Hc}_1 \ot {\Hc}_1$ such that $\D \, \d 
= \d \ot 1 + 1 \ot \d + R$.}
\medskip
 
\noindent {\it Then $\Hc$ is a differential Hopf algebra.}

\medskip

\noindent {\it Proof.} Let $A^{1,\ify} = \{ a \in A \, ; \ h(a) \in A \quad 
\fl \, h \in {\Hc}_1 \}$. By hypothesis if $u$ is invertible in $A$ and $u \in 
A^{1,\ify}$ one has $u^{-1} \in A^{1,\ify}$. It is enough, modulo some simple 
properties of the resolvent ([Co]) to show that the same holds in $A^{\ify}$. 
One has using a) that $A^{\ify} = \, \build \cap_{k}^{} \, \{ a \in A^{1,\ify} 
\, ; \ \d^j (a) \in A^{1,\ify} \quad 1 \leq j \leq k \} = \, \build 
\cap_{k}^{} \, A^{\ify ,k}$. Let $u \in A^{\ify}$ be invertible in $A$, then 
$u^{-1} \in A^{1,\ify}$. Let us show by induction on $k$ that $u^{-1} \in 
A^{\ify ,k}$. With $R = \sum \, h_i \ot k_i$ it follows from b) that
$$
\d (u \, u^{-1}) = \d (u) \, u^{-1} + u \, \d (u^{-1}) + \sum \, h_i (u) \, 
k_i (u^{-1}) \leqno (4)
$$
so that,
$$
\d (u^{-1}) = u^{-1} \, \d (1) - u^{-1} \, \d (u) \, u^{-1} - \sum \, u^{-1} 
h_i 
(u) 
\, k_i (u^{-1}) \, . \leqno (5)
$$
Since $\ve \, (\d) = 0$ one has $\d (1) = 0$ and the r.h.s. of (5) belongs to 
$A^{1,\ify}$. This shows that $u^{-1} \in A^{\ify ,1}$. Now by a), $A^{\ify 
,1}$ is stable by the action of ${\Hc}_1$ and by b) it is a subalgebra of 
$A^{1,\ify}$. Thus (5) shows that $\d (u^{-1}) \in A^{\ify ,1}$, i.e. that 
$u^{-1} \in A^{\ify ,2}$. Now similarly since,
$$
A^{\ify ,2} = \{ a \in A^{\ify ,1} \, ; \ \d (a) \in A^{\ify ,1} \} \leqno (6)
$$
we see that $A^{\ify ,2}$ is a subalgebra of $A^{\ify ,1}$ stable under the 
action of ${\Hc}_1$, so that by (5) we get $\d (u^{-1}) \in A^{\ify ,2}$ and 
$u^{-1}\in A^{\ify ,3}$. The conclusion follows by induction.~\xx

\medskip

\noindent {\bf Proposition 3.} {\it The Hopf algebra $\Hc$ of Section 2 is 
differential.}

\medskip

\noindent {\it Proof.} Let ${\Hc}^1 \sbs {\Hc}$ be the inductive limit of the 
Hopf subalgebras ${\Hc}_n$. For each $n$ the inclusion of ${\Hc}_n$ in 
${\Hc}_{n+1}$ fulfills the hypothesis of Lemma 2 so that ${\Hc}^1 = 
\cup_n \, {\Hc}_n$ is differential.

\smallskip

\noindent Let then ${\Hc}^2 \sbs {\Hc}$ be generated by ${\Hc}^1$ and $Y$, 
again the inclusion ${\Hc}^1 \sbs {\Hc}^2$ fulfills the hypothesis of the 
Lemma 2 since $\D Y = Y \ot 1 + 1 \ot Y$ while $h \ra Yh - hY$ is a derivation 
of ${\Hc}^1$. Finally $\Hc$ is generated over ${\Hc}^2$ by $X$ and one checks 
again the hypothesis of Lemma 2, since in particular $XY = (Y+1) \, X$.~\xx

\medskip

\noindent It is not difficult to give examples of Hopf algebras which are not 
differential such as the Hopf algebra generated by $\t$ with
$$
\D \, \t = \t \ot \t \, . \leqno (7)
$$

\vglue 1cm
\noindent {\bf IX. The index formula}
\smallskip
We use the notations of section II, so that ${\Ac} = C_c^{\ify} (F^+) \semi \G$ 
is the crossed product of the positive frame bundle of a flat manifold $M$ by a 
pseudogroup $\G$. We let $v$ be the canonical ${\rm Diff}^+$ invariant volume 
form on $F^+$ and,
$$
L^2 = L^2 (F^+ ,v) \, , \leqno (1)
$$
be the corresponding Hilbert space.
\smallskip
\noindent The canonical representation of ${\Ac}$ in $L^2$ is given by
$$
(\pi (f \, U_{\psi}^*) \, \xi) (j) = f(j) \, \xi (\wt{\psi} (j)) \quad \fl \, j 
\in F^+ \, , \ \xi \in L^2 \, , \ f \, U_{\psi}^* \in {\Ac} \, . \leqno (2)
$$
The invariance of $v$ shows that $\pi$ is a unitary representation for the 
natural involution of $\Ac$. When no confusion can arise we shall write simply  
$ a \xi$ instead of $\pi (a) \, \xi$ for $a \in {\Ac}$, $\xi \in L^2$.
\smallskip
\noindent The flat connection on $M$ allows to extend the canonical action of 
the group $G_0 = GL^+ (n,\Rb)$ on $F^+$ to an action of the affine group,
$$
G_1 = \Rb^n \semi G_0 \leqno (3)
$$
generated by the vector fields ${X}_i$ and ${Y}_{\ell}^k$ of section I.
\smallskip
\noindent This representation of $G_1$ admits the following compatibility with 
the representation $\pi$ of $\Ac$,
\medskip
\noindent {\bf Lemma 1.} {\it Let $a \in {\Ac}$ one has,}
\noindent 1) ${Y}_{\ell}^k \, \pi (a) = \pi (a) \, {Y}_{\ell}^k + \pi 
(Y_{\ell}^k (a))$
\noindent 2) ${ X}_i \, \pi (a) = \pi (a) \, {X}_i + \pi (X_i (a)) + \pi 
(\d_{ij}^k (a)) \, { Y}_k^j$.
\medskip
\noindent {\it Proof.} With $a = f \, U_{\psi}^* \in {\Ac}$ and $\xi \in 
C_c^{\ify} (F^+)$ one has ${ Y}_{\ell}^k \, \pi (a) \, \xi = {Y}_{\ell}^k 
(f \, \xi \circ \psi) = Y_{\ell}^k (f) \, \xi \circ \psi + f \, { Y}_{\ell}^k 
(\xi \circ \psi)$. Since $Y_{\ell}^k$ commutes with diffeomorphisms the last 
term is $({Y}_{\ell}^k \, \xi) \circ \psi$ which gives 1).
\smallskip
\noindent The operators $X_i$ acting on $\Ac$ satisfy (11) of section II which 
one can specialize to $b \in C_c^{\ify} (F^+)$, to get
$$
X_i (f \, U_{\psi}^* \, b) = X_i (f \, U_{\psi}^*) \, b + f \, U_{\psi}^* \, 
X_i 
(b) + \d_{ij}^k (f \, U_{\psi}^*) \, Y_k^j (b) \, . \leqno (4)
$$
One has $X_i (f \, U_{\psi}^* \, b) = X_i (f (b \circ \psi)) \, U_{\psi}^* = 
{ X}_i (\pi (a) b) \, U_{\psi}^*$, $X_i (f \, U_{\psi}^*) \, b = \pi (X_i 
(a))$ $b \, U_{\psi}^*$, $f \, U_{\psi}^* \, X_i (b) = \pi (a) \, {X}_i (b) 
\, U_{\psi}^*$ and $\d_{ij}^k (f \, U_{\psi}^*) \, Y_k^j (b) = \pi (\d_{ij}^k 
(a)) \, {Y}_k^j (b) \, U_{\psi}^*$, thus one gets 2).~\xx
\medskip
\noindent {\bf Corollary 2.} {\it For any element $Q$ of ${\Uc} ({\bf G}_1)$ 
there exists finitely many elements 
$Q_i \in {\Uc} ({\bf G}_1)$ and $h_i \in {\Hc}$ 
such that
$$
Q \, \pi (a) = \sum \, \pi (h_i (a)) \, Q_i \qquad \fl \, a \in {\Ac} \, .
$$
}

\smallskip
\noindent {\it Proof.} This condition defines a subalgebra of ${\Uc} ({\bf 
G}_1)$ and we just checked it for the generators.~\xx
\medskip
\noindent Let us now be more specific and take for $Q$ the hypoelliptic 
signature operator on $F^+$. It is not a scalar operator but it acts in the 
tensor product
$$
{\Hc}_0 = L^2 (F^+ , v) \ot E \leqno (5)
$$
where $E$ is a finite dimensional representation of $SO(n)$ specifically given 
by
$$
E = \L \, {\bf P}_n \ot \L \, \Rb^n \ , \ {\bf P}_n = S^2 \, \Rb^n \, . \leqno (6)
$$
The operator $Q$ is the graded sum,
$$
Q = (d_V^* \, d_V - d_V \, d_V^*) \op (d_H + d_H^*) \leqno (7)
$$
where the horizontal (resp. vertical) differentiation $d_H$ (resp. $d_V$) is a 
matrix in the ${X}_i$ resp. ${Y}_{\ell}^k$.  When $n$ is equal to 1 or 2 modulo 4 one has to replace
$F^+$ by its product by $ S^1$ so that the dimension of the vertical fiber is even (it is then $1+{n(n+1) \over 2}$ )
and the vertical signature operator makes sense.  The longitudinal part is not elliptic but only 
transversally elliptic with respect to the action of $SO(n)$.  
Thus, to get a 
hypoelliptic operator, one restricts $Q$ to the Hilbert space
$$
{\Hc}_1 = (L^2 (F^+ , v) \ot E)^{SO(n)} \leqno (8)
$$
and one takes the following subalgebra of $\Ac$,
$$
{\Ac}_1 = ({\Ac})^{SO(n)} = C_c^{\ify} (P) \semi \G \ , \ P = F^+ / SO(n) \, . 
\leqno (9)
$$
Let us note that the operator $Q$ is in fact the image under the
right regular representation of the affine group $G_1$
of a (matrix-valued)
hypoelliptic symmetric element in ${\Uc}({\bf G}_1)$. 
By an easy adaptation
of a theorem of Nelson and Stinespring, it then follows that 
$Q$ is essentially selfadjoint
(with core any dense, $G_1$-invariant subspace of 
the space of $C^{\ify}$-vectors of the right regular
representation of $G_1$). 
\smallskip
\noindent We let $\tau$ be the trace on ${\Ac}_1$ which is dual to the invariant volume 
from $v_1$ on $P$, that is,
$$
\tau \, (f \, U_{\psi}^*) = 0 \quad \hbox{if} \quad \psi \ne 1 \, , \ \tau (f) 
= 
\int_P f \, dv_1 \, . \leqno (10)
$$
Also we adapt the result of the previous sections to the relative case, and use 
the action of ${\Hc}$ on ${\Ac}$ to get a characteristic map,
$$
HC^* ({\Hc} , SO(n)) \ra HC^* ({\Ac}_1) \leqno (11)
$$
associated to the trace $\tau$.
\medskip
\noindent {\bf Proposition 3.} {\it Let us assume that the action of $\G$ on 
$M$ 
has no degenerate fixed point. Then any cochain on ${\Ac}_1$ of the form,
$$
\vp (a^0 , \ldots , a^n) = \, \int \!\!\!\!\!\! - \, a^0 [Q,a^1]^{(k_1)} \ldots 
[Q , a^n]^{(k_n)} \, (Q^2)^{-{1 \over 2} (n + 2\vert k \vert)} \qquad \fl \, a^j 
\in {\Ac}_1
$$
(with $T^{(k)} = [Q^2 , \ldots [Q^2 , T] \ldots ]$), is in the range of the 
characteristic map.
}

\medskip
\noindent {\it Proof.} By Corollary 2 one can write each operator 
$a^0 [Q,a^1]^{(k_1)} \ldots [Q,a_n]^{(k_n)}$ in the form 
$$
\sum_{\a} \, a^0 \, h_1^{\a} (a^1) \ldots h_n^{\a} (a^n) \, Q_{\a} \qquad h_j 
\in {\Hc} \ , \ Q_{\a} \in {\Uc} (G_1) \, , \leqno (12)
$$
thus we just need to understand the cochains of the form,
$$
\int \!\!\!\!\!\! - \, a^0 \, h_1 (a^1) \ldots h_n (a^n) \, R = \,  \int 
\!\!\!\!\!\! - \, a \, R \leqno (13)
$$
where $R$ is pseudodifferential in the hypoelliptic calculus and commutes with 
the affine subgroup of ${\rm Diff}$,
$$
R \, U_{\psi} = U_{\psi} \, R \qquad \fl \, \psi \in G'_1 \sbs {\rm Diff} \, . 
\leqno (14)
$$
Since $R$ is given by a smoothing kernel outside the diagonal and the action of 
$\G$ on $F^+ = F^+ (M)$ is free by hypothesis, one gets that
$$
\int \!\!\!\!\!\! - \, f \, U_{\psi}^* \, R = 0 \qquad \hbox{if} \ \psi \ne 1 
\, 
. \leqno (15)
$$
Also by (14) the functional ${\displaystyle \int \!\!\!\!\!\! -} \, f \, R$ is 
$G'_1$ invariant and is hence proportional to
$$
\int f \, dv_1 \, . \leqno (16)
$$
We have thus proved that $\vp$ can be written as a finite linear combination,
$$
\vp (a^0 , \ldots , a^n) = \sum_{\a} \, \tau (a^0 \, h_1^{\a} (a^1) \ldots 
h_n^{\a} (a^n)) \, . \leqno (17)
$$
Since we restrict to the subalgebra ${\Ac}_1 \sbs {\Ac}$, the cochain
$$
c = \sum \, h_1^{\a} \ot \ldots \ot h_n^{\a} \leqno (18)
$$
should be viewed as a {\it basic} cochain in the cyclic complex of ${\Hc}$, 
relative to the subalgebra ${\Hc}_0 = {\Uc} ({\bf O}(n))$.~\xx
\medskip
\noindent Now by theorem 10 of section VII one has an isomorphism,
$$
 H^* ({\Ac}_n , SO(n)) \build \ra_{}^{\t} HC^* ({\Hc} , SO(n)) \leqno (19)
$$
where the left hand side is the relative Lie algebra cohomology of the Lie 
algebra of formal vector fields.
\smallskip
\noindent Let us recall the result of Gelfand-Fuchs (cf.~[G-F], [G]) 
which allows to 
compute the left hand side of (19).
\smallskip
\noindent One lets $G_0 = GL^+ (n,\Rb)$ and ${\bf G}_0$ its Lie algebra, 
viewed 
as a subalgebra of ${\bf G}_1 \sbs {\Ac}_n$. One then views (cf.~[G])
the natural projection,
$$
\pi : {\Ac}_n \ra {\bf G}_0 \leqno (20)
$$
as a connection 1-form. Its restriction to ${\bf G}_0 \sbs {\Ac}_n$ is the 
identity map and,
$$
\pi \, ([X_0 , X]) = [X_0 , \pi (X)] \qquad \fl \, X_0 \in {\bf G}_0 \, , \ X 
\in {\Ac}_n \, . \leqno (21)
$$
The curvature of this connection is
$$
\Om = d\pi + {1 \over 2} \, [\pi , \pi] \leqno (22)
$$
is easy to compute and is given by 
$$
\Om (X,Y) = [Y_1 , X_{-1}] - [X_1 , Y_{-1}] \, , \qquad \fl \, X,Y \in {\Ac}_n 
\, , \leqno (23)
$$
in terms of the projections $X \ra (X)_j$ associated to the grading of the Lie 
algebra ${\Ac}_n$.
\smallskip
\noindent It follows from the Chern-Weil theory that one has a canonical map 
from the Weil complex, $S^* ({\bf G}_0) \wdg \L^* ({\bf G}_0) = W ({\bf G}_0)$,
$$
W ({\bf G}_0) \build \ra_{}^{\vp} C^* ({\Ac}_n) \, . \leqno (24)
$$

For $\xi \in {\bf G}_0^*$ viewed as an odd element $\xi^- \in \L^* \, {\bf 
G}_0$, one has
$$
\vp (\xi^-) \in {\Ac}_n^* \ , \quad \vp (\xi^-) = \xi \circ \pi \, . \leqno 
(25)
$$
For $\xi \in {\bf G}_0^*$ viewed as an even element $\xi^+ (\in S^* \, {\bf 
G}_0)$, one has
$$
\vp (\xi^+) \in \L^2 \, {\Ac}_n^* \ , \quad \vp (\xi^+) = \xi \circ \Om \, , 
\leqno (26)
$$
moreover the map $\vp$ is an algebra morphism, which fixes it uniquely. By 
construction $\vp$ vanishes on the ideal $J$ of $W = W ({\bf G}_0)$ generated 
by 
${\displaystyle \sum_{r>n}} \, S^r ({\bf G}_0^*)$. By definition one lets
$$
W_n = W / J \leqno (27)
$$
be the corresponding differential algebra and $\wt{\vp}$ the quotient map,
$$
\wt{\vp} : H^* (W_n) \ra H^* ({\Ac}_n) \, . \leqno (28)
$$
It is an isomorphism by [G-F].
\smallskip
\noindent The discussion extends to the relative situation and yields a 
subcomplex
$WSO(n)$ of the elements of $W_n$ which are basic relative to the action of 
$SO(n)$. Again by [G-F], the morphism $\wt{\vp}$ gives an isomorphism,
$$
\wt{\vp} : H^* (WSO(n)) \ra H^* ({\Ac}_n ,\ SO(n)) \, . \leqno (29)
$$
A concrete description of $H^* (WSO(n))$ is obtained (cf.~[G]) as a small 
variant of $H^* (WO(n))$, i.e. the orthogonal case. The latter is the 
cohomology 
of the complex
$$
E (h_1 , h_3 , \ldots , h_m) \ot P (c_1 , \ldots ,c_n) \ , \ d \leqno (30)
$$
where $E (h_1 , h_3 , \ldots , h_m)$ is the exterior algebra in the generators 
$h_i$ of dimension $2i-1$, (m is the largest odd integer less than n) and $i$ 
odd $\leq n$, while $P (c_1 , \ldots ,c_n)$ is 
the polynomial algebra in the generators $c_i$ of degree $2i$ {\it truncated} 
by 
the ideal of elements of weight $>2n$. The coboundary $d$ is defined by,
$$
dh_i = c_i \ , \ i \ \hbox{odd} \ , \ dc_i = 0  \ \hbox{for all} \ i \, . 
\leqno (31)
$$

\noindent One lets $p_i = c_{2i}$ be the Pontrjagin classes, they are non 
trivial cohomology classes for $2i 
\leq n$. One has $H^* (WSO(n)) = H^* (WO(n))$ for $n$ odd, while for $n$ even 
one has
$$
H^* (WSO (n)) = H^* (WO(n)) \, [\chi] / ( \chi^2 - c_n ) \, . \leqno (32)
$$
Let us now recall the index theorem of [C-M] for spectral triples $({\Ac} , 
{\Hc} , D)$ whose dimension spectrum is discrete and simple, which is the case 
(cf.~[C-M]) for the transverse fundamental class, (we treat the {\it odd} case)
\medskip
\noindent {\bf Theorem 4.} a) {\it The equality ${\displaystyle \int 
\!\!\!\!\!\! -} \, P = {\rm Res}_{z=0} \, \hbox{\rm Trace} \, (P \vert D 
\vert^{-z})$ defines a trace on the algebra generated by ${\Ac}$, $[D, {\Ac}]$ 
and $\vert D \vert^z$, $z \in \Cb$.}
\noindent b) {\it The following formula only has a finite number of non zero 
terms and defines the components $(\vp_n)_{n=1,3\ldots}$ of a cocycle in the 
$(b,B)$ bicomplex of ${\Ac}$,
$$
\vp_n (a^0 , \ldots , a^n) = \sum_k \, c_{n,k} \, \int \!\!\!\!\!\! - \, a^0 
[D,a^1]^{(k_1)} \ldots [D,a^n]^{(k_n)} \, \vert D \vert^{-n-2 \vert k \vert}
$$
$\fl \, a^j \in {\Ac}$ where one lets $T^{(k)} = \nb^k (T)$, $\nb (T) = D^2 T - 
TD^2$ and where $k$ is a multiindex, $c_{n,k} = (-1)^{\vert k \vert} \, 
\sqrt{2i} \, (k_1 ! \ldots k_n !)^{-1} \, (k_1 + 1)^{-1} \ldots (k_1 + k_2 + 
\ldots + 
k_n + n)^{-1} \, \G \left( \vert k \vert + {n \over 2} \right)$, $\vert k \vert 
= k_1 + \ldots + k_n$.}
\noindent c) {\it The pairing of the cyclic cohomology class $(\vp_n) \in HC^* 
({\Ac})$ with $K_1 ({\Ac})$ gives the Fredholm index of $D$ with coefficients 
in 
$K_1 ({\Ac})$.}
\medskip
\noindent Let us remark that the theorem is unchanged if we replace everywhere 
the operator $D$ by
$$
Q = D \, \vert D \vert \, . \leqno (33)
$$
This follows directly from the proof in~[C-M].
\smallskip
\noindent In our case the operator $Q$ is differential, given by (7), so that 
by 
proposition 3 we know that the components $\vp_n$ belong to the range of the 
characteristic map. Since the computation is local we thus get, collecting 
together the results of this paper,
\medskip
\noindent {\bf Theorem 5.} {\it There exists for each $n$ a universal 
polynomial 
${\wt L}_n \in H^* (WSO_n)$ such that,
$$
Ch_* ( Q ) = \underline{\t} \, ({\wt L}_n) .
$$
}

\noindent Here $\t$ is the isomorphism of theorem 10 section VII,
in its relative version $(19)$ and 
$\underline{\t}$ denotes the composition of $\t$ with the 
relative characteristic map $(11)$ associated to the 
action of the Hopf algebra ${\Hc}$ on ${\Ac}$.
\medskip
\noindent One can end the computation of ${\wt L}_n$ by evaluating the index on the 
range of the assembly map,
$$
\mu : K_{*,\tau} (P {\textstyle \semi_{\G}} \, E \, \G) \ra K ({\Ac}_1) \, , 
\leqno (34)
$$
provided one makes use of the conjectured (but so far only partially
verified, cf.~[He], [K-T]) injectivity of the natural map,
$$
H_d^* (\G_n , \Rb) \ra H^* (B\G_n , \Rb) \leqno (35)
$$
from the smooth cohomology of the Haefliger groupoid $\G_n$ to its real 
cohomology.
\smallskip
\noindent One then obtains that ${\wt L}_n$ is the product of the 
usual $L$-class by another universal expression 
in the Pontrjagin classes $p_i$, accounting for the 
cohomological analogue of the $K$-theory Thom isomorphism 
$$
\beta : K_{*} ( C_0(M) {\textstyle \semi} \, \G) \ra 
K_{*} ( C_0 (P) {\textstyle \semi} \, \G)   
\leqno (36)
$$
of [Co1, \S V]. This can be checked
directly in small dimension. It is noteworthy also that the first Pontrjagin 
class $p_1$ already appears with a non zero coefficient for $n=2$.

\vglue 1cm
\bigskip

\noindent {\bf Appendix: the one-dimensional case}

\bigskip
In the one dimensional case the operator $Q$ is readily reduced to the 
following operator on the product $\wt{F}$
of the frame bundle $F^+$ by an auxiliary $S^1$ 
whose corresponding periodic coordinate is called $\a$ (and whose role is as mentionned above 
 to make the vertical fiber even dimensional).
$$
Q=Q_V+Q_H . \leqno (1)
$$
We work with 2 copies of $L^2 (\wt{F}, e^s d\a ds dx)$ and the
following gives the vertical operator $Q_V$,
$$
Q_V = \left[ \matrix{-\part_{\a}^2 + \part_s (\part_s +1)
&-2\part_{\a} \part_s - \part_{\a} \cr -2 \part_{\a} \part_s -
\part_{\a} &\part_{\a}^2 -\part_s (\part_s +1) \cr} \right] \,
. \leqno (2)
$$

The horizontal operator $Q_H$ is given by
$$
Q_H = {1 \over i} \, e^{-s} \part_{x} \g_2 \leqno (3)
$$
where $\g_2 = \left[ \matrix{ 0 &-i \cr i &0 \cr}
\right]$ anticommutes with $Q_V$.

\smallskip

We use the following notation for $2 \ts 2$ matrices
$$
\g_1 = \left[ \matrix{ 0 &1 \cr 1 &0 \cr} \right] \ , \ \g_2 =
\left[ \matrix{ 0 &-i \cr i &0 \cr} \right] \ ,\ \g_3 =
\left[ \matrix{ 1 &0 \cr 0 &-1 \cr} \right] \, . \leqno (4)
$$
We can thus write the full operator $Q$ acting in 2 copies of
$L^2 (\wt{F}, e^s d\a ds dx)$ as
$$
Q=( -2 \part_{\a} \part_s - \part_{\a}) \g_1 + {1 \over i} \, e^{-s} \part_x  \g_2 + \left( -\part_{\a}^2 + \part_s (\part_s +1) \right)
\g_3 \, . \leqno (5)
$$
 
\noindent {\bf Theorem A.} {\it Up to a coboundary $ch_*{Q}$ is
equal to twice the transverse fundamental class $[\wt{F}]$.}

\bigskip

The factor 2 is easy to understand since it is the local index
of the signature operator along the fibers of $\wt{F}$.

\smallskip

The formulas of theorem 4 for the components
$\vp_1$ and $\vp_3$ of the character are,
$$
\leqalignno{
\vp_1 (a^0 ,a^1) = & \ \sqrt{2i} \G \left( {1\over 2} \right)
\int 
\!\!\!\!\!\! -(a^0 [Q, a^1] (Q^2)^{-1/2}) &(6) \cr
- & \ \sqrt{2i} {1 \over 2} \, \G \left( {3 \over 2} \right)
\int 
\!\!\!\!\!\! -(a^0 \nb [Q,a^1] (Q^2)^{-3/2}) \cr
+ & \ \sqrt{2i} {1 \over 2} \, {1 \over 3} \, \G \left( {5 \over
2} \right) \int 
\!\!\!\!\!\! -(a^0 \nb^2 [Q,a^1] (Q^2)^{-5/2}) \cr
- & \ \sqrt{2i} {1 \over 2\cdot 3 \cdot 4} \, \G \left( {7 \over
2} \right) \int 
\!\!\!\!\!\! -(a^0 \nb^3 [Q,a^1] (Q^2)^{-7/2}) \cr
}
$$
$$
\leqalignno{
\vp_3 (a^0 , a^1 ,a^2 ,a^3) = & \ \sqrt{2i} {1 \over 3i} \, \G
\left( {3 \over 2} \right) \int 
\!\!\!\!\!\! -(a^0 [Q,a^1] [Q,a^2] 
[Q,a^3] (Q^2)^{-3/2} ) \cr
- & \ \sqrt{2i} {1 \over 2 \cdot 3 \cdot 4} \, \G \left( {5
\over 2} \right) \int 
\!\!\!\!\!\! -(a^0 \nb [Q,a^1] [Q,a^2] [Q,a^3]
(Q^2)^{-5/2} ) &(7) \cr
- & \ \sqrt{2i} {1 \over 3 \cdot 4} \, \G \left( {5 \over 2}
\right) \int 
\!\!\!\!\!\! -(a^0 [Q,a^1] \nb ([Q,a^2] [Q,a^3] (Q^2)^{-5/2} )
\cr
- & \ \sqrt{2i} {1 \over 2 \cdot 4} \, \G \left( {5 \over 2}
\right) \int 
\!\!\!\!\!\! -(a^0 [Q,a^1] [Q,a^2] \nb [Q,a^3] (Q^2)^{-5/2} )
\, . \cr
}
$$
The computation gives the following result,
$$
\vp_1 (a^0 ,a^1) = 0 \qquad \fl \, a^0 , a^1 \in \Ac \leqno
(8)
$$
in fact, each of the 4 terms of (6) turns out to be 0.
$$
\vp_3 = 2
(\wt{[F]} + b\psi), \leqno (9)
$$
where $\wt{[F]}$ is the tranverse fundamental class ([Co]), i.e. the extension to the crossed product of 
the following invariant cyclic 3-cocycle on the algebra $C_c^{\ify} (\wt{F}) $,
$$
\mu (f^0 ,f^1 ,f^2 ,f^3) = \int_{\wt{F}} f^0 df^1 \wdg df^2 \wdg df^3
\, . \leqno (10)
$$

We shall now give the explicit form of both $b\psi$ and $\psi$
with $B\psi =0$.

\smallskip
We let $\tau$ be the trace on $\Ac$ given by the measure
$$
f \in C_c^{\ify} (\wt{F)} \ra \int_{\wt{F}} f e^s d\a ds dx \, . \leqno
(11)
$$
This measure is invariant under the action of $\Diff^+$ and
thus gives a dual trace on the crossed product.

\smallskip
The two derivations $\part_{\a}$ and $\part_s$ of $C_c^{\ify}
(\wt{F})$ are invariant under the action of $\Diff^+$ and we denote
by the same letter their canonical extension to $\Ac$,
$$
\part_{\a} (f U_{\psi}^*) = (\part_{\a} f) U_{\psi}^* \ , \ \part_s
(f U_{\psi}^*) = (\part_s f) U_{\psi}^* \, . \leqno (12)
$$

We let $\d_1$ be the derivation of $\Ac$ defined in section II.
\smallskip
By construction both $\d_1$ and $\tau$ are invariant under
$\part_{\a}$ but neither of them is invariant under $\part_s$.

But the following
derivation $\part_u : \Ac \ra \Ac^*$ commutes with both
$\part_{\a}$ and $\part_s$,
$$
\lgl \part_u (a), b \rgl = \tau ( \d_1 (a) b ) \, . \leqno (13)
$$
We then view $\part_{\a}$ , $\part_s$ and $\part_u$ as three commuting derivations, where
$\part_u$ cannot be iterated and 
use the notation $(us ,\a ,s)$ for the cochain
$$
a^0 , a^1 , a^2 ,a^3 \ra \lgl a^0 , (\part_u \part_s a^1)
(\part_{\a} a^2)( \part_s a^3) \rgl \, .
$$
\smallskip

The formula for $b\psi$ is then the following,
$$
\eqalign{
{1 \over 8} \, (- & \ (u,\a ,s) + (u,s,\a) + (\a ,s,u) - (s,\a
,u)) \cr
{1 \over 2} \, ( - & \ (us,\a ,s) + (\a ,us,s) + (\a ,s,us) \cr
+ & \ (us,s,\a) + (s,us,\a) - (s,\a ,us) \qquad + (\a , u\a
,\a)) \cr
- & \ (u,\a s,s) - (s,u\a ,s) - (s,\a s ,u) \cr
{1 \over 4} \, ( - & \ (u,\a ,ss) + (ss,u,\a) - (ss,\a ,u) \cr
+ & \ (u,ss,\a) + (\a ,u,ss) + (\a ,ss,u) \cr
- & \ (u,\a \a \a) + (\a \a ,u, \a) - (\a \a ,\a ,u) \cr
- & \ (u,\a \a ,\a) + (\a ,u,\a \a) - (\a , \a \a ,u)) \, . \cr
}
$$

This formula is canonical and a possible choice of $\psi$ is
given by,
$$
{1 \over 8} \, ((\a ,su) - (\a u,s) - (s,\a u) + (su,\a))
$$
$$
{1 \over 4} \, (-(\a \a u ,\a) - {1 \over 3} (\a \a \a ,u) +
\a u , \a \a ))
$$
$$
{1 \over 4} ((u\a s,s) - (\a ,ssu) + {1 \over 2} (\a ss,u) +
(\a s ,us) + {1 \over 2} (\a u ,ss) - {1 \over 2} (ssu,\a)) \,
.
$$
The natural domain of $\psi$ is the algebra $C^3$ of 3 times
differentiable elements of $\Ac$ where the derivation
$\part_u$ is only used once.

\vglue 1cm  

\centerline{\bf Bibliography}

\bigskip

\item{[B-S]} S. Baaj and G. Skandalis : Unitaires multiplicatifs et dualit\'e 
pour les produits crois\'es de $C^*$-alg\`ebres, {\it Ann. Sci. Ec. Norm. 
Sup.}, 4 s\'erie, t.~26, 1993, (425-488). 

\item{[C-E]} H. Cartan and S. Eilenberg : Homological algebra, Princeton 
University Press (1956).

\item{[C-M1]} A. Connes and H. Moscovici : Cyclic cohomology, the Novikov 
conjecture and hyperbolic groups, {\it Topology} {\bf 29} (1990), 345-388.

\item{[C-M2]} A. Connes and H. Moscovici : The local index formula in 
noncommutative geometry, GAFA {\bf 5} (1995), 174-243.

\item{[Co]} A. Connes : Noncommutative geometry, Academic Press (1994).

\item{[Co1]} A. Connes : Cyclic cohomology and the transverse fundamental class  of a foliation, {\it Geometric methods in operator algebras}, (Kyoto, 1983), 
52-144; {\it Pitman Res. Notes in Math.} {\bf 123}, Longman, Harlow (1986).

\item{[Co2]} A. Connes : $C^*$-alg\`ebres et g\'eom\'etrie diff\'erentielle,
{\it C.R. Acad. Sci. Paris}, Ser.~A-B {\bf 290} (1980), A599-A604.
  
\item{[Dx]} J. Dixmier : Existence de traces non normales, 
{\it C.R. Acad. Sci. Paris}, Ser.~A-B {\bf 262} (1966), A1107-A1108.

\item{[G-F]} I. M. Gelfand and D.B. Fuks : Cohomology of the Lie algebra 
of formal 
vector fields, {\it Izv. Akad. Nauk SSSR} {\bf 34} (1970), 322-337.

\item{ \ } Cohomology of Lie algebra of vector fields with nontrivial 
coefficients, {\it Funct. Anal.} {\bf 4} (1970), 10-45.

\item{ \ } Cohomology of Lie algebra of tangential vector fields, {\it Funct. 
Anal.} {\bf 4} (1970), 23-31.

\item{[G]} C. Godbillon, Cohomologies d'alg\`ebres de Lie
de champs de vecteurs formels, 
{\it S\'eminaire Bourbaki (1971/1972), Expos\'e No.~421},
Lecture Notes in Math., Vol. {\bf 383}, 69-87, Springer, Berlin 1974.

\item{[H]} A. Haefliger, Sur les classes caract\'eristiques des feuilletages, 
{\it S\'eminaire Bourbaki (1971/1972), Expos\'e No.~412},
Lecture Notes in Math., Vol. {\bf 317}, 239-260, Springer, Berlin 1973.

\item{[He]} J. L. Heitsch : Independent variation of secondary classes, {\it 
Ann. of Math.} {\bf 108} (1978), 421-460.

\item{[H-S]} M. Hilsum and G. Skandalis : Morphismes $K$-orient\'e
d'espaces de feuilles et fonctorialit\'e en th\'eorie de Kasparov, 
{\it  Ann. Sci. \'Ecole Norm. Sup. (4)} {\bf 20} (1987), 325-390.

\item{[Ho-Se]} G. Hochschild and J.-P. Serre : Cohomology of Lie algebras, 
{\it  Ann. of Math.} {\bf 57} (1953), 591-603.

\item{[K]} G.I. Kac : Extensions of Groups to Ring Groups, {\it Math. USSR 
Sbornik}, {\bf 5} (1968), 451-474.

\item{[K-T]} F. Kamber and Ph. Tondeur : On the linear independence of certain
cohomology classes of $B{\Gamma}_q$, {\it Studies in algebraic topology}, 
Adv. in Math. Suppl. Studies, {\bf 5}, 213-263,
Academic Press, New York - London (1979).

\item{[M]} S. Majid : Foundations of Quantum Group Theory, Cambridge University Press (1995).

\item{[Ma]} Y. Manin, Quantum groups and noncommutative geometry, Centre 
Re\-cherche Math. Univ. Montr\'eal (1988).

\item{[M-S]} C. C. Moore and C. Schochet, Global analysis on foliated spaces,
{\it Math. Sci. Res. Inst. Publ.} {\bf 9}, Springer-Verlag New York Inc. 1988.

\bye